\documentclass[12pt,a4paper,dvipsnames]{amsart}

\usepackage[english]{babel}
\usepackage[T1]{fontenc}

\usepackage{etoolbox}
\apptocmd{\sloppy}{\hbadness 10000\relax}{}{}
\apptocmd{\sloppy}{\vbadness 10000\relax}{}{}

\usepackage{mathtools}
\usepackage{scalerel}

\usepackage{stackengine}

\usepackage[pagewise]{lineno}
\newcommand*\linenomathpatch[1]{%
  \cspreto{#1}{\linenomath}%
  \cspreto{#1*}{\linenomath}%
  \csappto{end#1}{\endlinenomath}%
  \csappto{end#1*}{\endlinenomath}%
}

\linenomathpatch{equation}
\linenomathpatch{gather}
\linenomathpatch{multline}
\linenomathpatch{align}
\linenomathpatch{alignat}
\linenomathpatch{flalign}

\usepackage[a4paper,top=2cm,bottom=2cm,left=1.5cm,right=1.5cm,marginparwidth=1.75cm]{geometry}

\usepackage{amsmath}
\usepackage{graphicx}
\usepackage[colorinlistoftodos]{todonotes}
\usepackage[pagebackref,hypertexnames=false, colorlinks, citecolor=red,linkcolor=blue, urlcolor=red]{hyperref}
\providecommand\phantomsection{}
\usepackage{diagbox,multirow} 
\usepackage{relsize}
\usepackage{comment}  

\usepackage{lmodern}
\usepackage[english]{babel}
\usepackage{amsmath}
\usepackage{amsfonts}
\usepackage{colonequals} 
\usepackage{amssymb}
\usepackage{amsthm}
\usepackage{stmaryrd}
\usepackage[all]{xy}

\usepackage{enumitem}
\setlist[enumerate]{topsep=1mm,parsep=0pt,partopsep=0mm,itemsep=1mm,leftmargin=12mm,labelsep=1mm} 

\usepackage{hyphenat}

\colorlet{Malte}{OliveGreen}
\colorlet{Takahiro}{red}
\colorlet{Michael}{blue}
\newcommand{\revise}[1]{{
    #1}}

\usepackage{pstricks}
\usepackage{epsfig}

\providecommand{\Asterisk}{\mathbin{\scalerel*{\ast}{\otimes}}}

\newtheorem{thm}{Theorem}[section] 
\newtheorem{lm}[thm]{Lemma}
\newtheorem{prop}[thm]{Proposition}

\newtheorem{cor}[thm]{Corollary}

\theoremstyle{definition} 
\newtheorem{defi}[thm]{Definition}
\theoremstyle{remark}
\newtheorem{rmq}[thm]{Remark}
\newtheorem{ex}[thm]{Example}

\numberwithin{equation}{section} 

\newcommand{\id}{\mathrm{id}}
\newcommand{\R}{\mathbb{R}}
\newcommand{\N}{\mathbb{N}}
\newcommand{\Z}{\mathbb{Z}}
\newcommand{\Alg}[1]{\operatorname{\ast-Alg}(#1)}
\newcommand{\RepCat}{\textnormal{$*$-}\mathfrak{Rep}}
\newcommand{\PreHilb}{\mathfrak{PreHilb}}
\newcommand{\Span}{\operatorname{span}}
\newcommand{\Lop}{L_{\mathrm a}}

\newcommand{\mon}{\mathbin{\scalerel*{\rhd}{\otimes}}}
\newcommand{\antimon}{\mathbin{\scalerel*{\lhd}{\otimes}}}
\newcommand{\bool}{\mathbin{
    \scalerel*{\diamond}{\otimes}}   
  }

\newcommand{\UU}{W} 
\newcommand{\TT}{T} 
\renewcommand{\SS}{S} 
\newcommand\AST{\ast} 
 
\newcommand{\Restr}[2]{{#1}{\restriction}_{#2}} 

\newcommand{\perm}[2]{\mathcal P^{#1}_{#2}}
\newcommand{\coef}[3]{c^{#1}_{#2}(#3)} 

\newcommand{\utprod}[2]{\mathop{{}_{#1}\otimes_{#2} }} 

\newcommand{\ufprod}[2]{\mathop{{}_{#1}{\overrightarrow\Asterisk}_{#2} }}

\newcommand{\twofacedprod}[3][.7]{
  \setstackgap{L}{#1\baselineskip}\mathbin{\scriptsize{{
    {
      \Vectorstack{{#2} {#3}}}}}}
\setstackgap{L}{\baselineskip}}

\newcommand{\uttprod}[4]{
  \twofacedprod{ {}_{#1}\otimes_{#2}} {{}_{#3}\otimes_{#4}}
}

\newcommand{\uffprod}[4]{\twofacedprod[.9]{
    \ufprod{#1}{#2}}{\ufprod{#3}{#4}}
}

\newcommand{\ubprod}[2]{\mathop{{}_{#1}{\overleftarrow\Asterisk}_{#2} }} 

\newcommand{\ubfprod}[4]{\twofacedprod[.9]{
    \ufprod{#1}{#2}}{\ubprod{#3}{#4}}
 } 

 \newcommand{\ubffprod}[4]{\twofacedprod[.9]{
    \ubprod{#1}{#2}}{\ufprod{#3}{#4}}} 

\newcommand{\ubfbfprod}[4]{\twofacedprod[.9]{
    \ubprod{#1}{#2}}{\ubprod{#3}{#4}}}

 \newcommand{\deffreeprod}[2]{ {\twofacedprod[.9]{\ufprod{#1}{#2}}{\ufprod{}{}}}}

\newcommand{\defbifreeprod}[2]{
     \twofacedprod[.9]{\ufprod{#1}{#2}}{\ubprod{}{}}
   }

\newcommand{\freeboolprod}{\twofacedprod{\Asterisk}{\bool}}

\newcommand{\boolfreeprod}{\twofacedprod{\bool}{\Asterisk}}

\newcommand{\boolboolprod}{\twofacedprod{\bool}{\bool}}

\newcommand{\fllift}%
{{\protect\overrightarrow{\lambda}}}
\newcommand{\frlift}%
{{\protect\overrightarrow{\rho}}}
\newcommand{\bfllift}%
{{\protect\overleftarrow{\lambda}}}
\newcommand{\bfrlift}%
{{\protect\overleftarrow{\rho}}}

\newcommand{\allift}{\overset{\antimon}{\lambda}{}} 
\newcommand{\arlift}{\overset{\antimon}{\rho}{}} 

\newcommand{\longrarrow}[1]{\ar@<0.5ex>[rr]^-{#1} }
\newcommand{\longlarrow}[1]{\ar@<0.5ex>[ll]^-{#1}}
\newcommand{\longdarrow}[1]{\ar@<0.5ex>[dd]^-{#1}}  
\newcommand{\longuarrow}[1]{\ar@<0.5ex>[uu]^-{#1}}  

\newcommand{\rarrow}[1]{\ar@<0.5ex>[r]^-{#1} }
\newcommand{\larrow}[1]{\ar@<0.5ex>[l]^-{#1}}
\newcommand{\darrow}[1]{\ar@<0.5ex>[d]^-{#1}}
\newcommand{\uarrow}[1]{\ar@<0.5ex>[u]^-{#1}}

\newcommand{\ruarrow}[1]{\ar@<0.5ex>[ru]^-{#1}}
\newcommand{\rdarrow}[1]{\ar@<0.5ex>[dr]^-{#1}}
\newcommand{\luarrow}[1]{\ar@<0.5ex>[lu]^-{#1}}
\newcommand{\ldarrow}[1]{\ar@<0.5ex>[ld]^-{#1}}

\newcommand{\rddarrow}[1]{\ar@<0.5ex>[rdd]^(.7){#1}}
\newcommand{\ruuarrow}[1]{\ar@<0.5ex>[ruu]^(.7){#1}}
\newcommand{\lddarrow}[1]{\ar@<0.5ex>[ldd]^(.7){#1}}
\newcommand{\luuarrow}[1]{\ar@<0.5ex>[luu]^(.7){#1}}

\newcommand{\subcase}[1]{\phantomsection\par\vspace{2mm}\noindent\textbf{[#1]}}
\newcommand{\subcasenpds}[1]{\phantomsection\par\vspace{4mm}\noindent\textbf{#1}.}

\newcommand{\subcasens}[1]{\phantomsection\par\noindent\textbf{[#1]}} 

\newcommand{\subcasenp}[1]{\phantomsection\par\vspace{2mm}\noindent\textbf{#1:}} 

\newcommand{\subcasenps}[1]{\phantomsection\par\noindent\textbf{#1:}} 

\title[Towards a classification of multi-faced independence]
{Towards a classification of multi-faced independence:\\ 
A representation-theoretic approach}

\author{Malte Gerhold}
\address[Malte Gerhold]{Department of Mathematical Sciences, NTNU Trondheim, Norway}
\address{Institute of Mathematics and Computer Science, University of Greifswald, Germany}
\thanks{M.G.\ was supported by the German Research Foundation (DFG) grant no.\ 397960675. The work of M.G.\ was partially carried out during the tenure of an ERCIM `Alain Bensoussan' Fellowship Programme and partially as a guest researcher at Saarland University in the scope of the SFB-TRR 195.}

\author{Takahiro Hasebe}

\address[Takahiro Hasebe]{Department of Mathematics, Hokkaido University Kita 10, Nishi 8, Kita-Ku, Sapporo, Hokkaido, Japan}
\thanks{T.H.\ is supported by JSPS Grant-in-Aid for Young Scientists 19K14546.}
\email{thasebe@math.sci.hokudai.ac.jp}

\author{Micha\"{e}l Ulrich} 
\address[Micha\"{e}l Ulrich]{Laboratoire de Math\'{e}matiques de Besan\c{c}on, Université de Bourgogne-Franche-Comté,16, route de Gray, 25000 Besan\c{c}on, France}

\email{michael.ulrich@univ-fcomte.fr}

\thanks{Final version published at \url{https://doi.org/10.1016/j.jfa.2023.109907}}
\thanks{\textcopyright 2023. This manuscript version is made available under the CC-BY-NC-ND 4.0 license\\ \url{https://creativecommons.org/licenses/by-nc-nd/4.0/}.}
\begin{document}
\maketitle

\begin{abstract}
We attack the classification problem of multi-faced independences, the first non-trivial example being Voiculescu's bi-freeness. While the present paper does not achieve a complete classification, it formalizes the idea of lifting an operator on a pre-Hilbert space in a ``universal'' way to a larger product space, which is key for the construction of (old and new) examples. It will be shown how universal lifts can be used to construct very well-behaved (multi-faced) independences in general. Furthermore, we entirely classify universal lifts to the tensor product and to the free product of pre-Hilbert spaces. Our work brings to light surprising new examples of two-faced independences. Most noteworthy, for many known two-faced independences, we find that they admit continuous deformations within the class of two-faced independences, showing in particular that, in contrast with the single faced case, this class is infinite (and even uncountable).

\end{abstract}

\section{Introduction}

\subsection{Background}
Free independence (or freeness) for non-commutative random variables was introduced by Voiculescu \cite{Voiculescu85}, originally in order to understand free group factors better, especially to solve the isomorphism problem (see e.g.\ \cite{SpeicherCours20182019}, especially p.\ 4 and p.\ 90, and \cite[Preface]{VoiculescuDykemaNica}). It has since attracted much attention because of its various connections and applications, in particular to operator algebras and random matrices \cite{MingoSpeicher,NicaSpeicher,VoiculescuDykemaNica}. 
 Free independence also revealed the fact that a natural notion of independence is not unique in non-commutative probability (because there is already ``tensor independence'' which is more or less the standard independence for quantum mechanical systems). Later, Speicher and Woroudi \cite{SpeicherWoroudi97} formulated boolean independence (which was implicitly discussed earlier by von Waldenfels \cite[Section II.2]{vWaldenfels73}, \cite{vWaldenfels75} and Bo\.zejko \cite{Bozejko86}), and then Muraki formulated monotone and antimonotone independence \cite{Muraki01}. Many attempts were made to interpolate, unify or generalize those notions of independence (e.g.\ \cite{BozejkoLeinertSpeicher,Mlotkowski04,Lenczewski10,Wysoczanski10,Muraki13t,JekelLiu_tree_independence,Skoufranis_AIHP}).  
 
 As more and more notions of independence were proposed, it was a natural direction of research to axiomatize the notions of independence and even classify them --- see e.g.\ \cite{Speicher97,BenGhorbalSchurmann, Muraki02, Muraki03,Muraki13,Lehner04,GerholdLachs15,HasebeLehner17} for such attempts. 
 Among them, the present paper is closest to \cite{Muraki13}, which verified that there are only five notions of independence (tensor, free, boolean, monotone, antimonotone) that satisfy a set of axioms including certain universality, associativity, and positivity conditions (see Definition \ref{defi:universal_product} with $m=1$ and Remark \ref{rmq:UP-linear-functionals} below).

For reference, these five notions of independence are listed below. Let $( A, \phi)$ be a $\ast$-probability space (i.e.\ $A$ is a $\ast$-algebra and $\phi$ is a state on it, see Definition \ref{defi:state} below for details) and $ A_1$ and $ A_2$ be $\ast$-subalgebras of $ A$. 
  
 \begin{enumerate}[label=\rm(\roman*)]
 \item $ A_1$ and $ A_2$ are \emph{tensor independent} if  for every $n\in\N, k_1,k_2, \dots, k_n \in \{1,2\}$ with $k_1 \ne k_2 \ne \cdots \ne k_n$ (meaning that $k_i \ne k_{i+1}$ for all $i\in \{1,2,\dots, n-1\}$) and every $(a_1,a_2,\dots, a_n) \in  A_{k_1} \times  A_{k_2} \times\cdots \times  A_{k_n}$, we have
 \[
 \phi(a_1 a_2 \dotsm a_n) = \phi(a_1 a_3 a_5 \dotsm ) \phi(a_2 a_4 a_6 \dotsm ). 
 \]
  \item Let $A\langle 1 \rangle$ be the unitization of $A$ with unit denoted by $1$ and $\tilde \phi\colon A\langle 1 \rangle \to \mathbb C$ be the unital extension of $\phi$. Let $\tilde A_i$ be the $\ast$-subalgebra of $A\langle 1 \rangle$ generated by $1$ and $A_i$ for $i \in \{1,2\}$. We say that $A_1$ and $ A_2$ are \emph{freely independent (or free)} if  for every $n\in\N, k_1,k_2, \dots, k_n \in \{1,2\}$ with $k_1 \ne k_2 \ne \cdots \ne k_n$ and every $(a_1,a_2,\dots, a_n) \in  \tilde A_{k_1} \times \tilde A_{k_2}\times\cdots \times  \tilde A_{k_n}$ with $\tilde\phi(a_i)=0$ for all $i \in \{1,2,\dots,n\}$, we have
 \[
\tilde \phi(a_1 a_2 \dotsm a_n) = 0. 
 \]
 
\item The ordered pair $(A_1,A_2)$ is \emph{monotonically independent} if for every $n\in\N, k_1,k_2, \dots, k_n \in \{1,2\}$ with $k_1 \ne k_2 \ne \cdots \ne k_n$ and every $(a_1,a_2,\dots, a_n) \in  A_{k_1} \times  A_{k_2} \times\cdots \times  A_{k_n}$, we have
 \[
 \phi(a_1 a_2 \dotsm a_n) = \phi(a_p)\phi(a_1 a_2 \dotsm a_{p-1} a_{p+1}  \dotsm a_n) 
 \]
 whenever $p$ satisfies $k_p=2$. 
 
 \item The ordered pair $(A_1,A_2)$ is \emph{antimonotonically independent} if its flip $(A_2,A_1)$ is monotonically independent. 
 
\item  $ A_1$ and $ A_2$ are \emph{boolean independent} if  for every $n\in\N, k_1,k_2, \dots, k_n \in \{1,2\}$ with $k_1 \ne k_2 \ne \cdots \ne k_n$ and every $(a_1,a_2,\dots, a_n) \in  A_{k_1} \times  A_{k_2} \times\cdots \times  A_{k_n}$, we have
 \[
 \phi(a_1 a_2 \dotsm a_n) = \phi(a_1)\phi(a_2) \dotsm \phi(a_n). 
 \]
 \end{enumerate}
 Among the five notions, monotone and antimonotone independences have a non-symmetric nature between the subalgebras $A_1$ and $A_2$. This is why one needs to formulate the independence for ordered pairs.  
We note that for elements $a_1,a_2 \in A$ the independence of $a_1$ and $a_2$ (or the ordered pair $(a_1,a_2)$) is defined by applying the above definitions to the $\ast$-subalgebras $A_1$ and $A_2$ generated by $a_1$ and $a_2$, respectively.  

All independences we consider can be phrased as coincidence of (non-commutative) ``joint distribution'' with a certain product of ``marginal distributions'', let us roughly sketch what we mean by this.
A ($*$-algebraic) \emph{random variable} is a $*$-homomorphism $j\colon B\to A$, where $B$ is a $*$-algebra and $(A,\Phi)$ is a $*$-probability space. The most important special case is that the $*$-homomorphism is of the form $j_{a_1,\ldots, a_m}\colon \mathbb C\langle{x_1,\ldots, x_m}\rangle\to A$, $x_k\mapsto a_k$, for self-adjoint elements $a_1,\ldots a_m\in A$, where $\mathbb C\langle{x_1,\ldots, x_m}\rangle$ denotes the $*$-algebra of polynomials in $m$ non-commuting self-adjoint indeterminates without constant term; in this special case we also speak of $(a_1,\ldots, a_m)\in A^m$ as a random vector (or random variable if $m=1$) and we think of it as a non-commutative analogue of random vector in $\mathbb R^m$. The independences discussed in this article arise in the following way: A product operation $\odot$ is fixed which maps states $\varphi_k$ on $*$-algebras\footnote{Actually, we will need to deal with $*$-algebras $B_k$ that are \emph{multi-faced}, an extra structure defined in Section \ref{sec:states} reflecting that we want to interpret each random variable $j_k\colon B_k\to A$ as a (not necessarily real valued) non-commutative random vector of fixed size, but in essence, the described relationship between independences and products of states remains the same.} $B_k$ ($k\in\{1,2\})$, to a state $\varphi_1\odot \varphi_2$ on the non-unital free product $*$-algebra $B_1\sqcup B_2$; random variables $j_k\colon B_k\to A$ are called \emph{independent} (with respect to $\odot$) if the \emph{joint distribution} $\Phi\circ(j_1\sqcup j_2)$ coincides with the product $(\Phi\circ j_1)\odot(\Phi\circ j_2)$ of \emph{marginal distributions} $\Phi \circ j_k$. One can specialize to embeddings $j_k=\iota_k\colon A_k\hookrightarrow A$ of $*$-subalgebras to recover the five definitions of independence of $*$-subalgebras given above (when the product $\odot$ is chosen correspondingly). One can specialize to random variables $a_k\in A$ as explained above (i.e.\ $j_{a_k}\colon \mathbb C\langle x\rangle \to A, x\mapsto a_k$) to derive the correct definition of independence for elements of $A$. We usually think of an independence and the underlying product of states as the same thing.

Recent progress in non-commutative probability includes the discovery of new notions of multi-faced independence, i.e.\ independence for ``non-commutative random vectors''. Actually, the present paper was initiated with the motivation to classify multi-faced independence (we could not complete classification, but obtained some new examples). The first substantial development in this direction is due to Voiculescu, who introduced the notion of bi-freeness \cite{Voiculescu14}. Then Gerhold introduced bi-monotone independence \cite{Gerhold17p}, Liu introduced free-boolean independence \cite{Liu19} and also free-free-boolean independence \cite{Liu_preprint}. 
On the other hand, Gu and Skoufranis defined bi-boolean independence \cite{GuSkoufranis} and Gu, Hasebe and Skoufranis defined another notion of bi-monotone independence \cite{GuHasebeSkoufranis20} (different from Gerhold's aforementioned definition); however, both of them do not respect the positivity of states as noted in \cite{GuHasebeSkoufranis20}. This means that if we drop the positivity from the definition of independence, then we get various more notions of multi-faced independence. 
An axiomatic (or categorical) formulation of multi-faced independence is discussed in \cite{ManzelSchurmann17} (see also \cite{GerholdLachsSchurmann22}, \cite[Definition 3.3]{Gerhold21p} or Definitions  \ref{defi:universal_product} and \ref{defi:multi-faced_independence} below), which includes all the above mentioned examples.

\subsection{Lifts, products of representations, and products of states}

Each notion of independence among the five, i.e.\ tensor, free, boolean, monotone and antimonotone independence, has a canonical realization by ``lifting up'' (adjointable) operators on (pre-) Hilbert spaces $H_1$ and $H_2$ equipped with a unit vector $\Omega$ (sometimes called vacuum vector) to operators on a larger space $H_1\boxtimes H_2$ for some monoidal product $\boxtimes$\footnote{A monoidal product is, roughly speaking, a product which up to natural isomorphism is associative and has a unit. As we will only deal with very specific examples (the tensor product and the free product of pre-Hilbert spaces), we will not give a precise definition here. Definition \ref{defi:monoidal_product} and the following discussion contain the relevant references and some details.}. 
For example, a canonical operator model for tensor independence can be provided by the tensor lift $(\lambda^{\rm tensor},\rho^{\rm tensor})$ consisting of mappings $\lambda_{H_1,H_2}^{\rm tensor}\colon \Lop(H_1) \to \Lop(H_1\otimes H_2)$ and $\rho_{H_1,H_2}^{\rm tensor}\colon \Lop(H_2) \to \Lop(H_1\otimes H_2)$ defined by 
\begin{align}\label{eq:tensor_lift}
\lambda_{H_1,H_2}^{\rm tensor}(T_1)= T_1 \otimes \id \qquad \text{and} \qquad \rho_{H_1,H_2}^{\rm tensor}(T_2) = \id\otimes T_2
\end{align}
for $T_k\in \Lop(H_k), k=1,2$. Then the operators $\lambda_{H_1,H_2}^{\rm tensor}(T_1)$ and $\rho_{H_1,H_2}^{\rm tensor}(T_2)$ are tensor independent with respect to the vacuum state (i.e.\ the vector state assciated with $\Omega\otimes \Omega$) on $\Lop(H_1\otimes H_2)$. This model is the non-commutative analogue of the canonical construction of independent random variables on the product of probability spaces.  Similarly,  boolean,  monotone and antimonotone independence can be canonically realized by the respective lifts   
\begin{align}\label{eq:boole_lift}
\lambda_{H_1,H_2}^{\rm boole}(T_1)= T_1 \otimes P_{\Omega} \qquad \text{and} \qquad
\rho_{H_1,H_2}^{\rm boole}(T_2) = P_{\Omega}\otimes T_2, 
\end{align}
\begin{align}\label{eq:monotone_lift}
\lambda_{H_1,H_2}^{\rm mono}(T_1)= T_1 \otimes P_{\Omega} \qquad \text{and} \qquad 
\rho_{H_1,H_2}^{\rm mono}(T_2) = \id\otimes T_2, 
\end{align}
\begin{align}\label{eq:antimonotone_lift}
\lambda_{H_1,H_2}^{\rm amono}(T_1)= T_1 \otimes \id \qquad \text{and} \qquad
\rho_{H_1,H_2}^{\rm amono}(T_2) = P_{\Omega}\otimes T_2,   
\end{align}
where $P_\Omega$ is the orthogonal projection onto $\mathbb C \Omega$ in the space $H_1$ or $H_2$, respectively. 
On the other hand, free independence is constructed on the free product Hilbert space $H_1 \AST H_2$. 
There are canonical unitaries $U_{H_1,H_2} \colon H_1 \AST H_2 \to H_1 \otimes H(1)$ and $V_{H_1,H_2} \colon H_1 \AST H_2 \to H_2 \otimes H(2)$ with certain complementary spaces $H(k), k=1,2$; then the (left) free lift $(\fllift^{\rm free}, \frlift^{\rm free})$ consists of $\fllift^{\rm free}_{H_1,H_2}\colon \Lop(H_1) \to \Lop(H_1\AST H_2)$ and $\frlift^{\rm free}_{H_1,H_2}\colon \Lop(H_2) \to \Lop(H_1\AST H_2)$ defined by 
\begin{align}\label{eq:free_lift}
\fllift^{\rm free}_{H_1,H_2} (T_1) 
&= (U_{H_1,H_2})^* (T_1 \otimes \id)U_{H_1,H_2} \quad \text{and} \\
\frlift^{\rm free}_{H_1,H_2} (T_2) 
&= (V_{H_1,H_2})^* (T_2 \otimes \id)V_{H_1,H_2}, \notag 
\end{align}
which are actions on $H_1\AST H_2$ ``from the left side''\footnote{
We are aware that the notation $(\fllift,\frlift)$ is somewhat confusing in conjunction with the free product of Hilbert spaces, which has an intrinsic left-right symmetry.
A common notation in free probability is $(\lambda_1, \lambda_2)$ rather than something like $(\fllift,\frlift)$, and is  $(\rho_1,\rho_2)$ rather than something like $(\bfllift, \bfrlift)$ for the analogous lifts ``from the right''. Despite the possible confusion for free probabilists, we decided to choose the present notation for two reasons: first, we wish to avoid heavy subscripts (or superscripts); second, the present notation is natural for the tensor product. The point is that the free lift consists of actions from the left (or right) side only, while the tensor lift consists of the action on the left component and the action on the right component of $H_1\otimes H_2$.  
}. 
The operators $\fllift^{\rm free}_{H_1,H_2}(T_1) $ and $\frlift^{\rm free}_{H_1,H_2} (T_2)$ are then freely independent with respect to the vacuum state on $\Lop(H_1\AST H_2)$.

With those canonical operator models, we associate products of ($\ast$-) representations. Once we choose one of the five lifts $(\lambda,\rho)$ with the corresponding monoidal product $\boxtimes$, we can associate to any two representations $\pi_k\colon A_k\to \Lop(H_k)$  the representation $\pi_1 \odot \pi_2\colon A_1 \sqcup A_2 \to \Lop(H_1 \boxtimes H_2)$ defined by 
\[
\pi_1 \odot \pi_2 = (\lambda_{H_1,H_2}\circ \pi_1) \sqcup (\rho_{H_1, H_2}\circ \pi_2), 
\]
where $A_1 \sqcup A_2$ is the non-unital free product of algebras.

Finally, any one of the five products of representations yields a notion of \emph{universal independence}, or more precisely, a \emph{universal product of states} on the free product of algebras (see Definition \ref{defi:universal_product}, case $m=1$), in the following way. Given two states $\phi_k$ on $A_k$, take representations $\pi_k$ of $A_k$ on pre-Hilbert spaces $H_k$ with fixed unit vector $\Omega$ such that $\phi_k(\cdot) = \langle \Omega, \pi_k(\cdot)\Omega\rangle$ (e.g.\ GNS representations), and define 
\begin{equation}\label{eq:well-defined}
   (\phi_1 \odot \phi_2)(\cdot)= \langle \Omega, (\pi_1\odot \pi_2)(\cdot)\Omega \rangle
\end{equation}
on $A_1 \sqcup A_2$. Of course, one has to check that this is well-defined, i.e.\ the right hand side of \eqref{eq:well-defined} does not depend on the choice of $\pi_1$ and $\pi_2$ with $\phi_k(\cdot) = \langle \Omega, \pi_k(\cdot) \Omega \rangle$.

To summarize, any of the five notions of independence can be obtained along the scheme 
\begin{equation} \label{scheme}
\text{lift $\to$ product of representations $\to$ universal product of states.}  
\end{equation}

Deriving a product of linear functionals from a product of representations has of course been utilized before, mainly because of its merit to prove positivity. For example, the free product of representations is used to construct the reduced free product of $C^*$-algebras \cite{Avitzour82,Voiculescu85}. Furthermore, Franz used the operator model of monotone independence to give a simple proof of associativity of the monotone product of linear functionals \cite{Franz01}. 
However, so far, there seems to be no contribution to the problem of classifying independences which substantially uses the representation-theoretic viewpoint.

\subsection{Main results}

While universal products of states are rather well understood in the literature, lifts and products of representations are mostly used as ad-hoc tools and have not been investigated so deeply from an axiomatic viewpoint. The first main objective of this paper is to formulate a suitable set of axioms for lifts and products of representations in order that the scheme \eqref{scheme} works well. 
We provide a solution to this problem by introducing notions of \emph{universal lifts} and \emph{universal products of ($\ast$-) representations} referring to a monoidal product (Definitions \ref{defi:universal_lift} and \ref{defi:universal_product_representation}). 
We prove that those two notions are actually in bijection with each other (Theorem \ref{thm:correspondence}) and they give rise to universal products of states (Theorem \ref{thm:repr_state}), thus completing the scheme \eqref{scheme} into the following:
\begin{equation} \label{new_scheme}
\text{universal lift $\rightleftarrows$ universal product of representations $\to$ universal product of states.}  
\end{equation}
Note that the proofs are based only on axioms and are independent of any classification results for those objects.

We then classify the universal lifts to the tensor product (Theorem \ref{thm:tensor_lift}). 
The complete list is given by   
 \begin{align}\label{eq:main_universal_tensor_lift}
      \lambda_{H_1,H_2}^\gamma (X)=X\otimes P_{\Omega} + X_\gamma \otimes P_{\Omega^\perp} \qquad \text{and} \qquad \rho_{H_1,H_2}^\delta(Y) = P_\Omega \otimes Y + P_{\Omega^\perp} \otimes Y_\delta, 
\end{align} 
where $\gamma,\delta \in \mathbb T \cup\{0\}$ are parameters such that either $\gamma=\delta\in\mathbb T$ or at least one of the two vanishes, and where
\begin{equation}\label{eq:hom_gamma} 
\Lop(H)\ni T= 
 \begin{pmatrix} 
 \tau & (t')^* \\ 
 t & \hat T \end{pmatrix} 
 \mapsto 
 T_\gamma \colonequals
 \begin{pmatrix}
 |\gamma|\tau & (\gamma t')^* \\ 
 \gamma t & |\gamma|\hat T 
 \end{pmatrix} 
 \in \Lop(H)
\end{equation}
is a $\ast$-homomorphism defined through the matrix form of operators based on the decomposition $H = \mathbb C\Omega \oplus \hat H$. 
Of course, this classification contains the tensor, monotone, antimonotone and boolean lifts as special cases.

We also classify the universal lifts to the free product (Theorem \ref{thm:free_lift}). The classification list is given by $\{(\fllift^\gamma,\frlift^\delta)\}_{(\gamma,\delta) \in \mathbb T^2 \cup \{(0,0)\}}$ and $\{(\bfllift^\gamma,\bfrlift^\delta)\}_{(\gamma,\delta) \in \mathbb T^2 \cup \{(0,0)\}}$, where  
 \begin{align}
      \fllift_{H_1,H_2}^\gamma (X)
      &=(U_{H_1,H_2})^*(X\otimes P_{\Omega} + X_\gamma \otimes P_{H(1) \ominus (\mathbb C\Omega)})U_{H_1,H_2} \quad \text{and}  \label{eq:main_universal_free_lift} \\
      \frlift_{H_1,H_2}^\delta(Y) &= 
      (V_{H_1,H_2})^* (Y \otimes P_\Omega  + Y_\delta \otimes P_{H(2) \ominus (\mathbb C\Omega)})V_{H_1,H_2}, \notag 
\end{align} 
and $\bfllift_{H_1,H_2}^\gamma,\bfrlift_{H_1,H_2}^\delta$ are analogous actions from the right side. The definition is apparently motivated by the tensor case, but the allowed set of parameters $(\gamma,\delta)$ is the set $\mathbb T^2 \cup\{(0,0)\}$, which is, somewhat surprisingly, different from the tensor case. 

We only study lifts to the tensor product and to the free product. Having a prescribed space where to lift makes a big difference when one tries to determine all universal lifts. On the one hand, this is what allows us to actually classify all universal lifts to the tensor and free products. On the other hand, this restriction might well be the reason why we cannot get a perfect correspondence between universal lifts and universal products of states in scheme \eqref{new_scheme}.  

The second main contribution is the construction of multi-faced universal products of states (roughly, independence for random vectors), based on the established scheme \eqref{new_scheme} and the new families of universal lifts \eqref{eq:main_universal_tensor_lift} and \eqref{eq:main_universal_free_lift} (Sections \ref{sec:tensor_lift} and \ref{sec:free_lift}). Note that the new universal lifts above can not  produce any new single-faced  universal products of states because the latter objects are already classified into the five kinds in \cite[Theorem 3.7]{Muraki13}.\footnote{\label{fn:Muraki}\revise{\cite[Theorem 3.7]{Muraki13} only states that there are exactly five \emph{non-degenerate} positive universal products. For the proof, Muraki refers to Ben Ghorbal and Sch{\"u}rmann \cite[Theorem 2.5]{BenGhorbalSchurmann05}, who show that positivity implies Muraki's condition N4 (sometimes called \emph{stochastic independence} or \emph{factorization on length 2}); for some reason, Ben Ghorbal and Sch{\"u}rmann formulate the statement under the condition of non-degenerateness, however, looking at their proof and the definition of non-degenerateness, it becomes apparent that this condition is superfluous. The \emph{degenerate} product is not positive, every positive universal product fulfills N4 and, thus, Muraki's theorem yields that there are exactly five positive universal products.}} However, in the multi-faced setting, our universal lifts turn out to yield new universal products of states. What happens is that once we focus on a single component (also called a face) of random vectors, the independence relation satisfied by them should be one of the five kinds, but there is a possibility that a new ``universal correlation'' might appear between different components of random vectors. Indeed, we find such a new type of correlation in which the deformation parameters $\gamma$ and $\delta$ are involved. 

The multi-faced universal products of states thus obtained can be non-symmetric (such as monotone independence in the single-faced case). We identify which ones are symmetric and which ones are not (Propositions \ref{prop:symmetry_tensor} and \ref{prop:symmetry_free}). 

\begin{rmq}
The possible existence of continuous families of two-faced independences has first been recognized by Var\v{s}o during his PhD studies \cite{Varso} on a combinatorial approach to multi-faced symmetric independence. Var\v{s}o gives a set of axioms for classes of (multi-faced) set partitions which assure that there is a (necessarily unique) independence such that its moment-cumulant-relation is determined by the class of partitions in question \cite[Def.\ 3.4.9 and Theorem 3.4.32]{Varso}. In the single-faced case, every positive and symmetric universal product (tensor, free, and boolean) is such a \emph{partition induced} universal product (the moment-cumulant relations are given by all partitions, all noncrossing partitions, and all interval partitions, respectively); this is closely related to Speicher's approach to classification using the \emph{highest coefficients} in \cite{Speicher97}, which Speicher proves must be 0 or 1 and, therefore, are determined by the class of partitions for which they do not vanish. When trying to prove the same in the two-faced case, Var\v{s}o recognized that analogous arguments do not force highest coefficients to be 0 or 1. The class of partitions for which the highest coefficients do not vanish must fulfill Var\v{s}o's axioms, but for some classes of partitions a free choice of one parameter $q\in\mathbb C$ seemed possible, which then determines all other highest coefficients; this is the case for the classes of all partitions, of all noncrossing partitions and of all bi-noncrossing partitions, corresponding to tensor, free, and bi-free independence, respectively \cite[Remark 5.2.29]{Varso}. It is not decided by Var\v{s}o whether there really exist (associative and positive) universal products for $q\neq 1$, only $|q|>1$ is ruled out by positivity.\footnote{During the revision process of this article, Gerhold and Var\v{s}o \cite{GerholdVarso23p} improved the results and techniques of \cite{Varso} and showed that there exists an associative universal product in the mentioned cases if and only if $|q|\in\{0,1\}$.} Still, Var\v{s}o's discoveries were very motivating for us, making it seem likely that our representation theoretic approach might enable us to find previously unknown independences. In the symmetric case, we indeed find deformations within the positive universal products of the same three products as Var\v{s}o's considerations suggest, namely tensor, free, and bi-free. Our results show that the deformed product is positive as long as the deformation parameter is of modulus 1. The case $|q|<1$, still open in \cite{Varso}, was settled in the negative in \cite{GerholdVarso23p}. The question of positivity for the new examples of partition induced universal products appearing in Var\v{s}o's work (most notably a mixture of tensor and free independence) remains open. On the other hand, our method is not restricted to the symmetric case and also brings to light new non-symmetric independences.
\end{rmq}

\subsection{General structure of the article}

Sections \ref{sec:lift_repr} -- \ref{sec:free_lift} are devoted to the formulation of new concepts and proofs of the main results mentioned above. On the other hand, Section \ref{sec:setup} collects basic notions and notations. In this paper we mostly work on algebraic objects which are not very commonly used in functional analysis, e.g.\ pre-Hilbert spaces rather than Hilbert spaces, and (possibly non-unital) $\ast$-algebras rather than $C^*$-algebras. Accordingly, there are several technical issues which the reader might not be familiar with; therefore we will include in Section \ref{sec:setup} some elementary remarks which might help reading.

\section{Setup, notation and remarks}\label{sec:setup}

This paper focuses on possibly non-unital $\ast$-algebras over $\mathbb C$ and their  $\ast$-representations on complex pre-Hilbert spaces with a fixed unit vector. This rather algebraic setup causes some technical differences from the more common setting of $\ast$-representations of $C^*$-algebras on Hilbert spaces in functional analysis.

To begin, as $\ast$-algebras over $\mathbb C$ are not assumed to be unital, one needs to be careful about positivity. The following definition is based on \cite{Gerhold21p,Lachs15}. Note that the term ``state'' is called  ``strongly positive linear functional'' in \cite{Lachs15} and ``restricted state'' in \cite{Gerhold21p}, but we simply call it ``state'' in this paper. 

\begin{defi}\label{defi:state}
Let $A$ be a $\ast$-algebra and $\phi$ be a linear functional on $A$. 

\begin{enumerate}[label=\rm(\arabic*)]
\item\label{item:positive} $\phi$ is said to be \emph{positive} if $\phi(a^*a) \ge 0$ for all $a\in A$. 
\item\label{item:restr} $\phi$ is called a \emph{state} if its unital extension to the unitized $\ast$-algebra $A\langle 1 \rangle=\mathbb C1\oplus A$ is positive. 
\end{enumerate}

The pair $(A,\phi)$ is called a \emph{$\ast$-probability space} if $\phi$ is a state. 
\end{defi}
If $A$ is unital and $\phi(1_A)=1$, then the above conditions \ref{item:positive} and \ref{item:restr} are equivalent --- see Remark \ref{rmq:restricted_state}. 
On the other hand, on a non-unital $\ast$\hyp{}algebra $A$ there may exist positive linear functionals which do not have a positive extension to $A\langle 1\rangle$; in particular they are not states. 

Be aware that the convention has the slightly confusing consequence that linear functionals on unital algebras can be states even if they are not unital. The motivation behind this terminology is that the algebras should always be thought of as ideals in their unitization and the linear functionals as the restriction of their unital extensions. Of course, the restriction of a state to an ideal which has an internal unit may or may not be unital.

\begin{ex}
Let $A= \mathbb C[x]_0$ be the non-unital polynomial algebra equipped with a linear functional defined by $\phi(x)= \alpha$ with some $\alpha \in \mathbb R\setminus\{0\}$ and $\phi(x^n)=0$ for all $n\ge2$. This is obviously positive (and also hermitian), but any extension $\Phi\colon\mathbb C[x] \to \mathbb C$ satisfies that $\Phi((\lambda 1 +x)^\ast(\lambda1+x))= |\lambda|^2\Phi(1) +(\lambda+\bar \lambda)\alpha$. This is non-positive for some $\lambda \in \mathbb R$ sufficiently close to $0$ with $\lambda\alpha<0$. 
\end{ex} 

\revise{Let $K$ be a set. The entries of a tuple $\mathbf k\in K^n$ of \emph{length} $n$ are written $k_1,\ldots, k_n$. A tuple $\mathbf k=(k_1,\ldots, k_n)\in K^n$ is called \emph{alternating} if $k_1\neq\ldots\neq k_n$, which is an abbreviation for $k_i\neq k_{i+1}$ for all $i=1,\ldots, n-1$. The set of all alternating tuples of length $n$ is denoted $K^{n}_{\mathrm{alt}}$.
}

The free product of $\ast$-algebras is defined as the $*$-algebra with underlying vector space
\[
  A_1 \sqcup A_2 =
  \bigoplus_{n\in \N, \mathbf k\in \{1,2\}^n_{\mathrm{alt}}}
  A_{k_1} \otimes A_{k_2} \otimes \cdots \otimes A_{k_n}
\]
and with multiplication and involution determined by ($a_i\in A_{k_i}, b_j\in A_{\ell_j}$ for $i=1,\ldots,n$, $j=1,\ldots,m$)
\begin{align*}
    (a_1 \otimes \cdots \otimes a_n)(b_1\otimes\cdots \otimes b_m)&=\begin{cases}
        a_1 \otimes \cdots \otimes a_nb_1\otimes\cdots \otimes b_m& \text{if }k_n=\ell_1,\\
        a_1 \otimes \cdots \otimes a_n\otimes b_1\otimes\cdots \otimes b_m & \text{if }k_n\neq \ell_1, 
    \end{cases} \\
    (a_1 \otimes \cdots \otimes a_n)^* &= a_n^* \otimes \cdots \otimes a_1^*.
\end{align*} 
This object can be characterized by the following universal property in the category of $\ast$-algebras with $\ast$-homomorphisms as arrows: for any $j_k\colon A_k \to B$ $(k=1,2)$ there exists a unique $j \colon A_1 \sqcup A_2 \to B$ such that $j \circ i_{k}^{A_1,A_2} = j_k$ for $k=1,2$, where $i_k^{A_1,A_2}\colon A_k \hookrightarrow A_1 \sqcup A_2$ are the canonical embeddings. 
The unique $\ast$-homomorphism $j$ is denoted by $j_1\sqcup j_2$.  
In particular, for two $\ast$-homomorphisms $j_k\colon A_k \to B_k$, the $\ast$-homomorphism $(i_1^{B_1,B_2}\circ j_1) \sqcup (i_2^{B_1,B_2} \circ j_2)\colon A_1\sqcup A_2 \to B_1 \sqcup B_2$ can be defined, which will be also denoted by $j_1 \sqcup j_2$ for simplicity (the intended codomain should be obvious from the context).  
We often identify $A_1 \sqcup (A_2\sqcup A_3)$ and $(A_1 \sqcup A_2)\sqcup A_3$ and denote them as $A_1\sqcup A_2 \sqcup A_3$, and also identify $A \sqcup \{0\} = A = \{0\} \sqcup A$. 

\revise{For a $\ast$-algebra $A$ and a subset $M\subset A$, we denote by $\Alg{M}$ the $\ast$-algebra generated by $M$. Note that we do not automatically include $1_A$ in $\Alg{M}$ even if $A$ is unital. If $M=\{a_1,\ldots, a_n\}$ is a finite set, we simply write $\Alg{a_1,\ldots, a_n}$ instead of $\Alg{\{a_1,\ldots, a_n\}}$.}

For pre-Hilbert spaces $H$ and $G$ over $\mathbb C$, an operator $T\colon H\to G$ is called \emph{adjointable} if there exists a (necessarily unique) formal \emph{adjoint} $T^*\colon G\to H$, i.e.\ $\langle T x,y\rangle=\langle x,T^* y\rangle$ for all $x\in H,y\in G$. 
We denote by $\Lop(H,G)$ the set of adjointable linear operators from $H$ to $G$.  In particular, $\Lop(H,H)$ forms a $\ast$-algebra and will be denoted by $\Lop(H)$. In this paper we always assume that every pre-Hilbert space is equipped with a specified unit vector denoted $\Omega$. 

A $\ast$\hyp{}representation of a $\ast$-algebra $A$ is a $\ast$-homomorphism $\pi\colon A\to \Lop(H)$, where $H$ is a pre-Hilbert space. Let $(A,\phi)$ be a $\ast$-probability space. If $A$ and $\phi$ are unital, then by the standard arguments $A$ admits a $\ast$-representation $\pi\colon A \to \Lop(H)$ such that $\pi(A)\Omega=H$ and $\phi(a)=\langle \Omega, \pi(a)\Omega\rangle$ for all $a\in A$; see e.g.\ \cite[Theorem 1.19]{HoraObata07}. Such a representation is called GNS representation and is unique up to unitary equivalence.

\begin{rmq}\label{rmq:restricted_state}
    A linear functional $\phi$ on a $\ast$-algebra $A$ is a state if and only if there exists a $\ast$\hyp{}representation $\pi\colon A\to \Lop(H)$ such that $\phi(a) = \langle \Omega, \pi(a)\Omega\rangle$ for all $a\in A$. 
    For the proof of the sufficiency, one needs to just observe that the unital extension of any $\ast$\hyp{}representation of $A$ to the unitized $\ast$-algebra $A\langle1\rangle$ is a $\ast$\hyp{}representation. For the necessity, given a state on $A$, extend it to a positive linear functional on $A\langle 1 \rangle$, take a GNS representation of the extended functional, and restrict it to the original $\ast$-algebra $A$. 

    The above arguments also imply that a positive unital linear functional on a unital $\ast$-algebra is a state.   

\end{rmq}

One technical issue about pre-Hilbert spaces is orthogonal decomposition. For a closed subspace $E$ of a pre-Hilbert space $H$, it may or may not be the case that $E \oplus E^\perp$ is equal to $H$, where $E^\perp =\{x\in H: \langle x, y\rangle=0 \text{~for all~}y\in E \}$. This is closely related to adjointability of linear mappings; in fact, the orthogonal decomposition $H=E \oplus E^\perp$ holds if and only if the embedding $W\colon E \hookrightarrow H$ is adjointable. Note that if $E$ is of finite dimensions, then $H=E \oplus E^\perp$ holds because we can explicitly construct the orthogonal projection $P_E$ onto $E$ by taking an orthonormal basis of $E$, and then $E^\perp = (\id-P_E)H$. When $H$ decomposes as $H=E\oplus E^\perp$, we will also write $H \ominus E$ instead of $E^\perp$. The orthogonal decomposition $H = \mathbb C \Omega \oplus (H \ominus \mathbb C \Omega)$  will often be used, and usually $H \ominus \mathbb C \Omega$ will be written as $\hat H$ for brevity.

\begin{ex}
Let $H= C[-1,1]$ equipped with $L^2$ norm and the unit vector $\Omega\equiv 1$, and let $E= \{f \in H: f=0 \textnormal{~on~} [-1,0]\}$. Then $E$ is a closed subspace  and $E^\perp=\{f \in H: f=0 \textnormal{~on~} [0,1]\}$. The function $\Omega$ cannot be written as a sum of functions in $E$ and $E^\perp$.  
\end{ex}

We denote by $\PreHilb$ the category of pre-Hilbert spaces (with a fixed unit vector $\Omega$ as before). The arrows of this category are possibly non-adjointable $\Omega$-preserving isometries, i.e.\  norm-preserving linear maps $\UU\colon H\to G$ with $H\ni\Omega\mapsto \Omega\in G$. Then the arrows also preserve the orthogonal complements to $\Omega$, i.e.\ $\UU(\hat H) \subset \hat G$. To see this, note that $\langle \UU x, \UU y \rangle = \langle x,y\rangle$ for all $x,y \in H$ by the polarization identity; then, for each $x \in \hat H=H\ominus \mathbb C \Omega$, we have 
\[
\langle \UU x, \Omega \rangle = \langle \UU x, \UU\Omega \rangle  = \langle x, \Omega \rangle =0, 
\]
as desired. 

The reader might wonder why we are not assuming the existence of adjoints for arrows in $\PreHilb$ although most of the objects we discuss admit the structure of adjoints (e.g.\ $\ast$-algebras and $\ast$-representations). The main reason is that some crucial intertwiners for $\ast$-representations can be non-adjointable; see the proof of Lemma \ref{lm:repr_state} and Remark \ref{rmq:nonadjointabe}.

Here we prove one fact which will be used in Section \ref{sec:lift_repr}. 

\begin{prop}\label{pd-restriction}
Let $H$ and $G$ be pre-Hilbert spaces,  $\UU\colon H\rightarrow G$ be a possibly non-adjointable isometry, and $\SS\in \Lop(G)$. Then the following are equivalent.
\begin{enumerate}[label=\rm(\arabic*)]
    \item $\UU(H)$ is invariant for $\SS$ and $\SS^*$.
    \item\label{item:restriction2} There is a $\TT\in \Lop(H)$ with $\UU \TT=\SS \UU$ and $\UU \TT^* = \SS^*\UU$.
    \item\label{item:restriction3} There is a $\ast$-homomorphism $\pi\colon \Alg{\SS}\to \Lop(H)$ such that $\UU\pi(\SS)= \SS \UU$ and $\UU\pi(\SS^*)=\SS^*\UU$. 
\end{enumerate}
Moreover, the operator $T$ in \ref{item:restriction2} and  $\ast$-homomorphism $\pi$ in \ref{item:restriction3} are unique. 
\end{prop}

\begin{proof}
\subcasens{(1)$\implies$(2)} Simply define $\TT\colonequals \UU^{-1}\SS \UU$ and check that for all $x,y\in H$,
    \[\langle x,\TT y \rangle= \langle x, \UU^{-1} \SS \UU y\rangle = \langle \UU x,\SS \UU y\rangle = \langle \SS^* \UU x, \UU y\rangle = \langle \UU^{-1}\SS^*\UU x, y\rangle,\]
    so that $\TT^*=\UU^{-1}\SS^*\UU$ and $\UU\TT^*=\SS^*\UU$.  
\subcase{(2)$\implies$(3)} Assume $\UU\TT=\SS\UU$ and $\UU\TT^*=\SS^*\UU$. Then $\pi(p(\SS,\SS^*))\colonequals p(\TT,\TT^*)$ (for $p$ a non-commutative polynomial) well-defines a $*$-homomorphism as sought; indeed if $p(\SS,\SS^*)=0$, then $0=p(\SS,\SS^*)\UU=\UU p(\TT,\TT^*)$, so $p(\TT,\TT^*)=0$ ($\UU$ is an injective map). 
    \subcase{(3)$\implies$(1)} If $\UU\pi(\SS)=\SS\UU$, then $\SS(\UU x)=\UU\pi(\SS)x\in \UU(H)$ for all $x\in H$; if $\UU\pi(\SS^*)=\SS^*\UU$, then $\SS^*(\UU x)=\UU\pi(\SS^*)x\in \UU(H)$.
    
    \subcase{Uniqueness} The uniqueness of $T$ and $\pi$ is obvious from the formulas $WT=SW$ and $W\pi(S)=S W$ and the injectivity of $W$. 
\end{proof}

\section{Universal lifts and products of representations}\label{sec:lift_repr}

This section is devoted to introducing universal lifts and universal products of $\ast$-representations referring to a fixed monoidal product on $\PreHilb$. Actually, we will work only on the tensor product $\otimes$ and the free product $\AST$ in Sections \ref{sec:tensor_lift} and \ref{sec:free_lift}, but in the present section and in Section \ref{sec:states}, we treat more general monoidal products for possible future research. 

The following definition is just for reference; below it we sketch the structure that will be used. 

\begin{defi}\label{defi:monoidal_product}
    A \emph{monoidal product with embeddings} on a category $\mathcal C$ is a bifunctor $\boxtimes\colon \mathcal C\times \mathcal C\to \mathcal C$ (together with natural associativity and unitality isomorphisms) such that $(\mathcal C,\boxtimes)$ is a monoidal category in the sense of \cite[Section VII.1]{MacLane} whose unit object $E$ is an initial object in $\mathcal C$. 
\end{defi}

As a consequence \cite[Theorem 3.5(b)]{GerholdLachsSchurmann22}, a monoidal product with embeddings gives rise to (unique) canonical morphisms $I_{A,B}\colon A\cong A\boxtimes E\to A\boxtimes B$, $J_{A,B}\colon B\cong E\boxtimes B\to A\boxtimes B$ which constitute \emph{compatible inclusions} in the sense of \cite[Definition 3.4]{GerholdLachsSchurmann22}.

For convenience of the reader, we summarize what this means in practice on the category $\PreHilb$.
\begin{itemize}
    \item (Arrow/functoriality) Any two $\Omega$-preserving isometries $\UU_k\colon H_k \to G_k$ $(k\in\{1,2\})$ yield an $\Omega$-preserving isometry $\UU_1 \boxtimes \UU_2 \colon H_1 \boxtimes H_2 \to G_1 \boxtimes G_2$.
    \item (Embedding) Any pair $(H_1,H_2)$ admits specified $\Omega$-preserving isometries $I_{H_1,H_2}\colon H_1 \to {H_1 \boxtimes H_2}$ and $J_{H_1,H_2}\colon H_2 \to {H_1 \boxtimes H_2}$. Through these isometries we regard $H_1$ and $H_2$ as subspaces of $H_1\boxtimes H_2$. 
\item (Unitality) The embeddings $I_{H_1, \mathbb C \Omega}\colon H_1 \to H_1 \boxtimes {\mathbb C\Omega}$ and $J_{\mathbb C \Omega, H_2} \colon H_2 \to {\mathbb C\Omega \boxtimes H_2}$ are isomorphisms (i.e.\ $\Omega$-preserving unitaries).  Through these isomorphisms we will   identify  $H \boxtimes \mathbb C\Omega = H = \mathbb C \Omega \boxtimes H$. 
    
\item (Associativity) Any triple $(H_1,H_2,H_3)$ admits a specified isomorphism $A_{H_1,H_2,H_3}\colon H_1 \boxtimes (H_2 \boxtimes H_3) \to {(H_1 \boxtimes H_2) \boxtimes H_3}$. Through this arrow we identify the two objects and simply write them as $H_1 \boxtimes H_2 \boxtimes H_3$.

\end{itemize}

Note that the following proposed axioms defining universal lifts are all satisfied by the five lifts \eqref{eq:tensor_lift} -- \eqref{eq:free_lift}. This fact can be directly checked, but with Theorems \ref{thm:tensor_lift} and \ref{thm:free_lift} we actually prove it for the wider families \eqref{eq:main_universal_tensor_lift} and \eqref{eq:main_universal_free_lift}, respectively.

\begin{defi} \label{defi:universal_lift} Let $\boxtimes$ be a monoidal product with embeddings on the category $\PreHilb$. 
    A \emph{left universal lift} (for $\boxtimes$) is a family of $*$-homomorphisms \(\lambda_{H_1,H_2}\colon \Lop(H_1)\to \Lop(H_1\boxtimes H_2)\) with the following properties. 
  
\begin{itemize}
    \item (Left associativity) For any triple of pre-Hilbert spaces $(H_1,H_2,H_3)$ we have \[
    \lambda_{H_1,H_2\boxtimes H_3}=\lambda_{H_1\boxtimes H_2,H_3}\circ\lambda_{H_1,H_2}.
    \]
    \item (Left universality of pre-Hilbert spaces)
 For any $\TT \in \Lop(H_1), \SS \in \Lop(G_1)$ and arrows $\UU_k\colon H_k \to G_k$ in the category $\PreHilb$ $(k=1,2)$ such that   
\[  
 \xymatrix{
    H_1\ar[rr]^{\TT}\ar[d]_{\UU_1}
    && H_1\ar[d]^{\UU_1}\\G_1\ar[rr]_{\SS}&&G_1
    } 
    \qquad \raisebox{-18pt}{\text{and}} \qquad  \xymatrix{
    H_1\ar[rr]^{\TT^*}\ar[d]_{\UU_1}&& H_1\ar[d]^{\UU_1}\\G_1\ar[rr]_{\SS^*}&&G_1,
    } \qquad \raisebox{-18pt}{\text{commute,}}
\]   
we have that
\[
    \xymatrix{
    H_1\boxtimes H_2\ar[rr]^{\lambda_{H_1,H_2}(\TT)}\ar[d]_{\UU_1\boxtimes \UU_2}&& H_1\boxtimes H_2\ar[d]^{\UU_1\boxtimes \UU_2}\\G_1\boxtimes G_2\ar[rr]_{\lambda_{G_1,G_2}(\SS)}&&G_1\boxtimes G_2.
    } \qquad \raisebox{-18pt}{\text{commutes.}}
    \]

 \item (Left restriction) for any pre-Hilbert spaces $H_1, H_2$ and any $T \in \Lop(H_1)$ we have 
 \[
 \Restr{\lambda_{H_1,H_2}(\TT)}{H_1} = \TT.
 \]
\end{itemize}

A \emph{right universal lift} for $\boxtimes$ is a family of $*$-homomorphisms \(\rho_{H_1,H_2}\colon \Lop(H_2)\to \Lop(H_1\boxtimes H_2)\) such that the family with reversed indices $(\rho_{H_1,H_2})_{H_2,H_1}$ is a left universal lift for the opposite monoidal product $\boxtimes^{\mathrm{op}}$. In particular, the associativity condition explicitly reads \[
\rho_{H_1\boxtimes H_2,H_3}=\rho_{H_1,H_2\boxtimes H_3}\circ\rho_{H_2,H_3}
\]
and will be referred to as the \emph{right associativity}.  

A \emph{universal lift} for $\boxtimes$ is a pair $(\lambda,\rho)$ consisting of a left universal lift and a right universal lift such that 
 \begin{itemize}
     \item (Middle associativity) for any triple of pre-Hilbert spaces $(H_1,H_2,H_3)$ we have \[
     \rho_{H_1,H_2\boxtimes H_3}\circ\lambda_{H_2,H_3}=\lambda_{H_1\boxtimes H_2,H_3}\circ\rho_{H_1,H_2}.
     \]
 \end{itemize}
The three associativity conditions (left, right, middle) together will simply be called the \emph{associativity} of $(\lambda,\rho)$. 
\end{defi}

\begin{rmq}
The left restriction property can be replaced with:  
\begin{itemize}
\item (Left unitality) for any $T \in \Lop(H_1)$, we have $\lambda_{H_1,\mathbb C\Omega}(T)=T$ under the identification $H_1\boxtimes\mathbb C\Omega = H_1$. 
\end{itemize}

The left restriction property immediately implies the left unitality; for the converse, use the universality axiom for the isometries $W_1= \mathrm{id} \in \Lop(H_1)$ and $W_2 \colon \mathbb C \Omega \hookrightarrow H_2$ to get $\lambda_{H_1,H_2}(T) \circ(W_1 \boxtimes W_2)= (W_1 \boxtimes W_2) \circ \lambda_{H_1,\mathbb C \Omega}(T) = T$ as a map from $H_1$ to $H_1\subset H_1\boxtimes H_2$. The left hand side is exactly the restriction $\Restr{\lambda_{H_1,H_2}(T)}{H_1}$ because $W_1 \boxtimes W_2\colon H_1 = H_1 \boxtimes \mathbb C\Omega \to H_1 \boxtimes H_2$ is the natural embedding. 
\end{rmq}

\begin{rmq}\label{rem_lifts}
For the tensor product $\boxtimes = \otimes$, a classification result and its proof in Section \ref{sec:tensor_lift} show that the left associativity axiom for left universal lifts actually follows from the other two axioms, the left universality of pre-Hilbert spaces and the left restriction property. By symmetry, a similar fact holds for the right universal lifts. For the free product $\boxtimes = \AST$, however,  the left associativity cannot be dropped --- see Example \ref{ex:monotone_lift}.
\end{rmq}

\begin{defi} \label{defi:universal_product_representation} Let $\boxtimes$ be a monoidal product with embeddings on $\PreHilb$. 
    A \emph{universal product of $\ast$\hyp{}representations} (for $\boxtimes$) is a rule that, given two $\ast$-representations of $*$-algebras on pre-Hilbert spaces $\pi_1\colon A_1\rightarrow \Lop(H_1)$ and $\pi_2\colon A_2\rightarrow \Lop(H_2)$, produces a $*$\hyp{}representation $\pi_1\odot\pi_2\colon A_1\sqcup A_2\rightarrow \Lop(H_1\boxtimes H_2)$, such that the following axioms are verified. 
\begin{itemize}
    \item (Associativity) For any three  $*$\hyp{}representations $\pi_{1}, \pi_2,\pi_3$ we have 
        \[
        \pi_1\odot(\pi_2\odot\pi_3)=(\pi_1\odot\pi_2)\odot\pi_3. 
        \]
        
    \item (Universality of $\ast$-algebras) For any $*$-homomorphisms $j_{k}\colon B_{k}\rightarrow A_{k}$ ($k=1,2$), we have 
    \[
    (\pi_1\odot\pi_2)\circ(j_1\sqcup j_2)=(\pi_1\circ j_1)\odot(\pi_2\circ j_2).
    \]
        
    \item (Universality of pre-Hilbert spaces)  
        For any arrows $\UU_{k}\colon H_{k}\rightarrow G_{k}$ in $\PreHilb$ and any $\ast$\hyp{}representations $\pi_k\colon A_k \to \Lop(H_k)$ and $\sigma_k\colon A_k \to \Lop(G_k)$  such that 
    \[  
        \xymatrix{
        H_k\ar[rrr]^{\pi_k(a)}\ar[d]_{\UU_k}&&& H_k\ar[d]^{\UU_k}\\G_k\ar[rrr]_{\sigma_k(a)}&&&G_k
        }
        \quad \raisebox{-18pt}{\text{commutes for all $a\in A_k$ and $k\in\{1,2\}$,}}
    \]
(i.e.\ $W_k$ is an \emph{intertwiner}), we have that 
    \[  
        \xymatrix{
        H_1\boxtimes H_2\ar[rrr]^{\pi_1\odot \pi_2(a)}\ar[d]_{\UU_1\boxtimes \UU_2}&&& H_1\boxtimes H_2\ar[d]^{\UU_1\boxtimes \UU_2}\\G_1\boxtimes G_2\ar[rrr]_{(\sigma_1\odot \sigma_2) (a)}&&&G_1\boxtimes G_2
        }
          \quad \raisebox{-18pt}{\text{commutes for all $a\in A_1\sqcup A_2$. }}
        \]
        
    \item (Restriction) For any two $\ast$-representations $\pi_k\colon A_k\to\Lop(H_k) ~(k\in\{1,2\})$  we have 
    \[
    \Restr{(\pi_1\odot\pi_2)(a)}{H_k}=\pi_k(a)
    \]
    for any $a\in A_k$ and $k\in\{1,2\}$. 
    \end{itemize}
\end{defi}

\begin{rmq}
    Denote by $\RepCat$ the category with
  
    \begin{description}[labelindent=20pt]
    \item[objects] $*$-representations $\pi\colon A\to \Lop(H)$, written as triple $(\pi,A,H)$, 
    \item[morphisms $(\pi,A,H)\to(\sigma,B,G)$] pairs $(j,W)$ where $j\colon A \to B$ is a $*$-algebra homomorphism and $\UU\colon H\to G$ is an isometry such that 
    \[  
        \xymatrix{
        H\ar[rrr]^{\pi(a)}\ar[d]_{\UU}&&& H\ar[d]^{\UU}\\G\ar[rrr]_{\sigma(j(a))}&&&G
        }
        \quad \raisebox{-18pt}{\text{commutes for all $a\in A$.}}
    \]
    \end{description}  
    The two universality conditions in Definition \ref{defi:universal_product_representation} are equivalent to demanding that the product of morphisms defined as $(j_1,\UU_1)\odot (j_2,\UU_2)\colonequals(j_1\sqcup j_2,\UU_1\boxtimes \UU_2)$ is again a morphism, i.e.\ 
    \begin{multline*}
    \xymatrix{
        H_k\ar[rr]^{\pi_k(a)}\ar[d]_{\UU_k}&& H_k\ar[d]^{\UU_k}\\G_k\ar[rr]_{\sigma_k(j_k(a))}&&G_k
        }
        \quad \raisebox{-18pt}{commutes for all $ a\in A_k,k\in\{1,2\}$}\\
        \raisebox{-18pt}{$\implies$} \quad
        \xymatrix{
        H_1\boxtimes H_2\ar[rrr]^{\pi_1\odot \pi_2(a)}\ar[d]_{\UU_1\boxtimes \UU_2}&&& H_1\boxtimes H_2\ar[d]^{\UU_1\boxtimes \UU_2}\\G_1\boxtimes G_2\ar[rrr]_{(\sigma_1\odot \sigma_2) (j_1\sqcup j_2(a))}&&&G_1\boxtimes G_2
        }
          \quad \raisebox{-18pt}{commutes for all $a\in A_1\sqcup A_2$.}
    \end{multline*}
    Therefore, a universal product of $*$-representations can be viewed as a bifunctor which turns $\RepCat$ into a monoidal category with embeddings.
\end{rmq}

The main result of this section is the following. 

\begin{thm}\label{thm:correspondence} Let $\boxtimes$ be a monoidal product with embeddings on $\PreHilb$. There is a one-to-one correspondence between universal lifts for $\boxtimes$ and universal products of $\ast$\hyp{}representations for $\boxtimes$ given by
    \begin{align*}
        \odot\mapsto (\lambda^{\odot},\rho^{\odot})\quad &\text{with}\quad \begin{cases}\lambda^{\odot}_{H_1,H_2}\colonequals 
        \Restr{\id_{\Lop(H_1)}\odot \id_{\Lop(H_2)}}{\Lop(H_1)},\\ \rho^{\odot}_{H_1,H_2}\colonequals \Restr{\id_{\Lop(H_1)}\odot \id_{\Lop(H_2)}}{\Lop(H_2)}
        \end{cases}
        \end{align*}
for pre-Hilbert spaces $H_1, H_2$, and 
\begin{align*}
        (\lambda, \rho)\mapsto{{}_\lambda\odot_\rho} \quad &\text{with}\quad \pi_1\mathbin{{}_\lambda\odot_\rho}\pi_2\colonequals (\lambda_{H_1,H_2}\circ\pi_1)\sqcup(\rho_{H_1,H_2}\circ\pi_2)
    \end{align*}
    for $\ast$\hyp{}representations $\pi_k\colon A_k \to \Lop(H_k)$ $(k=1,2)$.
\end{thm}

\begin{proof}     First, note that the canonical embedding $A_1= A_1 \sqcup \{0\}\hookrightarrow A_1\sqcup A_2$ can be written as $\mathrm {id}_{A_1}\sqcup 0_{\{0\}\to A_2}$ where $0_{A\to B}$ denotes the trivial $*$-homomorphism $A\to B$ for $\ast$-algebras $A,B$. Therefore, given a universal product of $\ast$\hyp{}representations $\odot$, it follows from universality that
    \[\Restr{\pi_1\odot \pi_2}{A_1}
        = (\pi_1\odot \pi_2)\circ (\mathrm {id}_{A_1}\sqcup 0_{\{0\}\to A_2}) 
        = \pi_1 \odot 0_{\{0\}\to \Lop(H_2)}.
    \]
    Analogously $\Restr{\pi_1\odot \pi_2}{A_2}=0_{\{0\}\to \Lop(H_1)}\odot \pi_2$. As a special case, we find 
\begin{align}
   \lambda^{\odot}_{H_1,H_2}=\Restr{\id_{\Lop(H_1)}\odot \id_{\Lop(H_2)}}{\Lop(H_1)} 
   &= \id_{\Lop(H_1)}\odot 0_{\{0\}\to \Lop(H_2)}  \quad \text{and} \label{eq:left_lift}\\ 
    \rho^{\odot}_{H_1,H_2}=\Restr{\id_{\Lop(H_1)}\odot \id_{\Lop(H_2)}}{\Lop(H_2)}
    &=0_{\{0\}\to \Lop(H_1)}\odot \id_{\Lop(H_2)}.\label{eq:right_lift}
\end{align}
    
    Assume for the moment that the given maps do indeed give universal lifts from universal products and vice versa. Then they are easily seen to be inverse to each other. Indeed, by the preparatory comments above,
    \begin{align*}
      \MoveEqLeft \pi_1\mathbin{{}_{\lambda_{\odot}}\odot_{\rho_{\odot}}}\pi_2=\left(\lambda^{\odot}_{H_1,H_2}\circ\pi_1\right)\sqcup\left( \rho^{\odot}_{H_1,H_2}\circ \pi_2\right)\\
        &=
        \left[\left(\id_{\Lop(H_1)}\odot 0_{\{0\}\to \Lop(H_2)}\right)\circ \left(\pi_1\sqcup 0_{\{0\}\to \{0\}}\right)\right]
        \sqcup
        \left[\left(0_{\{0\}\to \Lop(H_1)}\odot \id_{\Lop(H_2)}\right)\circ \left(0_{\{0\}\to \{0\}}\sqcup \pi_2\right)\right]\\
    &=\left(\pi_1\odot 0_{\{0\}\to \Lop(H_2)}\right)\sqcup\left(0_{\{0\}\to \Lop(H_1)}\odot \pi_2\right) 
        = \left( \Restr{\pi_1\odot \pi_2}{A_1}\right)\sqcup\left( \Restr{\pi_1\odot \pi_2}{A_2}\right) \\
        &= \pi_1\odot \pi_2.  
    \end{align*}
 The last equality holds because both representations of $A_1\sqcup A_2$ have the same restrictions to $A_1$ and $A_2$. The other direction is easier:
    \[
        \lambda_{H_1,H_2}^{{}_\lambda\odot_\rho}=\Restr{\lambda_{H_1,H_2} \sqcup \rho_{H_1,H_2}}{\Lop(H_1)} = \lambda_{H_1,H_2},
        \quad
        \rho_{H_1,H_2}^{{}_\lambda\odot_\rho}=\Restr{\lambda_{H_1,H_2} \sqcup \rho_{H_1,H_2}}{\Lop(H_2)} = \rho_{H_1,H_2}.
    \]

It remains to prove that a universal product of $\ast$\hyp{}representations yields a universal lift and vice versa. We begin with the former claim. Except for the middle associativity, it suffices to focus on the left lift by symmetry. 

\subcasenpds{From universal product to universal lift}
Let $\odot$ be a fixed universal product of $*$\hyp{}representations. We will prove that $(\lambda^{\odot},\rho^{\odot})$, as defined in the statement of the theorem, is a universal lift.

\subcase{Left universality of pre-Hilbert spaces}
 Let $H_1,H_2,G_1,G_2$ be four pre-Hilbert spaces, $\TT \in \Lop(H_1), \SS\in \Lop(G_1)$ and $\UU_k\colon H_k \to G_k$ be possibly non-adjointable $\Omega$-preserving isometries such that $\UU_1 \TT = \SS \UU_1$ and $\UU_1\TT^* = \SS^* \UU_1$. 
 Then let $A_1$ be the $\ast$-algebra generated by $\SS$ and $\sigma_1 \colon A_1 \hookrightarrow \Lop(G_1)$ be the embedding $\sigma_1 (a)=a, a\in A_1$. 
 By Proposition \ref{pd-restriction}, there is a $*$-homomorphism $\pi_1\colon A_1\to \Lop(H_1), \SS\mapsto \TT$ which fulfills $\UU_1 \pi_1 (a) = \sigma_1(a)\UU_1$ for all $a\in A_1$. 
 Also, let $A_2=\{0\}$ be the trivial algebra. 
 Trivially, $\UU_2 0_{A_2\to \Lop(H_2)}(b) = 0_{A_2\to \Lop(G_2)}(b) \UU_2$ for all $b \in A_2$. 
 The two universality axioms of the universal product $\odot$ then imply 
\begin{align*}
  (\UU_1\boxtimes \UU_2)\lambda^{\odot}_{H_1,H_2}(\TT) 
    &= (\UU_1\boxtimes \UU_2) [(\id_{\Lop(H_1)}\odot \id_{\Lop(H_2)}) (\pi_1(\SS))] \\
    &= (\UU_1\boxtimes \UU_2) [ (\id_{\Lop(H_1)}\odot \id_{\Lop(H_2)})\circ(\pi_1\sqcup 0_{A_2\to \Lop(H_2)})(\SS)]\\
    &= (\UU_1\boxtimes \UU_2)[  (\pi_1\odot 0_{A_2\to \Lop(H_2)}) (\SS)] \\
    &= [(\sigma_1\odot 0_{A_2\to \Lop(G_2)}) (\SS)] (\UU_1\boxtimes \UU_2)\\
    &= [(\id_{\Lop(G_1)}\odot \id_{\Lop(G_2)}) \circ(\sigma_1\sqcup 0_{A_2\to \Lop(G_2)}) (\SS)] (\UU_1\boxtimes \UU_2)\\
    &= [(\id_{\Lop(G_1)}\odot \id_{\Lop(G_2)}) (\SS)] (\UU_1\boxtimes \UU_2)\\
    &=\lambda^{\odot}_{G_1,G_2}(\SS) (\UU_1\boxtimes \UU_2). 
\end{align*}

\subcase{Left restriction property}
Let $T\in \Lop(H_1)$ and $h\in H_1$. We have:
\begin{align*}
[\lambda^{\odot}_{H_1,H_2}(T)]h&= [(\id_{\Lop(H_1)}\odot \id_{\Lop(H_2)})(T)]h= [\id_{\Lop(H_1)}(T)]h=Th. 
\end{align*}
Therefore, $\lambda^{\odot}$ satisfies the restriction property.

\subcase{Left associativity}
Note that, by universality of $\ast$-algebras and \eqref{eq:left_lift},
\begin{align*}
    \lambda^{\odot}_{H_1\boxtimes H_2,H_3}\circ\lambda^{\odot}_{H_1,H_2}
    &=  \lambda^{\odot}_{H_1\boxtimes H_2,H_3} \circ (\id_{\Lop(H_1)}\odot 0_{\{0\}\to \Lop(H_2)})\\
    &= (\id_{\Lop(H_1\boxtimes H_2)}\odot 0_{\{0\}\to \Lop(H_3)}) \circ [(\id_{\Lop(H_1)}\odot 0_{\{0\}\to \Lop(H_2)}) \sqcup 0_{\{0\}\to \{0\}}]\\
    &=[\id_{\Lop(H_1\boxtimes H_2)}\circ(\id_{\Lop(H_1)}\odot 0_{\{0\}\to \Lop(H_2)})]\odot (0_{\{0\}\to \Lop(H_3)}\circ 0_{\{0\}\to\{0\}})\\
    &=(\id_{\Lop(H_1)}\odot 0_{\{0\}\to \Lop(H_2)})\odot 0_{\{0\}\to \Lop(H_3)}. 
\end{align*}
Also due to universality of $\ast$-algebras, 
\[
0_{\{0\}\to \Lop(H_2\boxtimes H_3)}=0_{\{0\}\to \Lop(H_2)}\odot 0_{\{0\}\to \Lop(H_3)}.
\]
Combining the above two calculations with associativity of $\odot$ implies the left associativity:
\begin{align*}
\lambda^{\odot}_{H_1,H_2\boxtimes H_3}&= \id_{\Lop(H_1)}\odot \left(0_{\{0\}\to \Lop(H_2)}\odot 0_{\{0\}\to \Lop(H_3)}\right) \\
&=\left(\id_{\Lop(H_1)}\odot 0_{\{0\}\to \Lop(H_2)}\right)\odot 0_{\{0\}\to \Lop(H_3)} \\
&=\lambda^{\odot}_{H_1\boxtimes H_2,H_3}\circ\lambda^{\odot}_{H_1,H_2}.
\end{align*}

\subcase{Middle associativity} The middle associativity is similarly proved:
\begin{align*}
\rho^{\odot}_{H_1,H_2\boxtimes H_3}\circ\lambda^{\odot}_{H_2,H_3}
&= \rho^{\odot}_{H_1,H_2\boxtimes H_3}\circ\left(\id_{\Lop(H_2)}\odot 0_{\{0\}\to \Lop(H_3)}\right)\\
&= \left(0_{\{0\}\to \Lop(H_1)}\odot \id_{\Lop(H_2\boxtimes H_3)}\right)\circ\left[0_{\{0\}\to\{0\}}\sqcup\left(\id_{\Lop(H_2)}\odot 0_{\{0\}\to \Lop(H_3)}\right)\right]\\
&= 0_{\{0\}\to \Lop(H_1)}\odot\left[\id_{\Lop(H_2)}\odot 0_{\{0\}\to \Lop(H_3)}\right]\\
&= \left[0_{\{0\}\to \Lop(H_1)}\odot \id_{\Lop(H_2)}\right]\odot 0_{\{0\}\to \Lop(H_3)}\\
&= \left[\id_{\Lop(H_1\boxtimes H_2)}\odot 0_{\{0\}\to \Lop(H_3)}\right]\circ\left[\left(0_{\{0\}\to \Lop(H_1)}\odot \id_{\Lop(H_2)}\right)\sqcup 0_{\{0\}\to\{0\}}\right]\\
&= \lambda^{\odot}_{H_1\boxtimes H_2,H_3}\circ\rho^{\odot}_{H_1,H_2}. 
\end{align*}
This concludes the middle associativity and hence the first half of the proof.

\subcasenpds{From universal lift to universal product}
We now need to prove that $\mathbin{{}_\lambda\odot_\rho}$ constructed as $\pi_1\mathbin{{}_\lambda\odot_\rho}\pi_2\colonequals (\lambda\circ \pi_1) \sqcup (\rho\circ \pi_2)$ from a universal lift $(\lambda,\rho)$ is a universal product of $*$\hyp{}representations. As $\lambda, \rho$ and $\pi_k$ are $*$-homomorphisms, it is clear that $\pi_1\mathbin{{}_\lambda\odot_\rho}\pi_2$ is a well-defined $*$-homomorphism from $A_1\sqcup A_2$ to $\Lop(H_1\boxtimes H_2)$. We will therefore check the other axioms of a universal product of $\ast$\hyp{}representations.

\subcase{Universality of \texorpdfstring{$\ast$}{*}-algebras for the product} 
 Let $\pi_{k}\colon A_{k}\rightarrow \Lop(H_{k})$ ($k\in\{1,2\})$ be two $*$\hyp{}representations and let $j_{k}\colon B_{k}\rightarrow A_{k}$ be two $*$-homomorphisms. 
We have: 
\begin{align*}
\left(\pi_1\circ j_1\right)\mathbin{{}_\lambda\odot_\rho}\left(\pi_2\circ j_2\right)&=  (\lambda_{H_1,H_2}\circ \pi_1 \circ j_1 )\sqcup (\rho_{H_1,H_2}\circ \pi_2\circ j_2) \\
&= [(\lambda_{H_1,H_2}\circ\pi_1)\sqcup (\rho_{H_1,H_2}\circ\pi_2)]\circ(j_1\sqcup j_2)
= (\pi_1\mathbin{{}_\lambda\odot_\rho}\pi_2)\circ(j_1\sqcup j_2).
\end{align*}

\subcase{Universality of pre-Hilbert spaces for the product} Let $\pi_k\colon A_k \to \Lop(H_k)$ and $\sigma_k\colon A_k \to \Lop(G_k)$ be $\ast$\hyp{}representations and $\UU_{k}\colon (\pi_k,H_{k})\rightarrow (\sigma_k,G_{k})$ possibly non-adjointable $\Omega$-preserving isometric intertwiners ($k=1,2$). 
For every $a \in A_1$ we have $\UU_1\pi_1(a)= \sigma_1(a)\UU_1$ and $\UU_1\pi_1(a)^*= \sigma_1(a)^*\UU_1$, and hence, 
the universality for $\lambda$ implies that 
\begin{align*}
(\UU_1\boxtimes \UU_2)[(\pi_1\mathbin{{}_\lambda\odot_\rho} \pi_2)(a)] 
&= (\UU_1\boxtimes \UU_2)[ \lambda_{H_1,H_2}(\pi_1(a))]  \\
&= [\lambda_{G_1,G_2}(\sigma_1(a))] (\UU_1\boxtimes \UU_2) \\
&= [(\sigma_1\mathbin{{}_\lambda\odot_\rho} \sigma_2)(a)] (\UU_1\boxtimes \UU_2). 
\end{align*}
A similar reasoning for $\rho$ verifies the same identity on any $b \in A_2$. Since $\pi_1\mathbin{{}_\lambda\odot_\rho} \pi_2$ and $\sigma_1\mathbin{{}_\lambda\odot_\rho} \sigma_2$ are both homomorphisms, the identity holds on $A_1\sqcup A_2$.

\subcase{Associativity of the product} Let $\pi_{k}\colon A_{k}\rightarrow \Lop(H_{k})$ ($k\in\{1,2,3\}$) be three $*$\hyp{}representations. For $a\in A_1$ we have
\begin{align*}
\pi_1\mathbin{{}_\lambda\odot_\rho}(\pi_2\mathbin{{}_\lambda\odot_\rho}\pi_3)(a)&= \lambda_{H_1,H_2\boxtimes H_3}\circ\pi_1(a)\\
&= \lambda_{H_1\boxtimes H_2,H_3}\circ\lambda_{H_1,H_2}\circ\pi_1(a)\\
&= \lambda_{H_1\boxtimes H_2,H_3}\circ(\pi_1\mathbin{{}_\lambda\odot_\rho}\pi_2)(a)\\
&= (\pi_1\mathbin{{}_\lambda\odot_\rho}\pi_2)\mathbin{{}_\lambda\odot_\rho}\pi_3(a),
\end{align*}
where we used the left associativity of $\lambda$. Exactly the same reasoning using the right associativity shows that the relation remains true for all $c\in A_3$. Let us now check it for any $b\in A_2$:
\begin{align*}
\pi_1\mathbin{{}_\lambda\odot_\rho}(\pi_2\mathbin{{}_\lambda\odot_\rho}\pi_3)(b)&= \rho_{H_1,H_2\boxtimes H_3}\circ(\pi_2\mathbin{{}_\lambda\odot_\rho}\pi_3)(b)\\
&= \rho_{H_1,H_2\boxtimes H_3}\circ\lambda_{H_2,H_3}\circ\pi_2(b)\\
&= \lambda_{H_1\boxtimes H_2,H_3}\circ\rho_{H_1,H_2}\circ\pi_2(b)\\
&= \lambda_{H_1\boxtimes H_2,H_3}\circ(\pi_1\mathbin{{}_\lambda\odot_\rho}\pi_2)(b)\\
&= (\pi_1\mathbin{{}_\lambda\odot_\rho}\pi_2)\mathbin{{}_\lambda\odot_\rho}\pi_3(b),
\end{align*}
where we have used the middle associativity of the universal lift. As we deal only with $*$-homomorphisms this is enough to conclude the associativity of the universal product of $\ast$\hyp{}representations.

\subcase{Restriction property for the product} With the same notations as before and with $a\in A_1$ and $h\in H_1$ we have $(\pi_1\mathbin{{}_\lambda\odot_\rho}\pi_2)(a)h=\lambda_{H_1,H_2}(\pi_1(a))h=\pi_1(a)h$ by the axiom of left restriction for the lift. Therefore all the axioms for a universal product of $\ast$\hyp{}representations are indeed verified.
\end{proof}

\begin{rmq}
In the above proof, the other axioms are discussed independently of universality axioms of pre-Hilbert spaces. This means that, if we drop the universality axioms of pre-Hilbert spaces from lifts and products, there is still a one-to-one correspondence between lifts and products of $\ast$-representations. However, the universality of pre-Hilbert spaces will play an important role in constructing universal products of states --- see Section \ref{sec:states}.   
\end{rmq}

\section{Multi-faced universal products of states}\label{sec:states}

This section establishes a general method for constructing a notion of multi-faced independence from a given universal product of $\ast$-representations for a monoidal product. New examples of independence will be provided in Sections \ref{sec:tensor_lift} and \ref{sec:free_lift}. 
 
According to the axiomatic formulation in \cite{Muraki03} and \cite{BenGhorbalSchurmann}, independence can be interpreted as a product of states on the free product of $\ast$-algebras. The underlying intimate relation between monoidal categories with inclusions and independence relations has been first observed by Franz \cite{Franz06} and studied in detail by Gerhold, Lachs, and Sch\"{u}rmann \cite{GerholdLachsSchurmann22}.
Along the same lines, Manzel and Sch\"{u}rmann \cite{ManzelSchurmann17} presented an axiomatic approach to multi-faced independences based on universal products.  We will work with the definition of universal product of states on the \emph{free product of $m$-faced $\ast$-algebras} given by Gerhold in \cite[Remark 3.4]{Gerhold21p} (see Definition \ref{defi:universal_product} below).
In the case $m=1$, the universal products of states are classified in \cite{Muraki13} into the five kinds: tensor, free, boolean, monotone and antimonotone. 
The bi-free product \cite{Voiculescu14}, free-boolean product \cite{Liu19} and bi-monotone product \cite{Gerhold17p} are examples for $m=2$, and the free-free-boolean product in \cite{Liu_preprint} is an example for $m=3$. 
The product coming from the usual tensor independence for random vectors on $\mathbb C^m$ are also examples for any $m$. 
We will add more examples later in Subsections \ref{subsec:tensor_lift2} and \ref{subsec:free_lift2}.   

\begin{defi}\label{defi:universal_product}
    Let $m \in \N$. An \emph{$m$-faced $*$-algebra} is a $*$-algebra $A$ together with $\ast$-subalgebras $A^{(1)},\ldots, A^{(m)}$ freely generating $A$, which means that the canonical $*$-homomorphism $A^{(1)}\sqcup A^{(2)}\sqcup \cdots\sqcup A^{(m)}\to A$ is an isomorphism and which we indicate by writing $A=A^{(1)}\sqcup A^{(2)}\sqcup \cdots\sqcup A^{(m)}$. An \emph{$m$-faced $*$-homomorphism} is a $*$-homomorphism between $m$-faced $\ast$-algebras $j\colon A\to B$ with $j(A^{(k)})\subset B^{(k)}$. The free product $A_1\sqcup A_2$ of $m$-faced $*$-algebras is an $m$-faced $*$-algebra with $(A_1\sqcup A_2)^{(k)}=A_1^{(k)}\sqcup A_2^{(k)}$. 
    
    An \emph{$m$-faced universal product of states} is an operation $\odot$ which associates with two states $\phi_1,\phi_2$ on $m$-faced $\ast$-algebras $A_1,A_2$, respectively, a state $\phi_1\odot\phi_2$ on $A_1\sqcup A_2$ such that: 
    \begin{itemize}
        \item (Associativity) for any three states $\phi_{k}$ ($k\in\{1,2,3\})$, we have \[\phi_1\odot(\phi_2\odot\phi_3)=(\phi_1\odot\phi_2)\odot\phi_3;
        \]
        \item (Universality) for any two $m$-faced $*$-homomorphisms $j_{k}\colon B_{k}\rightarrow A_{k}$, we have 
        \[
        (\phi_1\circ j_1)\odot (\phi_2\circ j_2)=(\phi_1\odot\phi_2)\circ(j_1\sqcup j_2);
        \]
     \item (Restriction) for any two states $\phi_k\colon A_k \to \mathbb C$ $(k=1,2)$, we have \[
     \Restr{\phi_1\odot\phi_2}{A_k}=\phi_k
     \]
     for $k=1,2$.
    \end{itemize}
Moreover, an $m$-faced universal product $\odot$ is said to be \emph{symmetric} if for all states $\phi_{k}\colon A_k \to \mathbb C$ ($k\in\{1,2\})$ we have 
\[
\phi_1 \odot \phi_2 = (\phi_2 \odot \phi_1)\circ \Theta_{A_1,A_2}, 
\]
 where $\Theta_{A_1,A_2}$ is the natural identification mapping $A_1 \sqcup A_2 \to A_2 \sqcup A_1$. 
\end{defi}
\begin{rmq}
One can see that the restriction condition can be replaced by the following: 
\begin{itemize}
       \item (Unitality) let $\mathfrak{0}\colon\{0\}\to\{0\}$ be the trivial state. Then for any state $\phi$ on $A$ we have 
       \[
       \mathfrak{0}\odot\phi=\phi=\phi\odot\mathfrak{0} 
       \]
      under the natural identification $\{0\}\sqcup A = A =A\sqcup \{0\}$.
       \end{itemize}
For example, assuming the unitality condition, one retrieves the restriction property on $A_1$ by taking $j_1 = {\rm id}\colon A_1 \to A_1$ and $j_2 \colon \{0\}\to A_2$ and applying the universality.     
\end{rmq}

In Remark \ref{rmq:UP-linear-functionals} at the end of this section, we will compare Definition \ref{defi:universal_product} with the definitions of positive universal products used in \cite{ManzelSchurmann17} or \cite{Muraki13}.

The following definition shows how a notion of independence for non-commutative random vectors is related to a multi-faced universal product of states. 

\begin{defi}\label{defi:multi-faced_independence}
Let $m\in\N$. Let $\odot$ be an $m$-faced universal product of states and $(A,\phi)$ be a $\ast$-probability space. 

\begin{enumerate}[label=\rm(\arabic*)]
    \item 
Let $\mathbf A_k=(A_k^{(1)}, A_k^{(2)},\dots, A_k^{(m)})$ be an  $m$-tuple of $\ast$-subalgebras of $A$ for $k=1,2,\dots,n$. The sequence $(\mathbf A_1, \mathbf A_2, \dots, \mathbf A_n)$ is said to be \emph{$\odot$-independent} if for 
\begin{align*}
    A_k &\colonequals A_k^{(1)} \sqcup A_k^{(2)} \sqcup \cdots \sqcup A_k^{(m)} \quad \text{and}\\ 
    i_k&\colonequals i_k^{(1)}\sqcup i_k^{(2)}\sqcup \cdots \sqcup i_k^{(m)}\colon A_k \to A, 
\end{align*}
 where $i_k^{(\ell)}\colon A_k^{(\ell)} \hookrightarrow A$ is the embedding for $k=1,2,\dots, n$ and $\ell =1,2,\dots, m$, we have
 \[
 (\phi\circ i_1) \odot (\phi\circ i_2) \odot \cdots \odot (\phi\circ i_n) = \phi \circ (i_1 \sqcup i_2 \sqcup \cdots \sqcup i_n). 
 \]
 
\item Let $\mathbf a_k= (a_k^{(1)}, a_k^{(2)}, \dots, a_k^{(m)}), k=1,2,\dots, n,$ be $m$-tuples of elements of $A$ (also called random vectors). Then the sequence $(\mathbf a_1, \mathbf a_2, \dots, \mathbf a_n)$ is said to be \emph{$\odot$-independent} if the sequence $(\mathbf A_1, \mathbf A_2, \dots, \mathbf A_n)$ is $\odot$-independent, where 
\[
\mathbf A_k= \left(\Alg{a_k^{(1)}}, \Alg{a_k^{(2)}},\dots, \Alg{a_k^{(m)}}\right).
\]
 \end{enumerate}
The index set $\{1,2,\dots,n\}$ above can be generalized to an arbitrary linearly ordered set $I$ by applying the above definition(s) for each finite subset of $I$.      
\end{defi}

\begin{rmq}
If the $m$-faced universal product $\odot$ is symmetric, then one can see that $\odot$-independence of $(\mathbf A_1, \mathbf A_2, \dots, \mathbf A_n)$ implies that $(\mathbf A_{\sigma(1)}, \mathbf A_{\sigma(2)}, \dots, \mathbf A_{\sigma(n)})$ is also $\odot$-independent for any permutation $\sigma$. In this case the linear order on subalgebras does not matter and hence we can simply say that $\mathbf A_1, \mathbf A_2, \dots, \mathbf A_n$ are $\odot$-independent. Also, the index set $\{1,2,\dots,n\}$ can be generalized to any set. A similar remark applies to random vectors. 
\end{rmq}

Corresponding to the notions of multi-faced universal products of states, we extend the definition of universal products of representations to the multi-faced setting. 

\begin{defi}\label{defi:multifaced_universal_produict_repr} Let $m\in\N$ and $\boxtimes$ be a monoidal product with embeddings on $\PreHilb$. An \emph{$m$-faced universal product of $\ast$\hyp{}representations for $\boxtimes$} is a rule that, given two $*$\hyp{}representations of $m$-faced $*$-algebras\footnote{A $*$-representation of an $m$-faced $*$-algebra simply means a $*$-representation of the underlying $*$-algebra, ignoring the faces.} on pre-Hilbert spaces $\pi_k\colon A_k\to \Lop(H_k)$,
gives a $*$\hyp{}representation $\pi_1\odot\pi_2\colon A_1\sqcup A_2\to \Lop(H_1\boxtimes H_2)$ which fulfills associativity, universality of pre-Hilbert spaces,  the restriction property (as in the single-faced case --- see Definition \ref{defi:universal_product_representation}), and additionally:
    \begin{itemize}
        \item (Universality of $m$-faced $\ast$-algebras) for any $m$-faced $*$-homomorphisms $j_k\colon B_k\to A_k$ ($k\in\{1,2\}$), we have 
        \[
        (\pi_1\odot \pi_2)\circ (j_1\sqcup j_2)=(\pi_1\circ j_1)\odot (\pi_2\circ j_2).
        \]
    \end{itemize}
\end{defi}

Note that a $*$-representation $\pi\colon A\to \Lop(H)$ of an $m$-faced algebra $A$ can be identified with the collection $(\pi_i)_{i=1}^m$ of its restrictions to the faces $\pi^{(i)}:=\Restr{\pi}{A^{(i)}}$; indeed, $\pi$ can be reconstructed as $\pi=\pi^{(1)}\sqcup\cdots \sqcup \pi^{(m)}$. Furthermore, given an $m$-faced universal product of $\ast$\hyp{}representations $\odot$, one readily observes that $(\pi_1\odot \pi_2)^{(i)}$ only depends on $\pi_1^{(i)}$ and $\pi_2^{(i)}$. This allows, for each $i\in\{1,\ldots,m\}$, to write
$ (\pi_1\odot \pi_2)^{(i)} = \pi_1^{(i)}\odot_i \pi_2^{(i)} $ for uniquely determined single-faced universal products of $*$-representations $\odot_i$ (which follows easily from the universal property of the free product, we leave the details to the reader). Conversely, every $m$-tuple of single-faced universal products of $*$-representations allows to define an $m$-faced universal product of $*$-representations $ \pi_1\odot \pi_2 := (\pi_1^{(1)}\odot_1 \pi_2^{(1)})\sqcup\cdots\sqcup (\pi_1^{(m)}\odot_m \pi_2^{(m)}) $. In the following we identify $\odot$ with the collection $(\odot_i)_{i=1}^m$. With this convention, we obtain the following generalization of Theorem \ref{thm:correspondence} for free.


\begin{cor}\label{cor:multi-correspondence} Let $\boxtimes$ be a monoidal product with embeddings on $\PreHilb$. The one-to-one correspondence between single-faced universal lifts and universal products of $\ast$\hyp{}representations for $\boxtimes$ can be canonically extended to obtain one-to-one correspondences between
\begin{itemize}
    \item $m$-tuples of universal lifts for $\boxtimes$
    \item $m$-faced universal products of $\ast$\hyp{}representations for $\boxtimes$
\end{itemize}
given by
    \begin{align*}
        \odot=(\odot_i)_{i=1}^m\mapsto (\lambda^{(i)},\rho^{(i)})_{i=1}^m\quad &\text{with}\quad \begin{cases}\lambda^{(i)}_{H_1,H_2}\colonequals 
        \Restr{\id_{\Lop(H_1)}\odot_i \id_{\Lop(H_2)}}{\Lop(H_1)},\\ \rho^{(i)}_{H_1,H_2}\colonequals \Restr{\id_{\Lop(H_1)}\odot_i \id_{\Lop(H_2)}}{\Lop(H_2)}
        \end{cases}
        \end{align*}
for pre-Hilbert spaces $H_1, H_2$, and 
\begin{align*}
        (\lambda^{(i)}, \rho^{(i)})_{i=1}^m\mapsto(\odot_i)_{i=1}^m=\odot \quad &\text{with}\quad \pi_1\odot_i\pi_2\colonequals (\lambda^{(i)}_{H_1,H_2}\circ\pi_1)\sqcup(\rho^{(i)}_{H_1,H_2}\circ\pi_2)
    \end{align*}
    for (single-faced) $\ast$\hyp{}representations $\pi_k\colon A_k\to \Lop(H_k)$, $k=1,2$.
\end{cor}

Consequently, a classification of $m$-faced universal products of $\ast$\hyp{}representations is equivalent to a classification of single-faced universal products of $\ast$\hyp{}representations, or of universal lifts. Such a simple reduction to the single-faced case is impossible for $m$-faced universal products of states. There is however a close connection between universal products of $\ast$\hyp{}representations and of states. The following theorem states that every $m$-faced universal product of $\ast$\hyp{}representations gives rise to an $m$-faced universal product of states. It is an open question at the moment whether every  $m$-faced universal product of states can be acquired in that fashion. One problem might be that the definition of $m$-faced universal products of $\ast$\hyp{}representations depends on a chosen monoidal product on $\PreHilb$, and it is not clear whether other choices besides the tensor and the free product are necessary to succeed.

\begin{lm}\label{lm:repr_state} Let $m\in \N$ and $\boxtimes$ be a monoidal product with embeddings on $\PreHilb$. Let $\odot$ be an $m$-faced universal product of $\ast$\hyp{}representations for $\boxtimes$. 
For $k=1,2$, let $\pi_k\colon A_k \to \Lop(H_k)$ and $\sigma_k\colon A_k \to \Lop(G_k)$ be $\ast$\hyp{}representations of $m$-faced $\ast$-algebras $A_k$ such that $\langle \Omega, \pi_k(\cdot) \Omega\rangle = \langle \Omega, \sigma_k(\cdot) \Omega\rangle$ on $A_k$. Then 
\[\langle \Omega, (\pi_1\odot \pi_2)(\cdot)\Omega\rangle = \langle \Omega, (\sigma_1\odot \sigma_2)(\cdot)\Omega\rangle \quad \text{on}\quad A_1 \sqcup A_2.\]
\end{lm}

\begin{proof} Let $\tilde H_k \colonequals\mathbb C \Omega + \pi_k (A_k)\Omega \subset H_k$ and $\tilde \pi_k\colon A_k \to \Lop(\tilde H_k)$ be the restriction of $\pi_k$ to its invariant subspace $\tilde H_k$.  
Similarly, we define $\tilde \sigma_k\colon A_k \to \Lop(\tilde G_k)$. 

By the universality of pre-Hilbert spaces, we have for all $a \in A_1 \sqcup A_2$ 
\begin{align}
\tilde \pi_1\odot \tilde \pi_2(a) &=\Restr{\pi_1\odot \pi_2(a)}{\tilde H_1\boxtimes \tilde H_2} \quad \text{and} \label{tilde_pi} \\ 
\tilde \sigma_1\odot \tilde \sigma_2(a) &=\Restr{\sigma_1\odot \sigma_2(a)}{\tilde G_1\boxtimes \tilde G_2}. \label{tilde_sigma} 
\end{align}
To see this, let $V_k\colon \tilde H_k \hookrightarrow H_k$ be the embedding (note that this is an $\Omega$-preserving isometry but may not be adjointable --- see Remark \ref{rmq:nonadjointabe}). Because $V_k \tilde \pi_k(a) = \pi_k(a)V_k$ for all $a \in  A_k$ and $k\in\{1,2\}$, the universality of pre-Hilbert spaces implies 
\begin{equation}\label{tilde_pi2}
(V_1\boxtimes V_2) [\tilde \pi_1\odot \tilde \pi_2(a)] = [\pi_1\odot \pi_2(a)] (V_1\boxtimes V_2)
\end{equation}
for all $a \in  A_1 \sqcup  A_2$, i.e.\ \eqref{tilde_pi}. By the obvious symmetry we have \eqref{tilde_sigma} too. 
 
By the preceding arguments, it suffices to prove that 
\[
\langle \Omega, (\tilde\pi_1\odot \tilde\pi_2)(a)\Omega\rangle = \langle \Omega, (\tilde\sigma_1\odot \tilde\sigma_2)(a)\Omega\rangle
\]
for all $a \in  A_1 \sqcup  A_2$. 

Let $\UU_k\colon \tilde G_k \to \tilde H_k$ be defined by $\UU_k[\alpha \Omega +  \sigma_k(a)\Omega]= \alpha \Omega + \pi_k(a)\Omega$ for $\alpha\in\mathbb C$ and $a \in A_k$. This is a well-defined isometry (actually, unitary) because  
\begin{align*}
\|\alpha \Omega + \sigma_k(a)\Omega\|^2
&= |\alpha|^2 +  \alpha \langle \Omega,  \sigma_k(a^*) \Omega\rangle + \overline{\alpha} \langle \Omega, \sigma_k(a) \Omega\rangle + \langle \Omega, \sigma_k(a^*a) \Omega\rangle\\
&= |\alpha|^2 +  \alpha \langle \Omega,  \pi_k(a^*) \Omega\rangle + \overline{\alpha} \langle \Omega, \pi_k(a) \Omega\rangle + \langle \Omega, \pi_k(a^*a) \Omega\rangle\\
&= \|\alpha \Omega + \pi_k(a)\Omega\|^2.
\end{align*}
The computation 
\begin{align*}
\UU_k \tilde\sigma_k(a)[\alpha\Omega+\sigma_k(b)\Omega] 
&= \UU_k[\alpha \sigma_k(a)\Omega+\sigma_k(a b)\Omega] \\
&= \alpha \pi_k(a)\Omega+\pi_k(a b)\Omega \\
&= \tilde\pi_k(a)[\alpha\Omega+\pi_k(b)\Omega] \\
&= \tilde\pi_k(a)\UU_k[\alpha\Omega+\sigma_k(b)\Omega]
\end{align*}
yields that $\UU_k\tilde\sigma_k(a) = \tilde\pi_k(a) \UU_k$ for all $a \in A_k$, that is, $\UU_k$ is an ($\Omega$-preserving) intertwiner. 

For $a\in A_1\sqcup A_2$ we have 
\begin{align*}
\langle \Omega,(\tilde\pi_1\odot\tilde\pi_2)(a)\Omega\rangle
&= \langle \Omega,(\tilde\pi_1\odot\tilde\pi_2)(a)(\UU_1 \boxtimes \UU_2)\Omega\rangle \\
&=\langle(\UU_1\boxtimes \UU_2)\Omega,(\UU_1\boxtimes \UU_2)(\tilde\sigma_1\odot\tilde\sigma_2)(a)\Omega\rangle\\
&=\langle\Omega,(\tilde\sigma_1\odot\tilde\sigma_2)(a)\Omega\rangle, 
\end{align*}
where the last two lines are due to the universality of pre-Hilbert spaces of $\odot$ and the fact that $\UU_1\boxtimes \UU_2$ is an isometry which fixes $\Omega$. 
\end{proof}

\begin{rmq}\label{rmq:nonadjointabe} The isometries $V_k$ in Lemma \ref{lm:repr_state} are not adjointable in general. To see this, dropping the dependence on $k \in \{1,2\}$ for notational brevity, we take $A=C[0,1]$ regarded as a 1-faced $\ast$-algebra and $H=L^2[0,1]$ equipped with the standard inner product $\langle f, g\rangle= \int_0^1 \overline{f(x)} g(x)\,dx$ and the unit vector $\Omega\equiv 1$. Consider the $\ast$-representation $\pi\colon A \to \Lop(H)$, where $\pi(a)$ acts as the multiplication operator by $a \in A$. Then we have $\tilde H = \mathbb C \Omega+ \pi(A)\Omega = C[0,1]\subset H$. One can see that the embedding $V\colon \tilde H \hookrightarrow H$  is not adjointable. This is because if it had an adjoint $V^* \colon H \to \tilde H$ then $V^*V=\id_{\tilde H}$ and therefore $P\colonequals  VV^*$ would be a projection ($P=P^2=P^\ast$). Then for any element $x \in H$ the element $(\id_{H}-P)x$ would be orthogonal to $P(H)$. However, $P(H)=V( \tilde H) \subset H$ is dense in $H$, so that $(\id_H-P)x=0$, so $x = P x \in \tilde H$, a contradiction if we choose $x \in H \setminus \tilde H$. 
\end{rmq}

\begin{thm}\label{thm:repr_state}
    Let $m\in\N$ and $\boxtimes$ be a monoidal product with embeddings on $\PreHilb$. Let $\odot$ be an $m$-faced universal product of  $\ast$\hyp{}representations for $\boxtimes$. Then there exists an $m$-faced universal product of states, also denoted $\odot$, such that $\phi_k(\cdot)=\langle\Omega,\pi_k(\cdot)\Omega\rangle$ on $A_k$ for all $k\in\{1,2\}$ implies
    \[\phi_1\odot \phi_2 (\cdot)= \langle \Omega, (\pi_1\odot \pi_2)(\cdot)\Omega\rangle \quad \text{on}\quad A_1\sqcup A_2. \]
\end{thm}

\begin{proof}  For $k=1,2$, let $A_k$ be an $m$-faced $\ast$-algebra and $\phi_k\colon A_k\to \mathbb C$  be a state. Take a $\ast$\hyp{}representation $\pi_k\colon A_k \to \Lop(H_k)$ that realizes $\phi_k$, i.e.\  
$\phi_k(\cdot) = \langle \Omega, \pi_k(\cdot)\Omega \rangle$.  Note that, thanks to the assumption of $\phi_k$ being a state, such a $\ast$\hyp{}representation exists --- see Remark \ref{rmq:restricted_state}.
We then define a linear functional $\phi_1 \odot \phi_2$ on $A_1 \sqcup A_2$ by 
\begin{equation}\label{eq:def_universal_product}
\phi_1\odot \phi_2 (a)= \langle \Omega, (\pi_1\odot \pi_2)(a)\Omega\rangle
\end{equation}
    for all $a\in A_1\sqcup A_2$. 
This is a state --- see again Remark \ref{rmq:restricted_state}. By Lemma \ref{lm:repr_state} the definition of $\phi_1\odot \phi_2$ does not depend on a choice of $\pi_1$ and $\pi_2$. By definition, the last statement of the theorem clearly holds. It then remains to check the axioms for $\odot$ required in Definition \ref{defi:universal_product}.

\subcase{Universality} Let $j_{k}\colon B_{k}\rightarrow A_{k}$ be two  $m$-faced $*$-homomorphisms for $k=1,2$. Then, $\pi_k\circ j_k$ is a $*$\hyp{}representation of $B_k$ such that $\phi_k\circ j_k=\langle\Omega,(\pi_k\circ j_k)(\cdot)\Omega\rangle$. Therefore:
\begin{align*}
(\phi_1\odot\phi_2)\circ(j_1\sqcup j_2)(\cdot)&= \langle\Omega,(\pi_1\odot\pi_2)\circ(j_1\sqcup j_2)(\cdot)\Omega\rangle\\
&= \langle\Omega,(\pi_1\circ j_1)\odot(\pi_2\circ j_2)(\cdot)\Omega\rangle\\
&= (\phi_1\circ j_1)\odot(\phi_2\circ j_2)(\cdot). 
\end{align*}
This concludes the universality of the proposed universal product.

\subcase{Associativity} Let $\phi_{k}\colon A_{k}\rightarrow\mathbb{C}$ ($k \in \{1,2,3\}$) be three states and take $\ast$\hyp{}representations $\pi_{k}$ which realize $\phi_k$. We have, by the  definition \eqref{eq:def_universal_product}, 
\[
\phi_1\odot\phi_2(\cdot)=\langle\Omega,(\pi_1\odot\pi_2)(\cdot)\Omega\rangle  \qquad \text{and}\qquad \phi_2\odot\phi_3(\cdot)=\langle\Omega,(\pi_2\odot\pi_3)(\cdot)\Omega\rangle,
\]
and hence $\pi_1\odot\pi_2$ and $\pi_2\odot\pi_3$ are $*$\hyp{}representations which realize $\phi_1\odot\phi_2$ and $\phi_2\odot\phi_3$, respectively. We can then proceed as 
\begin{align*}
(\phi_1\odot\phi_2)\odot\phi_3 (\cdot)&= \langle\Omega,(\pi_1\odot\pi_2)\odot\pi_3(\cdot)\Omega\rangle\\
&= \langle\Omega,\pi_1\odot(\pi_2\odot\pi_3)(\cdot)\Omega\rangle\\
&= \phi_1\odot(\phi_2\odot\phi_3) (\cdot),  
\end{align*}
concluding the associativity.

\subcase{Restriction} For $k\in \{1,2\}$, let $\phi_k\colon A_k \to \mathbb C$ be a state and $\pi_k\colon A_k\rightarrow \Lop(H_k)$ be a $\ast$\hyp{}representation which realizes $\phi_k$. Then, for $a\in A_k$ we have 
\begin{align*}
\phi_1\odot\phi_2(a)&= \langle\Omega,(\pi_1\odot\pi_2)(a)\Omega\rangle\\
&= \langle\Omega, \pi_k(a)\Omega\rangle\\
&=\phi_k(a), 
\end{align*}
so that the restriction property holds true. 
\end{proof}

\begin{rmq}\label{rmq:UP-linear-functionals}
    A priori it is not clear whether \emph{universal products of states} as defined in Definition \ref{defi:universal_product} are the same as the positive universal products of linear functionals considered for example in \cite{ManzelSchurmann17}. The latter are defined for arbitrary linear functionals on algebras, the positivity condition states that the product of states on $*$-algebras is again a state on the free product $*$-algebra, i.e.\ a positive universal product of linear functionals ``restricts'' to a universal product of states. However, when we start with a universal product of states, in order to identify it with a positive universal product of linear functionals, we would have to extend it to arbitrary linear functionals. We do not know at the moment whether there is a direct way to do this for arbitrary $m$-faced universal products of states. However, we want to at least indicate how one could check that the examples of universal products of states discussed in this paper extend to a universal product of linear functionals. The representation theoretic approach in principle also works without positivity. 
    Every linear functional on an algebra can be realized as $\varphi=P_\Omega \pi(\cdot)\Omega$ (identifying $\lambda\Omega\in\mathbb C\Omega$ with $\lambda\in \mathbb C$) for a representation $\pi$ on a vector space $\mathbb C\Omega \oplus \hat V$. (There are some possible choices of what extra structure and conditions one wants to impose on $V$ and $\pi$, see e.g.\ \cite[Section 5]{GerholdLachs15}.) The free product and the tensor product are also defined for vector spaces with such a decomposition. For the concrete lifts we will exhibit in Sections \ref{sec:tensor_lift} and \ref{sec:free_lift}, it is not difficult to check that one can extend them and the corresponding universal products of representations to the non-*-case. It remains to check the analogue of Lemma \ref{lm:repr_state}, that 
    \[\varphi_1\odot \varphi_2(a)\colonequals P_\Omega \pi_1\odot\pi_2(a)\Omega\]
    where $\varphi_k=P_\Omega \pi_k(\cdot)\Omega$ is well-defined, i.e.\ the right hand side only depends on the $\varphi_k$, not on the chosen representations. Our concrete calculations of mixed moments in Subsection \ref{subsec:free_lift2} suggest that this would not be too hard to prove. However, since the paper is already quite long, and since it would be much more interesting to have a general argument, we postpone a detailed discussion of these issues.
\end{rmq}

\section{Universal lifts to the tensor product}\label{sec:tensor_lift}

In Subsection \ref{subsec:tensor_lift1} we classify  universal lifts to the tensor product $H_1\otimes H_2=\mathbb C \Omega \oplus [\hat H_1\oplus\hat H_2 \oplus (\hat H_1\otimes \hat H_2)]$, with the shorthand notation $\hat H_k\colonequals H_k\ominus \mathbb C\Omega$. When appropriate, we denote the unit vector in $H_1\otimes H_2$ as  $\Omega\otimes\Omega$ instead of $\Omega$. 
It is remarkable that the proofs do not need the left and right associativity axioms (cf.\ Remark \ref{rem_lifts}). In Subsection \ref{subsec:tensor_lift2} we apply Corollary~\ref{cor:multi-correspondence} and Theorem~\ref{thm:repr_state} to
the classified universal lifts to construct some new universal multi-faced products of states. We also classify the symmetric products among them.

\subsection{Classification of universal lifts}\label{subsec:tensor_lift1}

Let $\lambda$ be a left universal lift to the tensor product. 
We often use the following weaker form of universality  of $\lambda$: for any $\TT$ in $\Lop(H_1)$ and any $\Omega$-preserving \emph{adjointable} isometries $\UU_k\colon H_ k\to G_k$ $(k \in \{1,2\})$, we have 
\begin{align*}
\lambda_{G_1,G_2}(\UU_1 \TT \UU_1^*)\circ(\UU_1\otimes \UU_2)&= (\UU_1\otimes \UU_2)\circ\lambda_{H_1,H_2}(\TT).
\end{align*}

We will often use the notation $P_E \in \Lop(H)$ for the orthogonal projection onto a subspace $E$ of a pre-Hilbert space $H$ when it exists (typically, when $E$ is of finite dimension). If $E= \mathbb C \Omega$ then $P_{\mathbb C \Omega}$ will be abbreviated to $P_\Omega$.

For $x\in\hat H$ consider the operator
\begin{align*}
a_x^*\colon H\to H; \quad\Omega\mapsto x, \quad  \Restr{a_x^*}{\hat H}=0,
\end{align*}
whose adjoint is given by 
\begin{align*}
 a_x\colon x' \mapsto \langle x,x^\prime\rangle\Omega \quad \text{for all}\quad x' \in \hat H, \quad \Omega\mapsto0. 
\end{align*}

\begin{lm}\label{lm:unique}
Suppose that $\lambda$ and $\lambda'$ are left universal lifts for $\otimes$ such that $\lambda_{H_1,H_2}(a_x^*)=\lambda'_{H_1,H_2}(a_x^*)$ for all pre-Hilbert spaces $H_1$ and $H_2$  
and all $x \in \hat H_1$. Then $\lambda=\lambda'$. 
\end{lm}
\begin{proof}
We fix pre-Hilbert spaces $H_1,H_2$ and abbreviate $\lambda_{H_1,H_2}$ to $\lambda$ and similarly for $\lambda'$. When $\hat H_1\neq\{0\}$, all adjointable finite-rank operators are in the $*$-algebra generated by $\{a^*_x: x\in\hat H_1\}$, and therefore $\lambda_{H_1,H_2}(F) = \lambda_{H_1,H_2}'(F)$ for all finite-rank operators $F \in \Lop(H_1)$. 
For $T \in \Lop(H_1), x \in \hat H_1, y\in \hat H_2$, we get 
\begin{align}
\lambda(T)x\otimes y&=\lambda(TP_{\mathbb C\Omega+\mathbb Cx} + TP_{(\mathbb C\Omega+\mathbb Cx)^\perp})x\otimes y \label{eq:finite_rank} \\
&=\lambda(TP_{\mathbb C\Omega+\mathbb Cx})x\otimes y + \lambda(T)\lambda(P_{(\mathbb C\Omega+\mathbb Cx)^\perp})x\otimes y. \notag
\end{align}
The embedding $\UU_1\colon \mathbb C \Omega + \mathbb C x\hookrightarrow H_1$ satisfies $P_{(\mathbb C\Omega+\mathbb Cx)^\perp} \UU_1 = \UU_1 0_{\mathbb C \Omega + \mathbb C x\to \mathbb C \Omega + \mathbb C x}$; therefore the universality yields 
\begin{align*}
\lambda(P_{(\mathbb C\Omega+\mathbb Cx)^\perp})x\otimes y 
&= \lambda(P_{(\mathbb C\Omega+\mathbb Cx)^\perp}) (\UU_1\otimes \id)x\otimes y \\
&= (\UU_1\otimes \id)\lambda_{\mathbb C \Omega + \mathbb C x,H_2}(0_{\mathbb C \Omega + \mathbb C x\to \mathbb C \Omega + \mathbb C x}) x\otimes y=0.
\end{align*}
Combined with \eqref{eq:finite_rank} this implies 
$\lambda(T)x\otimes y= \lambda(TP_{\mathbb C\Omega+\mathbb Cx})x\otimes y$. The same reasoning gives $\lambda'(T)x\otimes y= \lambda'(TP_{\mathbb C\Omega+\mathbb Cx})x\otimes y$. Since $TP_{\mathbb C\Omega+\mathbb Cx}$ is a finite-rank operator, we have
$
\lambda(TP_{\mathbb C\Omega+\mathbb Cx})=\lambda'(TP_{\mathbb C\Omega+\mathbb Cx}).
$
Combining those facts we conclude that $\lambda(T)=\lambda'(T)$.

When $H_1= \mathbb C\Omega$, we fix any pre-Hilbert space $G_1$ of dimension $\ge2$ and the canonical embedding $\UU\colon \mathbb{C}\Omega\hookrightarrow G_1$, and then use the universality condition to get  
\[
\lambda(\id_{\mathbb{C}\Omega})=(\UU^*\otimes\id)\lambda_{G_1,H_2}(P_\Omega)(\UU\otimes\id).
\]
Since $P_\Omega$ is a finite-rank operator on $G_1$, the previous arguments imply $\lambda_{G_1,H_2}(P_\Omega)=\lambda_{G_1,H_2}'(P_\Omega)$, and hence $\lambda(\id_{\mathbb{C}\Omega})= \lambda'(\id_{\mathbb{C}\Omega})$. 
\end{proof}

In the following, let $\mathbb T$ denote the unit circle as a subset of the complex plane. 
\begin{lm}\label{lm:gamma}
    Let $\lambda$ be a left universal lift for $\otimes$. Then there exists a constant $\gamma\in\mathbb T \cup\{0\}$ such that for all $H_1,H_2$ 
    and all $x \in \hat H_1$
  \[\lambda_{H_1,H_2}(a_x^*)=a^*_x\otimes P_\Omega + \gamma a^*_x\otimes P_{\Omega^\perp}.\]
\end{lm}

\begin{proof} We fix pre-Hilbert spaces $H_1,H_2$ and abbreviate $\lambda_{H_1,H_2}$ to $\lambda$.
    To prove the statement, we evaluate $\lambda(a^*_x)$ on vectors of the form $\Omega\otimes \Omega$, $x'\otimes \Omega$, $\Omega\otimes y$ and $x'\otimes y$ with $x'\in \hat H_1,y\in \hat H_2$.

Using the restriction axiom, one readily obtains $\lambda(a_x^*)\Omega\otimes\Omega=x\otimes \Omega$ and  $\lambda(a_x^*)x^\prime\otimes\Omega=(a_x^*x^\prime)\otimes\Omega=0$, i.e.\ $\lambda(a^*_x)(\id\otimes P_\Omega)=a_x^*\otimes P_\Omega$. When $\hat H_2 =\{0\}$ this implies that $\lambda_{H_1,H_2}(a_x^*)=a_x^* \otimes P_\Omega$ and we are done. Hence, we may assume that $\hat H_2 \ne\{0\}$ below.

Next, universality implies that $\lambda(a_x^*)\Omega \otimes y\in\Span(\Omega\otimes\Omega, \Omega\otimes y, x\otimes\Omega, x\otimes y)$; indeed, with the inclusion maps $\UU_1\colon \Span(\Omega,x)\hookrightarrow H_1$ and $\UU_2\colon \Span(\Omega,y)\hookrightarrow H_2$, one sees that $a^*_x=\UU_1 \Restr{a^*_x}{\Span(\Omega,x)} \UU_1^*$  and, therefore,
\begin{align*}
\lambda(a^*_x)\Omega\otimes y
&=\lambda(\UU_1 \Restr{a^*_x}{\Span(\Omega,x)} \UU_1^*)(\UU_1\otimes \UU_2)(\Omega\otimes y) \\
&=(\UU_1\otimes \UU_2) \lambda_{\Span(\Omega,x), \Span(\Omega,y)}(\Restr{a^*_x}{\Span(\Omega,x)}) \Omega\otimes y \\
&\quad \in \Span(\Omega,x)\otimes \Span(\Omega,y).
\end{align*} Also, universality implies that the coefficients in 
\begin{equation}\label{eq:universal_coeff}
    \lambda(a^*_x)\Omega\otimes y = \alpha \Omega\otimes \Omega + \beta x\otimes \Omega + \delta \Omega\otimes y + \gamma x\otimes y
\end{equation}
cannot depend on unit vectors $x$ and $y$ and pre-Hilbert spaces $H_1, H_2$ with $\hat H_1 \ne\{0\}, \hat H_2 \ne \{0\}$; to see this, fix a reference pre-Hilbert space $E=\mathbb C \Omega\oplus \mathbb C e$ with a unit vector $e$, define two isometries which map $e$ to $x \in \hat H_1$ and $e$ to $y \in \hat H_2$, respectively, and apply the universality. 
Furthermore, if we replace $x$ and $y$ with the unit vectors $e^{i \theta}x$ and $e^{i \psi}y$ with $\theta,\psi \in\R$ respectively, then \eqref{eq:universal_coeff} is strengthened to 
\[
e^{i\theta}e^{i\psi}\lambda(a_x^*) \Omega \otimes y = \alpha  \Omega\otimes \Omega + \beta e^{i\theta} x\otimes \Omega + \delta  e^{i\psi} \Omega\otimes y + \gamma e^{i\theta} e^{i\psi} x\otimes y, \qquad \theta,\psi \in\R.  
\]
By the uniqueness of Fourier series expansion in two variables, we get $\alpha=\beta=\delta=0$, so that $\lambda(a_x^*)\Omega \otimes y=\gamma x\otimes y$ for all unit vectors $x$ and $y$. It is easy to see by linearity that this holds for all (not necessarily unit) vectors $x$ and $y$. 

Applying analogous reasoning to $\lambda(a_x^*)x^\prime\otimes y$ with the inclusion maps $\UU_1\colon \Span(\Omega,x,x')\hookrightarrow H_1$ and $\UU_2\colon \Span(\Omega,y)\hookrightarrow H_2$, we find $\lambda(a_x^*)x^\prime\otimes y\in\Span(\space{\Omega,x,x^\prime})\otimes \Span(\space{\Omega,y})$. 
The coefficients must be independent of $x,x^\prime,y$ in the unit sphere as long as $x \perp x'$. 
Because the result must be linear in $x,x^\prime,y$, this implies $\lambda(a_x^*)x^\prime\otimes y=0$ for all $x,x',y$ with $x \perp x'$. 
Similar arguments work for the case $x=x'$ and we arrive at the same conclusion $\lambda(a_x^*)x\otimes y=0$. 
Those two cases are combined to yield $\lambda(a_x^*)x^\prime\otimes y=0$ for all $x, x' \in \hat H_1$ and $y \in \hat H_2$.

The last two paragraphs show $\lambda(a^*_x)(\id\otimes P_{\Omega^\perp})=\gamma a_x^*\otimes P_{\Omega^\perp}$.
Combining everything, we have proved
\[
\lambda(a_x^*) = \lambda(a^*_x)(\id\otimes P_{\Omega}) + \lambda(a^*_x)(\id\otimes P_{\Omega^\perp}) = a^*_x\otimes P_\Omega + \gamma a^*_x\otimes P_{\Omega^\perp}.
\]

It remains to show that $\gamma \in \mathbb T \cup \{0\}$. For any unit vectors $x\in\hat H_1$ and $y \in \hat H_2$, we have 
\begin{align*}
   \lambda(a_xa_x^*)\Omega\otimes y 
   &=\lambda(a_x)\lambda(a_x^*)\Omega\otimes y \\
   &= (a_x\otimes P_\Omega + \overline{\gamma} a_x\otimes P_{\Omega^\perp})\gamma x \otimes y \\
   &= |\gamma|^2\Omega\otimes y, 
 \end{align*} 
 and hence $|\gamma|^2$ is an eigenvalue of $ \lambda(a_xa_x^*)$. Because $a_xa_x^*=P_\Omega$ is a projection, $\lambda(a_xa_x^*)$ is a projection as well, and consequently its eigenvalue $|\gamma|^2$ is 0 or 1, as claimed.
\end{proof}

In order to extend the formula for $\lambda_{H_1,H_2}(a_x^*)$ in Lemma \ref{lm:gamma} to $\lambda_{H_1,H_2}(X)$ for general operators $X \in \Lop(H_1)$, we introduce the linear mapping
\[
\Lop(H,G) \to \Lop(H,G), \qquad T\mapsto T_\gamma = \begin{pmatrix} |\gamma|\tau & (\gamma t')^* \\ \gamma t & |\gamma|\hat T \end{pmatrix}
\]
defined through the matrix form $T = \begin{pmatrix} \tau & (t')^* \\ t & \hat T \end{pmatrix}$ according to the decomposition $H = \mathbb C\Omega \oplus \hat H$ and $G = \mathbb C\Omega \oplus \hat G$, that is, $\tau \in \mathbb C, t \in \hat G, t' \in \hat H,$ and $\hat T \in \Lop(\hat H, \hat G)$.  This is a minor generalization of \eqref{eq:hom_gamma} to rectangular operators, which will be sometimes useful below. \revise{Note in particular that for $x\in\hat H$ the matrix form of $a_x^*$ is given by $a_x^*=
\begin{pmatrix}
  0&0\\x&0
\end{pmatrix}
$ and, thus, $\left(a_x^{*}\right)_\gamma=\gamma a_x^{*}$.}

\begin{lm}\label{lm:hom} Let $\gamma\in \mathbb T \cup\{0\}$. 

\begin{enumerate}[label=\rm(\arabic*)]

\item \label{item:hom1}
For $S \in \Lop(H, G)$ and $T \in \Lop(K, H)$ we have 
\begin{equation*}
(ST)_\gamma = S_\gamma T_\gamma\qquad \text{and} \qquad (T^\ast)_\gamma = (T_\gamma)^\ast. 
\end{equation*}

\item \label{item:hom2}
Let $P_\Omega \in \Lop(H_2)$ and $P_{\Omega^\perp} \in \Lop(H_2)$ be the orthogonal projections onto $\mathbb C \Omega$ and $\hat H_2$, respectively.  For $X \in \Lop(H_1)$ we have 
\[ 
(X \otimes P_\Omega)_\gamma = X_\gamma \otimes P_\Omega  \qquad \text{and}\qquad (X \otimes P_{\Omega^\perp})_\gamma = |\gamma|X \otimes P_{\Omega^\perp}. 
\]
\end{enumerate}
\end{lm}
\begin{proof}
\subcasenps{\ref{item:hom1}} The proof is simple computations together with the properties $|\gamma|^2=|\gamma|$ and $\gamma|\gamma| = \gamma$. 

\subcasenp{\ref{item:hom2}}
According to the decomposition $H_1 \otimes H_2 = \mathbb C \Omega \oplus [(\hat H_1 \otimes \Omega) \oplus (\Omega \otimes \hat H_2) \oplus (\hat H_1 \otimes \hat H_2)]$, the operator $X \otimes P_\Omega$ has the matrix form 
\[
X \otimes P_{\Omega} = 
\begin{pmatrix} 
\xi & (x' \otimes \Omega)^* \\ x\otimes \Omega & \hat X \otimes P_\Omega 
\end{pmatrix}
\]
where $X = \begin{pmatrix} \xi & (x')^* \\ x & \hat X \end{pmatrix}$. Therefore, we obtain  
\[ 
(X \otimes P_\Omega)_\gamma 
= 
\begin{pmatrix} 
|\gamma|\xi & (\gamma x' \otimes \Omega)^* \\ \gamma x\otimes \Omega & |\gamma|\hat X \otimes P_\Omega 
\end{pmatrix}
 = X_\gamma \otimes P_\Omega.
\]
On the other hand, the operator $X \otimes P_{\Omega^\perp}$ has the matrix form 
\[
X \otimes P_{\Omega^\perp} = 
\begin{pmatrix} 
0 & 0 \\ 0 & X'
\end{pmatrix}
\]
for some $X'$, and hence

\[ 
(X \otimes P_{\Omega^\perp})_\gamma 
= 
\begin{pmatrix} 
0 & 0 \\ 0 & |\gamma| X'
\end{pmatrix}
 = |\gamma|(X \otimes P_{\Omega^\perp})
\]
as desired.
\end{proof}

We are ready to state the classification result. \pagebreak[3]

\begin{thm}\label{thm:tensor_lift} \
\begin{enumerate}[label=\rm(\arabic*)]
    \item\label{item1}  The left universal lifts and right universal lifts for $\otimes$ are respectively classified into the one-parameter families $\{\lambda^\gamma\}_{\gamma\in \mathbb T \cup\{0\}}$ and $\{\rho^\gamma\}_{\gamma\in \mathbb T \cup\{0\}}$ defined by  
  \begin{align*}\label{universal_lift}
      \lambda_{H_1,H_2}^\gamma (X)=X\otimes P_{\Omega} + X_\gamma \otimes P_{\Omega^\perp} \qquad \text{and} \qquad \rho_{H_1,H_2}^\gamma(Y) = P_\Omega \otimes Y + P_{\Omega^\perp} \otimes Y_\gamma
      \end{align*}   
for $X \in \Lop(H_1)$ and $Y\in \Lop(H_2)$.

  \item\label{item3} The universal lifts for $\otimes$ are classified into the family  $\{(\lambda^{\gamma},\rho^{\delta})\}_{(\gamma, \delta) \in J_\otimes}$, where $J_\otimes=\{(\gamma,\delta) \in (\mathbb T \cup\{0\})^2: \gamma=\delta ~\text{or}~ \gamma=0 ~\text{or}~\delta=0\}$.
  
\end{enumerate}
\end{thm}

\begin{proof} \subcasenps{\ref{item1}} It suffices to work on the left lifts by symmetry. Since 
$
a^*_x=\begin{pmatrix} 0 & 0\\ x & 0 \end{pmatrix}\in \Lop(H_1)$, 
we find
  \[\lambda_{H_1,H_2}^\gamma(a^*_x)=\lambda_{H_1,H_2}^\gamma 
    \begin{pmatrix}
        0 & 0\\ x & 0
    \end{pmatrix}
    = a^*_x\otimes P_\Omega + \gamma a^*_x \otimes P_{\Omega^\perp}.\] 
According to Lemma \ref{lm:unique}, a left universal lift is uniquely determined by its action on operators of the form $a^*_x$. This implies, combined with Lemma \ref{lm:gamma}, that any left universal lift must coincide with a certain $\lambda^\gamma$. 
    
    It remains to prove that $\lambda^\gamma$ is really a left universal lift for all $\gamma \in \mathbb T\cup\{0\}$. To begin, $\lambda_{H_1,H_2}^\gamma$ is a $*$-homomorphism because $X\mapsto X_\gamma$ is a  $*$-homomorphism on $\Lop(H_1)$ by Lemma \ref{lm:hom} \ref{item:hom1}.   
For the left associativity, the application of Lemma \ref{lm:hom} \ref{item:hom2} implies
  \begin{align*}
      &\lambda^\gamma_{H_1\otimes H_2,H_3}(\lambda^\gamma_{H_1,H_2} (X))
      =\lambda^\gamma_{H_1\otimes H_2,H_3} (X\otimes P_\Omega + X_\gamma \otimes P_{\Omega^\perp})\\
      & = (X\otimes P_\Omega + X_\gamma \otimes P_{\Omega^\perp}) \otimes P_\Omega + (X\otimes P_\Omega + X_\gamma \otimes P_{\Omega^\perp})_\gamma \otimes P_{\Omega^\perp} \\
      &= X\otimes P_\Omega \otimes P_\Omega + X_\gamma\otimes P_{\Omega^\perp}\otimes P_{\Omega} + X_\gamma \otimes P_{\Omega}\otimes P_{\Omega^\perp} +  |\gamma| X_\gamma \otimes P_{\Omega^\perp} \otimes P_{\Omega^\perp}, 
\end{align*}
while 
\begin{align*}
      &\lambda^\gamma_{H_1,H_2\otimes H_3}(X)= X \otimes P_{\Omega\otimes \Omega}+ X_\gamma \otimes P_{(H_2\otimes H_3)\ominus \mathbb C\Omega}\\
      &=X\otimes P_\Omega\otimes P_\Omega + X_\gamma \otimes P_\Omega\otimes P_{\Omega^\perp} + X_\gamma \otimes P_{\Omega^\perp}\otimes P_\Omega +  X_\gamma \otimes P_{\Omega^\perp} \otimes P_{\Omega^\perp}.
  \end{align*}
They coincide thanks to $|\gamma| X_\gamma=X_\gamma$.  

 The restriction property is easy to see. It then remains to prove the universality of pre-Hilbert spaces. 
 By symmetry it suffices to work on the left universal lift. Let $T\in \Lop(H_1), S \in \Lop(G_1)$ and let $\UU_1\colon H_1 \to G_1$ and $\UU_2\colon H_2 \to G_2$ be possibly non-adjointable $\Omega$-preserving isometries such that $\UU_1 T = S \UU_1$ (the other assumption $\UU_1 T^* = S^* \UU_1$ is unnecessary below). 
 Key formulas are 
 \begin{equation}\label{eq:key}
 \UU_1 T_\gamma = S_\gamma \UU_1, \qquad \UU_2 P_\Omega  = P_\Omega \UU_2 \qquad \text{and}\qquad \UU_2 P_{\Omega^\perp} = P_{\Omega^\perp}\UU_2.
 \end{equation}
 The last two formulas are easy consequences of diagonality of the matrix form of $\UU_2$. The first formula follows from the arguments below: by Lemma \ref{lm:hom} \ref{item:hom1} one sees $(\UU_1)_\gamma T_\gamma = (\UU_1 T)_\gamma = (S \UU_1)_\gamma = S_\gamma (\UU_1)_\gamma$. The diagonality of $\UU_1$ implies that $(\UU_1)_\gamma = |\gamma|\UU_1$, which, together with $|\gamma|T_\gamma = T_\gamma$ and $|\gamma|S_\gamma = S_\gamma$,  implies the desired formula. 
 These formulas yield: 
 \begin{align*}
 \lambda_{G_1,G_2}^\gamma(S)(\UU_1\otimes \UU_2) &= S \UU_1 \otimes P_\Omega \UU_2 +  S_\gamma \UU_1 \otimes P_{\Omega^\perp} \UU_2 \\
 &= \UU_1 T \otimes \UU_2 P_\Omega + \UU_1 T_\gamma \otimes \UU_2 P_{\Omega^\perp} \\
 &= (\UU_1\otimes \UU_2)\lambda_{H_1,H_2}^\gamma(T).
 \end{align*}

 \subcasenp{\ref{item3}} For $\gamma,\delta \in \mathbb T \cup\{0\}$ and $Y \in \Lop(H_2)$, using again Lemma \ref{lm:hom} \ref{item:hom2} we get 
  \begin{align*}
      &\lambda^{\gamma}_{H_1\otimes H_2, H_3}(\rho^{\delta}_{H_1, H_2} (Y))
      = \lambda^{\gamma}_{H_1\otimes H_2, H_3}(P_\Omega\otimes Y + P_{\Omega^\perp}\otimes Y_\delta)\\
      &\ = (P_\Omega\otimes Y + P_{\Omega^\perp}\otimes Y_\delta)\otimes P_\Omega 
      + (P_\Omega\otimes Y + P_{\Omega^\perp}\otimes Y_\delta)_\gamma \otimes P_{\Omega^\perp} \\
      &\ = P_\Omega \otimes Y \otimes P_\Omega 
      + P_{\Omega^\perp}\otimes Y_\delta \otimes P_{\Omega} 
      + P_{\Omega}\otimes Y_\gamma \otimes P_{\Omega^\perp} 
      + |\gamma| P_{\Omega^\perp}\otimes Y_\delta\otimes  P_{\Omega^\perp}
\end{align*}
and 
\begin{align*}
      &\rho^{\delta}_{H_1, H_2 \otimes H_3}(\lambda^{\gamma}_{H_2, H_3} (Y)) 
      =\rho^{\delta}_{H_1, H_2 \otimes H_3}(Y \otimes P_\Omega + Y_\gamma \otimes P_{\Omega^\perp})\\
      &\ =P_\Omega\otimes ( Y \otimes P_\Omega + Y_\gamma \otimes P_{\Omega^\perp}) 
      + P_{\Omega^\perp} \otimes (Y \otimes P_\Omega + Y_\gamma \otimes P_{\Omega^\perp})_\delta \\
      &\ = P_\Omega\otimes  Y \otimes P_\Omega 
      + P_\Omega\otimes Y_\gamma \otimes P_{\Omega^\perp} 
      + P_{\Omega^\perp} \otimes Y_\delta \otimes P_\Omega 
      + |\delta|P_{\Omega^\perp} \otimes Y_\gamma \otimes P_{\Omega^\perp}.
  \end{align*}
Hence the middle associativity of the pair $(\lambda^{\gamma}, \rho^{\delta})$ holds if and only if $|\gamma|Y_\delta = |\delta|Y_\gamma$ for all $Y \in \Lop(H_2)$ and $H_2$, which holds if and only if $(\gamma,\delta) \in J_\otimes$. 
\end{proof}

\subsection{Multi-faced independence arising from the deformed universal lifts} \label{subsec:tensor_lift2}

According to Corollary \ref{cor:multi-correspondence} and Theorem \ref{thm:repr_state}, with any choice of $m$ universal lifts $(\lambda^{\gamma_k}, \rho^{\delta_k})$ from Theorem  \ref{thm:tensor_lift} (one for each face), we get an associated $m$-faced universal product of states. We investigate those universal products in the case of $m=1$ and $m=2$. 

We start with the case $m=1$, where all the universal products are known. Let us recall the definition of the associated universal product of states. Let $(\lambda^{\gamma},\rho^{\delta})$ with $(\gamma,\delta)\in J_\otimes$ be a universal lift and $\utprod{\gamma}{\delta}$ denote both associated universal product of $\ast$-representations and universal product of states.
For two $\ast$\hyp{}representations on pre-Hilbert spaces $\pi\colon A \to \Lop(H)$ and $\sigma\colon B \to \Lop(G)$ their product $\pi \utprod{\gamma}{\delta}\sigma \colon A\sqcup B \to \Lop(H\otimes G)$ is defined by  
\[
\pi \utprod{\gamma}{\delta}\sigma(c) =
\begin{cases}
\lambda_{H,G}^{\gamma}(\pi(c)) &\text{if} \quad c \in A, \\
\rho_{H,G}^{\delta}(\sigma(c)) &\text{if} \quad c \in B. 
\end{cases}
\]
For two states $\phi$ on $A$ and $\psi$ on $B$, their universal product is defined by 
\[
\phi\utprod{\gamma}{\delta}\psi(c) = \langle \Omega, (\pi\utprod{\gamma}{\delta}\sigma)(c)\Omega \rangle, \qquad c \in A \sqcup B, 
\]
where $\pi\colon A \to \Lop(H)$ and $\sigma\colon B \to \Lop(G)$ are any $\ast$\hyp{}representations on pre-Hilbert spaces such that $\phi(a) = \langle \Omega, \pi(a)\Omega \rangle$ and $\psi(b) = \langle \Omega, \sigma(b)\Omega\rangle$. 
For example, for $a_1, a_2 \in A$ and $b_1, b_2 \in B$, the value $\phi\utprod{\gamma}{\delta}\psi(a_1b_1a_2b_2)$ is computed by 
\[
\phi\utprod{\gamma}{\delta}\psi(a_1b_1a_2b_2) = \langle \Omega, \lambda_{H,G}^\gamma(\pi(a_1)) \rho^\delta_{H,G}(\sigma(b_1)) \lambda_{H,G}^\gamma(\pi(a_2)) \rho^\delta_{H,G}(\sigma(b_2))\Omega \rangle. 
\]
The classification of the products $\utprod{\gamma}{\delta}$ for states is provided below. 
\begin{prop} \label{prop:classification_product_of_states} Let $(\gamma,\delta)\in J_\otimes$.
\begin{enumerate}[label=\rm(\arabic*)]
\item $\utprod{\gamma}{\delta}$ is the tensor product if and only if $\gamma=\delta \in \mathbb T$.

\item $\utprod{\gamma}{\delta}$ is the antimonotone product if and only if $\gamma\in \mathbb T, \delta=0$.

\item $\utprod{\gamma}{\delta}$ is the monotone product if and only if $\gamma=0,\delta\in\mathbb T$.

\item $\utprod{\gamma}{\delta}$ is the boolean product if and only if $\gamma=\delta=0$.
\end{enumerate}
\end{prop}
\begin{proof}
We know that $\utprod{\gamma}{\delta}$ is a universal product, and hence it is one of the five universal products. Hence, computing some crucial moments will be enough to identify the product.  
We continue to use the notation from the previous paragraph and adopt the abbreviation $\langle T \rangle = \langle\Omega, T\Omega\rangle$ for operators $T$, $\langle a \rangle = \phi(a)$, $\langle b \rangle = \psi(b)$ for $a\in A$ and $b \in B$, $\lambda^\gamma=\lambda_{H,G}^\gamma$ and $\rho^\delta=\rho_{H,G}^\delta$. 

Let $X=\pi(a) \in \Lop(H)$ and $Y= \sigma(b) \in \Lop(G)$. Then a short calculation shows that 
\begin{align}\label{eq:t-lift-lr} \lambda^{\gamma}(X)\rho^{\delta}(Y)\Omega= X\Omega\otimes\langle Y\rangle\Omega+|\gamma|\langle X\rangle \Omega \otimes (Y-\langle Y\rangle) \Omega + \gamma (X-\langle X\rangle)\Omega \otimes (Y-\langle Y\rangle) \Omega 
\end{align}
and, analogously,
\begin{align}\label{eq:t-lift-rl}\rho^{\delta}(Y)\lambda^{\gamma}(X)\Omega =  \langle X\rangle\Omega\otimes Y\Omega+|\delta|(X-\langle X\rangle) \Omega \otimes \langle Y\rangle \Omega + \delta (X-\langle X\rangle)\Omega \otimes (Y-\langle Y\rangle) \Omega.
\end{align}

For later reference, we move to a more general situation than in the proposition. For $k=1,2$, Take $\gamma_k,\delta_k \in \mathbb T \cup\{0\}$ and denote  $X_k=\pi(a_k) \in \Lop(H)$ and $Y_k= \sigma(b_k) \in \Lop(G)$. Using \eqref{eq:t-lift-lr} and \eqref{eq:t-lift-rl}, we then compute the mixed moment
\begin{align}
&\langle \lambda^{\gamma_1}(X_1)\rho^{\delta_1}(Y_1) \lambda^{\gamma_2}(X_2)\rho^{\delta_2}(Y_2) \rangle \label{eq:tensor_lift_mixed_moments}=\langle\rho^{\delta_1}(Y_1^*)\lambda^{\gamma_1}(X_1^*)\Omega, \lambda^{\gamma_2}(X_2)\rho^{\delta_2}(Y_2)\Omega\rangle \\
&=\langle X_1 \rangle \langle Y_1\rangle  \langle X_2\rangle \langle Y_2\rangle + |\delta_1| \bigl( \langle X_1X_2\rangle-\langle X_1\rangle\langle X_2\rangle\bigr)\langle Y_1\rangle \langle Y_2\rangle \notag \\
&\quad+|\gamma_2| \langle X_1\rangle\langle X_2\rangle \bigl(\langle Y_1Y_2\rangle-\langle Y_1\rangle\langle Y_2\rangle\bigr)
+ \gamma_2\overline{\delta_1} \bigl( \langle X_1X_2\rangle-\langle X_1\rangle\langle X_2\rangle\bigr)\bigl( \langle Y_1Y_2\rangle-\langle Y_1\rangle\langle Y_2\rangle\bigr). \notag
\end{align}

Hence, specializing to the case $\gamma_1=\gamma_2=\gamma$ and $\delta_1=\delta_2=\delta$ and going back to the product of states, we obtain 
\[
\phi\utprod{\gamma}{\delta}\psi(a_1b_1a_2b_2)
=
\begin{cases}
 \langle a_1a_2\rangle\langle b_1b_2\rangle & \text{if} \quad \gamma=\delta \in \mathbb T, \\
 \langle a_1\rangle \langle a_2\rangle\langle b_1b_2\rangle & \text{if} \quad \gamma \in \mathbb T,\delta =0, \\
  \langle a_1a_2\rangle\langle b_1\rangle \langle b_2\rangle& \text{if}\quad \gamma=0, \delta\in\mathbb T, \quad \text{and} \\
  \langle a_1\rangle\langle a_2\rangle\langle b_1\rangle \langle b_2\rangle & \text{if}\quad \gamma=\delta=0.  
\end{cases}
\]
This implies the desired result. Note that the free product does not appear in the list, which would yield $\langle a_1a_2\rangle\langle b_1\rangle \langle b_2\rangle + \langle a_1\rangle \langle a_2\rangle\langle b_1b_2\rangle - \langle a_1\rangle\langle a_2\rangle\langle b_1\rangle \langle b_2\rangle$. 
\end{proof}

\begin{rmq}
The above proof is based on the strong result of Muraki \cite{Muraki13} that the positive universal products are classified into the five ones (see Footnote \ref{fn:Muraki}). However, it is also possible to give a self-contained proof by computing the moments $\phi\utprod{\gamma}{\delta}\psi(w)$ for all alternating words $w$ or by using 
\revise{the arguments preceding Theorem \ref{prop:tensor-tensor-parameters2} below}. 
\end{rmq}

As observed above, no continuous parameters appear in the vacuum expectations in the 1-faced case, in concordance with the fact that there are only five universal products. At the level of the product of $\ast$\hyp{}representations, there is no loss of parameters, but when we go down to the product of states, a large part of information gets lost. However, in the multi-faced case, continuous parameters still remain in some universal products of states, and thus a new notion of independence is obtained.  
To see this, let $(\gamma_k,\delta_k)\in J_\otimes$ for $k=1,2$ and  $(\lambda^{\gamma_1},\rho^{\delta_1})$ and $(\lambda^{\gamma_2},\rho^{\delta_2})$ be universal lifts for the face 1 and face 2, respectively, and let $\uttprod{\gamma_1}{\delta_1}{\gamma_2}{\delta_2}$ denote the associated two-faced universal product of states and also of $\ast$-representations.
Recall that, given two states $\phi$ on $A = A^{(1)}\sqcup A^{(2)}$ and $\psi$ on $B = B^{(1)}\sqcup B^{(2)}$, their universal product is defined by 
\[
\phi \uttprod{\gamma_1}{\delta_1}{\gamma_2}{\delta_2} \psi(c) = \langle \Omega, (\pi\uttprod{\gamma_1}{\delta_1}{\gamma_2}{\delta_2}\sigma)(c)\Omega \rangle, \qquad c \in A \sqcup B, 
\]
where $\pi\colon A \to \Lop(H)$ and $\sigma\colon B\to \Lop(G)$ are  any $\ast$\hyp{}representations on pre-Hilbert spaces such that $\phi(a) = \langle \Omega, \pi(a)\Omega \rangle$ and $\psi(b) = \langle \Omega, \sigma(b)\Omega \rangle$. Recall also that the universal product of $\ast$-representations is defined according to Corollary \ref{cor:multi-correspondence} by  
\[
\pi \uttprod{\gamma_1}{\delta_1}{\gamma_2}{\delta_2} \sigma(c) = 
\begin{cases}
(\lambda_{H,G}^{\gamma_1}\circ \pi) \sqcup (\rho_{H,G}^{\delta_1}\circ \sigma)(c) &\text{if} \quad c \in A^{(1)}\sqcup B^{(1)}, \\
(\lambda_{H,G}^{\gamma_2}\circ \pi) \sqcup ( \rho_{H,G}^{\delta_2}\circ \sigma)(c) &\text{if} \quad c \in A^{(2)}\sqcup B^{(2)}. 
\end{cases}
\]
 
\begin{ex}\label{ex:tensor_mixed_moments} Let $a_1 \in A^{(1)},b_1 \in B^{(1)} ,a_2 \in A^{(2)},b_2 \in B^{(2)}$. 
For notational simplicity, let $\langle T \rangle = \langle \Omega,T\Omega\rangle$ for operators $T$, $\langle a \rangle =\phi(a), \langle b \rangle =\psi(b)$ for $a \in A, b\in B$ and $\lambda^\gamma=\lambda_{H,G}^\gamma,  \rho^\delta=\rho_{H,G}^\delta$ as before. Then, according to the definitions, we have 
\begin{align*}
\phi\uttprod{\gamma_1}{\delta_1}{\gamma_2}{\delta_2}\psi(a_1b_1a_2b_2)
&=\langle \lambda^{\gamma_1}(\pi(a_1))\rho^{\delta_1}(\sigma(b_1))\lambda^{\gamma_2}(\pi(a_2))\rho^{\delta_2}(\sigma(b_2)) \rangle, 
\end{align*}
which was actually computed in \eqref{eq:tensor_lift_mixed_moments}. Hence, we obtain 
\begin{align}
&\phi\uttprod{\gamma_1}{\delta_1}{\gamma_2}{\delta_2}\psi(a_1b_1a_2b_2)   \label{eq:tensor_length_four} \\
&=\langle a_1 \rangle \langle b_1\rangle  \langle a_2\rangle \langle b_2\rangle + |\delta_1| \bigl( \langle a_1a_2\rangle-\langle a_1\rangle\langle a_2\rangle\bigr)\langle b_1\rangle \langle b_2\rangle+|\gamma_2| \langle a_1\rangle\langle a_2\rangle \bigl(\langle b_1b_2\rangle-\langle b_1\rangle\langle b_2\rangle\bigr)\notag \\
&\quad
+ \gamma_2\overline{\delta_1} \bigl( \langle a_1a_2\rangle-\langle a_1\rangle\langle a_2\rangle\bigr)\bigl( \langle b_1b_2\rangle-\langle b_1\rangle\langle b_2\rangle\bigr). \notag
\end{align}
\end{ex}
\begin{ex}\label{ex:tensor-tensor} Let $(a_1,a_2)$ and $(b_1,b_2)$ be  $\uttprod{1}{1}{1}{1}$-independent pairs of elements in a $\ast$-probability space $(A,\phi)$ in the sense of Definition \ref{defi:multi-faced_independence}. This actually means that $A\colonequals\Alg{a_1,a_2}$ and $B\colonequals\Alg{b_1,b_2}$ are tensor independent.
If the elements $a_i$ and $b_i$ are usual random variables defined on a probability space and $\phi$ is the expectation, then this is exactly the usual independence rule for random vectors in $\mathbb C^2$. 
\end{ex}

Let us investigate how the mixed moments depend on the parameters.
In order to compute mixed moments, we need to compute the product of lifted operators acting on the vacuum vector, 
\begin{equation}\label{eq:tensor_mixed_moments}
\kappa_1(T_1)\kappa_2(T_2) \cdots \kappa_n(T_n)\Omega\otimes\Omega, 
\end{equation}
where $\kappa_k$ is one of $\lambda^{\gamma_1},\lambda^{\gamma_2},\rho^{\delta_1},\rho^{\delta_2}$ and $T_k$ is an operator on $\Lop(H)$ or $\Lop(G)$ for $k\in \{1,2,\dots,n\}$.   
A careful observation of the computations of \eqref{eq:tensor_mixed_moments} gives rise to Figure \ref{figure:tensor_tensor}, which should be interpreted as follows. If $T=\begin{pmatrix} \tau & (t')^* \\ t & \hat T \end{pmatrix}$ is the matrix form of an operator $T$ on $\Lop(H)$ or $\Lop(G)$ and $v$ is a vector in one of the four spaces in the figure, two of the components of $T$ can be applied to the corresponding leg of $v$. Then, $\kappa(T)v$ is the weighted sum of those two terms, each multiplied with the parameter indicated in Figure \ref{figure:tensor_tensor}; the symbol $\times$ indicates that an edge is irrelevant for the application of $\kappa(T)$, so that there are in each case exactly two relevant edges starting from each vertex. 

\begin{ex}
    Let $\gamma_1,\gamma_2,\delta_1,\delta_2 \in \mathbb T \cup\{0\}$ and denote  $X_k=\begin{pmatrix} \xi_k & (x_k')^* \\ x_k & \hat X_k \end{pmatrix} \in \Lop(H)$ and $Y_k= \begin{pmatrix} \eta_k & (y_k')^* \\ y_k & \hat Y_k \end{pmatrix} \in \Lop(G)$ for $k=1,2$. We will now compute stepwise 
    \[\rho^{\delta_1}(Y_1) \lambda^{\gamma_2}(X_2)\rho^{\delta_2}(Y_2)\Omega\otimes\Omega.\]
To begin with, $\Omega\otimes\Omega$ is in the top-left space, the two arrows associated with $\rho^{\delta_2}(Y_2)$ are the loop and the downward arrow, both with weight 1. Accordingly,
\[\rho^{\delta_2}(Y_2)\Omega\otimes \Omega= 1\cdot\eta_2\Omega\otimes \Omega + 1\cdot \Omega\otimes y_2.\]
When we next apply $\lambda^{\gamma_2}(X_2)$, the application to the $\Omega\otimes\Omega$-component works analogously (this time with the loop and the rightward arrow), while for $\Omega\otimes y_2\in\Omega\otimes\hat G$, the relevant arrows are the loop with weight $|\gamma_2|$ and the rightward one with weight $\gamma_2$. Accordingly,
\[\lambda^{\gamma_2}(X_2)\rho^{\delta_2}(Y_2)\Omega\otimes\Omega=1\cdot\xi_2\eta_2\Omega\otimes\Omega + 1\cdot \eta_2 x_2\otimes \Omega + |\gamma_2|\cdot \xi_2\Omega\otimes y_2 + \gamma_2\cdot x_2\otimes y_2.\]
In the third step, we have components in all four spaces and find the result
\begin{align*}
\rho^{\delta_1}(Y_1)\lambda^{\gamma_2}(X_2)\rho^{\delta_2}(Y_2) \Omega\otimes\Omega 
&=1\cdot\eta_1  \xi_2 \eta_2  \Omega\otimes \Omega  
+ 1\cdot\xi_2 \eta_2  \Omega\otimes y_1 \notag \\
&\quad+ |\delta_1|\cdot \eta_1 \eta_2 x_2\otimes\Omega
+\delta_1\cdot \eta_2  x_2\otimes y_1 \notag \\
&\quad+1\cdot |\gamma_2| \xi_2 \langle y_1',y_2\rangle\Omega \otimes \Omega + 1\cdot |\gamma_2| \xi_2 \Omega \otimes (\hat Y_1 y_2) \notag \\
&\quad+ \overline{\delta_1}\cdot \gamma_2 \langle y_1',y_2\rangle x_2\otimes \Omega 
+ |\delta_1|\cdot\gamma_2 x_2\otimes (\hat Y_1y_2). \notag
\end{align*}
Applying $\lambda^{\gamma_1}(X_1)$ to the above vector and looking at the coefficient of $\Omega\otimes\Omega$, one deduces that
\begin{align*}
    &\langle \Omega\otimes \Omega, \lambda^{\gamma_1}(X_1) \rho^{\delta_1}(Y_1) \lambda^{\gamma_2}(X_2)\rho^{\delta_2}(Y_2)\Omega\otimes \Omega\rangle\\
    &
    =\xi_1\eta_1\xi_2\eta_2 + |\delta_1|\eta_1\eta_2\langle x_1',x_2\rangle + |\gamma_2|\xi_1\xi_2\langle y_1',y_2\rangle + \overline{\delta_1} \gamma_2 \langle y_1',y_2\rangle \langle x_1',x_2\rangle.
\end{align*}
When we write $\langle T \rangle$ for the vacuum expectation of an operator, we have $\xi_k=\langle X_k\rangle$, $\eta_k=\langle Y_k\rangle$, $\langle x_1',x_2\rangle=\langle X_1X_2\rangle-\langle X_1\rangle\langle X_2\rangle$, and $\langle y_1',y_2\rangle=\langle Y_1Y_2\rangle-\langle Y_1\rangle\langle Y_2\rangle$.
Therefore, we reproduced the result \eqref{eq:tensor_lift_mixed_moments} for the mixed moment $\langle \lambda^{\gamma_1}(X_1) \rho^{\delta_1}(Y_1) \lambda^{\gamma_2}(X_2)\rho^{\delta_2}(Y_2) \rangle $.
\end{ex}


A mixed moment $\langle \kappa_1(T_1)\cdots \kappa_n(T_n)\rangle$ can now be understood combinatorially: regard Figure \ref{figure:tensor_tensor} as a digraph with the four pre-Hilbert spaces as vertices and four weights (with defects denoted as $\times$) on each edge. On each edge, the 1st, 2nd, 3rd and 4th weight correspond to the lifts $\lambda^{\gamma_1},\lambda^{\gamma_2},\rho^{\delta_1}$ and $\rho^{\delta_2}$, respectively. Consider the set of directed paths of length $n$ which start and end in $\Omega\otimes \Omega$. Each edge on such a path associates the weight which  corresponds to the lifts  $\lambda^{\gamma_1},\lambda^{\gamma_2},\rho^{\delta_1},\rho^{\delta_2}$ appearing in \eqref{eq:tensor_mixed_moments}; for example, the $p$-th edge on a path should be equipped with the 3rd weight on it if $\kappa_{n-p+1}=\rho^{\delta_1}$. A path is relevant if none of its associated weights is $\times$, i.e.\ if all its edges are relevant for the application of the corresponding lifted operator $\kappa_k(T_k)$. The total weight of a relevant path is the product of all its associated weights. For each relevant path, apply the corresponding blocks of the operators $T_k$ to the appropriate tensor leg and multiply with the total weight of the path. The sum then yields the sought moment, expressed as a weighted sum of a tensor product of compositions of the blocks of the $T_k$. The total weights which appear in this process may include any monomial consisting of 
\begin{equation}\label{eq:monomial}
|\gamma_i|,|\delta_j|, \gamma_i \overline{\delta_j},\overline{\gamma_i} \delta_j~( i,j\in\{1,2\}), \gamma_1\overline{\gamma_2},\overline{\gamma_1}\gamma_2, \delta_1\overline{\delta_2},\overline{\delta_1}\delta_2.
\end{equation}

\begin{figure}[b]
 \[  
    \xymatrix{
    \mathbb C\Omega\otimes \Omega \longrarrow{1,1,{\times},{\times}} \darrow{{\times},{\times},1,1} \ar@(ur,ul)_{1,1,1,1}&& \hat H\otimes \Omega \longlarrow{1,1,{\times},{\times}} \darrow{{\times},{\times},\delta_1,\delta_2} \ar@(ur,ul)_{1,1,|\delta_1|,|\delta_2|} \\ \Omega\otimes \hat G \uarrow{{\times},{\times},1,1} \longrarrow{\gamma_1,\gamma_2,{\times},{\times}}\ar@(dr,dl)^{|\gamma_1|,|\gamma_2|,1,1} && \hat H \otimes \hat G  \uarrow{{\times},{\times},\overline{\delta_1},\overline{\delta_2}}  \longlarrow{\overline{\gamma_1},\overline{\gamma_2},{\times},{\times}} \ar@(dr,dl)^{|\gamma_1|,|\gamma_2|,|\delta_1|,|\delta_2| }
    } 
\] 
\caption{The way lifting operators by $\lambda^{\gamma_1},\lambda^{\gamma_2},\rho^{\delta_1},\rho^{\delta_2}$ to $H\otimes G$ creates parameters}\label{figure:tensor_tensor}
\end{figure}

\revise{Observe that the numbers listed in \eqref{eq:monomial} are all invariant under the simultaneous rotation 
\begin{equation}\label{eq:rotation1}
(\gamma_1,\delta_1,\gamma_2,\delta_2)\mapsto (\varepsilon\gamma_1,\varepsilon\delta_1,\varepsilon\gamma_2,\varepsilon\delta_2)
\end{equation}
for any $\varepsilon \in \mathbb T$. This implies that 
\begin{equation}\label{eq:rotation2}
    \uttprod{\gamma_1}{\delta_1}{\gamma_2}{\delta_2}  =\uttprod{\varepsilon\gamma_1}{\varepsilon\delta_1}{\varepsilon\gamma_2}{\varepsilon\delta_2}.   
\end{equation}
In particular, the parametrization of the constructed 2-faced product of states $\uttprod{\gamma_1}{\delta_1}{\gamma_2}{\delta_2}$ is not one-to-one. We therefore remove unnecessary parameters and provide a minimal set of parameters. It turns out that the above rotation \eqref{eq:rotation1} is the only source of redundancy. 
To see this, it is convenient to separately discuss products with fixed faces 1 and 2.

If both faces are the tensor product ($\gamma_1 = \delta_1 \in \mathbb T$ and $\gamma_2=\delta_2 \in \mathbb T$, see Proposition \ref{prop:classification_product_of_states}), then \eqref{eq:rotation2} implies $ \uttprod{\gamma_1}{\gamma_1}{\gamma_2}{\gamma_2}  = \uttprod{\zeta}{\zeta}{1}{1},$ where $\zeta := \gamma_1\overline{\gamma_2}$. One can see from \eqref{eq:tensor_length_four} that $\{\uttprod{\zeta}{\zeta}{1}{1}\}_{\zeta \in \mathbb T}$ is a one-to-one parametrization. 

If face 1 is the tensor product ($\gamma_1 = \delta_1 \in \mathbb T$) and face 2 is the monotone product ($\gamma_2=0, \delta_2 \in \mathbb T$), then \eqref{eq:rotation2} entails $ \uttprod{\gamma_1}{\gamma_1}{{}_{\phantom{1}}0}{\delta_2}  = \uttprod{\zeta}{\zeta}{0}{1},$ where $\zeta := \gamma_1\overline{\delta_2}$. To see that the parametrization $\{\uttprod{\zeta}{\zeta}{0}{1}\}_{\zeta \in \mathbb T}$ is one-to-one, we need to evaluate $b_1a_1b_2a_2$ instead of $a_1b_1a_2b_2$ in \eqref{eq:tensor_length_four}.  

For the other cases we can discuss analogously and obtain the following classification that has no redundant parameters. Note that as soon as the boolean face is involved, there is at most one non-zero parameter which can always be normalized without changing the universal product of states. Therefore, in these cases there is no deformation. 

  \begin{thm}
\label{prop:tensor-tensor-parameters2}
    The universal products of states $\uttprod{\gamma_1}{\delta_1}{\gamma_2}{\delta_2}$ are classified into the nine families of products with continuous parameter $\zeta\in\mathbb T$ and the  seven isolated products indicated in Table \ref{table:tensor_lift-products}. 

  \end{thm}
\begin{rmq}
For $\zeta=1$ the diagonal shows the trivial extensions of tensor, antimonotone, monotone and boolean products to the two-faced settings (in these cases, the faces do not matter for calculating mixed moments). 
\end{rmq}

}

\revise{
\begin{table}[b]\caption{Classification of two-faced universal products of states obtained from lifts to the tensor product.}
  \label{table:tensor_lift-products}
  \centering
  \begin{tabular}{|c|c|c|c|c|} \hline
   \multirow{2}{*}{\diagbox[width=2.8cm]{face 1}{face 2}} & tensor  & antimonotone  & monotone & boolean \\ 
    & $\otimes$ $(\lambda^{\alpha_2},\rho^{\alpha_2})$ & $\antimon$ $(\lambda^{\alpha_2},\rho^0)$& $\mon$ $(\lambda^0,\rho^{\alpha_2})$& $\bool$ $(\lambda^0,\rho^0)$\\ \hline
    tensor & \multirow{2}{*}{$ 
    \uttprod{\zeta}{\zeta}{1}{1}$}  & \multirow{2}{*}{$\twofacedprod{\otimes_\zeta}{\antimon_{\phantom{\zeta}}}:=\uttprod{\zeta}{\zeta}{1}{0}$}  & \multirow{2}{*}{$\twofacedprod{\otimes_\zeta}{\mon_{\phantom{\zeta}}}:=\uttprod{\zeta}{\zeta}{0}{1}$} & \multirow{2}{*}{$\twofacedprod{\otimes}{\bool}:=\uttprod{1}{1}{0}{0}$}  \\
    $\otimes$ $(\lambda^{\alpha_1},\rho^{\alpha_1})$
                                                          & &   & & \\ \hline
    antimonotone & \multirow{2}{*}{$\twofacedprod{\antimon_\zeta}{\otimes_{\phantom{\zeta}}}:=\uttprod{\zeta}{0}{1}{1}$}  & \multirow{2}{*}{$\twofacedprod{\antimon_\zeta}{\antimon_{\phantom{\zeta}}}:=\uttprod{\zeta}{0}{1}{0}$}  & \multirow{2}{*}{$\twofacedprod{\antimon_\zeta}{\mon_{\phantom{\zeta}}}:=\uttprod{\zeta}{0}{0}{1}$} & \multirow{2}{*}{$\twofacedprod{\antimon}{\bool}:=\uttprod{1}{0}{0}{0}$}  \\
    $\antimon$ $(\lambda^{\alpha_1},\rho^0)$
                                                          & &   & & \\ \hline
    monotone& \multirow{2}{*}{$\twofacedprod{\mon_\zeta}{\otimes_{\phantom{\zeta}}}:=\uttprod{0}{\zeta}{1}{1}$}  & \multirow{2}{*}{$\twofacedprod{\mon_\zeta}{\antimon_{\phantom{\zeta}}}:=\uttprod{0}{\zeta}{1}{0}$}  & \multirow{2}{*}{$\twofacedprod{\mon_\zeta}{\mon_{\phantom{\zeta}}}:=\uttprod{0}{\zeta}{0}{1}$} & \multirow{2}{*}{$\twofacedprod{\mon}{\bool}:=\uttprod{0}{1}{0}{0}$}  \\
     $\mon$ $(\lambda^0,\rho^{\alpha_1})$
                                                          & &   & & \\ \hline
    boolean & \multirow{2}{*}{$\twofacedprod{\bool}{\otimes}:=\uttprod{0}{0}{1}{1}$}  & \multirow{2}{*}{$\twofacedprod{\bool}{\antimon}:=\uttprod{0}{0}{1}{0}$}  & \multirow{2}{*}{$\twofacedprod{\bool}{\mon}:=\uttprod{0}{0}{0}{1}$} & \multirow{2}{*}{$\boolboolprod:=\uttprod{0}{0}{0}{0}$}  \\
    $\bool$ $(\lambda^0,\rho^0)$
                                                          & &   & & \\ \hline
    \end{tabular}\\[.5em]
($\zeta,\alpha_1,\alpha_2\in \mathbb T$, $\zeta=\alpha_1\overline{\alpha_2}$)
\end{table}
}

It is known that tensor, free and boolean products are symmetric while monotone and antimonotone products are not. 
We will classify all the two-faced symmetric products among $\uttprod{\gamma_1}{\delta_1}{\gamma_2}{\delta_2}$. 

\begin{lm}\label{lm:symmetry}
For all $(\gamma_1,\delta_1), (\gamma_2,\delta_2) \in J_\otimes$ and any two states $\phi$ on $A$ and $\psi$ on $B$,  we have $\phi \uttprod{\gamma_1}{\delta_1}{\gamma_2}{\delta_2} \psi =(\psi \uttprod{\delta_1}{\gamma_1}{\delta_2}{\gamma_2} \phi) \circ \Theta_{A,B}$. 
\end{lm}
\begin{proof}
Recall that $\Theta_{A,B}$ is the natural identification mapping $A \sqcup B \to B \sqcup A$, so we omit it in the following. We need to prove that $( \phi \uttprod{\gamma_1}{\delta_1}{\gamma_2}{\delta_2} \psi)(w) = (\psi \uttprod{\delta_1}{\gamma_1}{\delta_2}{\gamma_2} \phi)(w)$ for every word $w$. 
For notational simplicity, we work on the typical example $w = a_1 b_1 a_2 b_2$ where $a_1 \in A^{(1)}, a_2 \in A^{(2)}, b_1 \in B^{(1)}, b_2 \in B^{(2)}$. Take $\ast$-representations  $\pi \colon A \to \Lop(H)$ and $\sigma \colon B \to \Lop(G)$ which realize $\phi$ and $\psi$, respectively. Then we have 
\begin{align*}
    \phi \uttprod{\gamma_1}{\delta_1}{\gamma_2}{\delta_2} \psi(w) 
    &= \langle \Omega,  \pi \uttprod{\gamma_1}{\delta_1}{\gamma_2}{\delta_2} \sigma (w) \Omega\rangle_{H\otimes G} \\
    &= \langle \Omega, \lambda_{H,G}^{\gamma_1}(\pi(a_1)) \rho_{H,G}^{\delta_1}(\sigma(b_1)) \lambda_{H,G}^{\gamma_2}(\pi(a_2)) \rho_{H,G}^{\delta_2}(\sigma(b_2)) \Omega\rangle_{H\otimes G}.  
\end{align*}
We can switch the left and right lifts here by the symmetry 
\[
 \lambda_{H,G}^{\gamma}(X) = (F_{H,G})^*\rho_{G,H}^\gamma(X) F_{H,G}, 
\]
where $F_{H,G}$ is the flip operator from $H \otimes G$ to $G\otimes H$. Then we have
\begin{align*}
    \phi \uttprod{\gamma_1}{\delta_1}{\gamma_2}{\delta_2} \psi(w) 
    &= \langle \Omega, \rho_{G,H}^{\gamma_1}(\pi(a_1)) \lambda_{G,H}^{\delta_1}(\sigma(b_1)) \rho_{G,H}^{\gamma_2}(\pi(a_2)) \lambda_{G,H}^{\delta_2}(\sigma(b_2)) \Omega\rangle_{G\otimes H},   
\end{align*}
which coincides with $(\psi \uttprod{\delta_1}{\gamma_1}{\delta_2}{\gamma_2} \phi)(w)$ by definition. General words can be treated similarly. 
\end{proof}

\begin{prop}\label{prop:symmetry_tensor}
Let $(\gamma_1,\delta_1), (\gamma_2,\delta_2) \in J_\otimes$. The two-faced universal product $\uttprod{\gamma_1}{\delta_1}{\gamma_2}{\delta_2}$ of states is symmetric if and only if $\gamma_1=\delta_1$ and $\gamma_2=\delta_2$. 
\revise{In other words, the symmetric products are exactly those located at the four corners of Table \ref{table:tensor_lift-products}: $\twofacedprod[.8]{\otimes_\zeta}{\otimes_{\phantom{\zeta}}}~(\zeta\in\mathbb T), \twofacedprod[.6]{\bool}{\otimes}, \twofacedprod[.6]{\otimes}{\bool}, \boolboolprod$.} 
\end{prop}
\begin{proof}
  The ``if'' part follows directly from Lemma \ref{lm:symmetry}. For the ``only if'' part, observe that if one face is equipped with a  non-symmetric product then the whole two-faced product is also  non-symmetric. If $\gamma_1\neq\delta_1$, we 
  \revise{conclude from Proposition \ref{prop:classification_product_of_states}}
  that the first face is equipped with either the monotone or the antimonotone product, hence with a non-symmetric product. Analogously, $\gamma_2\neq\delta_2$ implies that the second face is equipped with a non-symmetric product. 
%
\end{proof}

\begin{rmq} 
A question is whether the two-faced universal product of states $\uttprod{\delta_1}{\gamma_1}{\delta_2}{\gamma_2}$ yields a convolution of probability measures on $\R^2$. In fact, a reasonable definition of convolution exists only for the standard tensor-tensor independence ($\gamma_1=\delta_1=\gamma_2=\delta_2 \in \mathbb T$) --- see Remark \ref{rmq:convolution} for further details. 
\end{rmq}

\section{Lifts to the free product}\label{sec:free_lift}

We investigate the deformation \eqref{eq:main_universal_free_lift} of the standard universal lift \eqref{eq:free_lift} to the free product, motivated by the tensor case (Theorem \ref{thm:tensor_lift}). Together with its version acting from the right side, we are able to even classify all the universal lifts to the free product. 

As a consequence, some new notions of multi-faced independence, including a deformation of bi-freeness, can be defined along the  lines of Subsection \ref{subsec:tensor_lift2}.

\subsection{Classification of universal lifts to the free product}\label{subsec:free_lift1}


We recall here that the free product of pre-Hilbert spaces is defined by 
\[
H_1 \AST H_2 = \mathbb C \Omega \oplus \bigoplus_{n \in \mathbb N, \revise{\mathbf k\in\{1,2\}^n_{\mathrm alt}}
} \hat H_{k_1} \otimes \hat H_{k_2} \otimes \cdots \otimes \hat H_{k_n}. 
\]
We abbreviate the space $\hat H_{k_1} \otimes \hat H_{k_2} \otimes \cdots \otimes \hat H_{k_n} $ to $\hat H_{\mathbf k}$. This serves as a monoidal product on $\PreHilb$; for example, possibly non-adjointable $\Omega$-preserving isometries $\UU_k\colon H_k \to G_k$ $(k=1,2)$ produce the arrow 
$\UU_1 \ast \UU_2\colon H_1 \AST H_2 \to G_1 \AST G_2$ defined by 
\[
\UU_1 \ast \UU_2= \mathrm{id}_{\mathbb C \Omega}\oplus \bigoplus_{n \in \mathbb N, \revise{\mathbf k\in\{1,2\}^n_{\mathrm alt}}}\UU_{k_1}\otimes \UU_{k_2} \otimes \cdots \otimes \UU_{k_n}. 
\]

Making use of the complementary spaces 
\[
H(\ell) = \mathbb C \Omega \oplus \bigoplus_{\substack{n \in \mathbb N, \revise{\mathbf k\in\{1,2\}^n_{\mathrm alt}} \\ k_1\neq \ell }} \hat H_{k_1,k_2,\dots, k_n} \subset H_1\AST H_2, 
\] 
we define canonical unitaries $U=U_{H_1,H_2} \colon H_1 \AST H_2 \to H_1 \otimes H(1)$ and $V=V_{H_1,H_2} \colon H_1 \AST H_2 \to H_2 \otimes H(2)$ by 
\begin{alignat*}{2}
U(\Omega) &= \Omega\otimes \Omega,  \\
U(x) &= x \otimes \Omega, & & x \in \hat H_1, \\
U(x \otimes w) &= x \otimes w, & & x \in \hat H_1, w \in \hat H_{2,1,2,\dots}, \\
U (w) &= \Omega \otimes w, &\qquad& w \in  \hat H_{2,1,2,\dots} 
\end{alignat*}
and similarly
\begin{alignat*}{2}
V(\Omega) &= \Omega\otimes \Omega,  \\
V(y) &= y \otimes \Omega, &\qquad& y \in \hat H_2, \\
V(y \otimes w') &= y \otimes w', &\qquad& y \in \hat H_2, w' \in \hat H_{1,2,1,\dots}, \\
V (w') &= \Omega \otimes w', &\qquad& w' \in  \hat H_{1,2,1,\dots}. 
\end{alignat*}
The free lift $(\fllift^{\rm free}, \frlift^{\rm free})$ is then defined as \eqref{eq:free_lift}, i.e., roughly speaking, by acting on $H_1*H_2$ from the left. 

By symmetry, the actions from the right side of the free product also work. This lift can be defined as 
\begin{align}
\bfllift^{\rm free}_{H_1,H_2}(T) 
&=(R_{H_1,H_2})^\ast \fllift^{\rm free}_{H_1,H_2}(T)R_{H_1,H_2},\qquad T \in \Lop(H_1)
\quad \text{and} \label{eq:free_lift_right}\\
\bfrlift^{\rm free}_{H_1,H_2}(S) 
&=(R_{H_1,H_2})^\ast \frlift^{\rm free}_{H_1,H_2}(S)R_{H_1,H_2}, \qquad S \in \Lop(H_2)
\end{align}
in terms of the unitary operator $R_{H_1,H_2} \in \Lop(H_1\AST H_2)$ which reverses the tensor components: $R_{H_1,H_2}(x_1 \otimes x_2 \otimes \cdots \otimes x_n) = x_n \otimes \cdots \otimes x_2 \otimes x_1$ and $R_{H_1,H_2}(\Omega)=\Omega$.
The above exposition is based on \cite{Voiculescu85} and \cite{Avitzour82}.

Considering the classification of universal lifts to the tensor product in Theorem \ref{thm:tensor_lift}, it is natural to ask whether the mappings 
$\fllift^\gamma_{H_1,H_2}\colon \Lop(H_1) \to \Lop(H_1\AST H_2)$ and $\frlift^{\delta}_{H_1,H_2}\colon \Lop(H_2) \to \Lop(H_1\AST H_2)$ defined by 
\begin{align*}
\fllift^\gamma_{H_1,H_2} (T) 
&= (U_{H_1,H_2})^* (T \otimes P_\Omega + T_\gamma \otimes P_{\Omega^\perp}) U_{H_1,H_2} \quad \text{and} \\
\frlift^\delta_{H_1,H_2} (S) &= (V_{H_1,H_2})^* (S \otimes P_\Omega + S_\delta \otimes P_{\Omega^\perp}) V_{H_1,H_2} 
\end{align*}
are universal lifts to the free product. By symmetry, we also consider the actions from the right side: 
 \begin{align*}
    \bfllift^\gamma_{H_1,H_2}(T)&= (R_{H_1,H_2})^\ast \fllift^{\gamma}_{H_1,H_2}(T)R_{H_1,H_2} \quad  \text{and} \\  
    \bfrlift^\delta_{H_1,H_2}(S)&= (R_{H_1,H_2})^\ast \frlift^{\delta}_{H_1,H_2}(S)R_{H_1,H_2}. 
    \end{align*}
Theorem \ref{thm:free_lift} below answers this question; the result is unexpected because the admissible set of parameters is different from the tensor case. 

Before we come to the main theorem of this section, however, we will illustrate how $\fllift^\gamma$ and $\frlift^\delta$ actually operate on elements of the free product $H_1\AST H_2$. 
Let $\hat T\in \Lop(\hat H_1),x,x'\in \hat H_1,y\in\hat H_2$, and $w\in \hat H_{1,2,1,\ldots}$ an alternating word of arbitrary length (including length 0, in which case $y\otimes w$ simply means $y$). Then
\begin{align*}
    \fllift^\gamma_{H_1,H_2}(\hat T)x&=\hat Tx, 
    &\fllift^\gamma_{H_1,H_2}(\hat T)(x\otimes y\otimes w)&=|\gamma|  (\hat Tx) \otimes y\otimes w,\\
    \fllift^\gamma_{H_1,H_2}(a^*_x)\Omega&=x,
    &\fllift^\gamma_{H_1,H_2}(a^*_x) y\otimes w&=\gamma x \otimes y\otimes w, \\
    \fllift^\gamma_{H_1,H_2}(a_x)x'&=\langle x,x'\rangle\Omega, 
    &\fllift^\gamma_{H_1,H_2}(a_x) x'\otimes y\otimes w&=\overline \gamma \langle x,x'\rangle y \otimes w, \\
    \fllift^\gamma_{H_1,H_2}(P_\Omega)\Omega &=\Omega,
    &\fllift^\gamma_{H_1,H_2}(P_\Omega) y\otimes w&= |\gamma| y \otimes w,  
\end{align*}
and analogous results hold for $\frlift^\delta_{H_1,H_2}$. This behaviour is summarized in Figure \ref{figure:free}.\footnote{Similar to Figure \ref{figure:tensor_tensor}, it is possible to attach two weights to each arrow, using again a defect symbol if an arrow is irrelevant for the lift. For our discussion the simplified version, where we just list the possible different weights that can occur, is sufficient. Later figures built on Figure \ref{figure:free} would otherwise become hardly readable.} 

\begin{figure}[tb]
\[  
    \xymatrix{
   & \ldarrow{1} \hat H_1 \rarrow{\text{$\delta$}} \ar@(u,l)_{1 \text{~or~} |\delta|}  &  \larrow{\text{$\overline{\delta}$}}  \hat H_2 \otimes \hat H_1  \rarrow{\text{$\gamma$}}  \ar@(ur,ul)_{|\gamma|\text{~or~}|\delta|} & \larrow{\text{$\overline{\gamma}$}} \hat H_1 \otimes \hat H_2 \otimes \hat H_1 \rarrow{\text{$\delta$}}  \ar@(ur,ul)_{|\gamma|\text{~or~}|\delta|}  & \larrow{\text{$\overline{\delta}$}} \hat H_2 \otimes \hat H_1 \otimes \hat H_2 \otimes \hat H_1  \rarrow{\text{$\gamma$}} \ar@(ur,ul)_{|\gamma|\text{~or~}|\delta|}  &\larrow{\text{$\overline{\gamma}$}} \cdots 
    \\
    \mathbb C \Omega \ruarrow{1} \rdarrow{1} \ar@(ul,dl)_{1}& & &  
    \\ 
    & \hat H_2 \rarrow{\text{$\gamma$}} \luarrow{1}  \ar@(l,d)_{1 \text{~or~} |\gamma|}  &  \larrow{\text{$\overline{\gamma}$}}  \hat H_1 \otimes \hat H_2  \rarrow{\text{$\delta$}}  \ar@(dl,dr)_{|\gamma|\text{~or~}|\delta|} & \larrow{\text{$\overline{\delta}$}} \hat H_2 \otimes \hat H_1 \otimes \hat H_2 \rarrow{\text{$\gamma$}} \ar@(dl,dr)_{|\gamma|\text{~or~}|\delta|}  & \larrow{\text{$\overline{\gamma}$}} \hat H_1 \otimes \hat H_2 \otimes \hat H_1 \otimes \hat H_2 \rarrow{\text{$\delta$}} \ar@(dl,dr)_{|\gamma|\text{~or~}|\delta|} &\larrow{\text{$\overline{\delta}$}} \cdots 
    } 
\] 
\caption{The way operators lifted by $(\fllift^{\gamma},\frlift^{\delta})$ create parameters as they act on $H_1\AST H_2$}\label{figure:free}
\end{figure}

\begin{thm}\label{thm:free_lift} \ 
\begin{enumerate}[label=\rm(\arabic*)]
    \item The left universal lifts to the free product are classified into the family $\{\fllift^\gamma\}_{\gamma \in \mathbb T}\cup \{\bfllift^\gamma\}_{\gamma \in \mathbb T} \cup \{\fllift^0\}$. Similarly, the right universal lifts to the free product are classified into the family $\{\frlift^\delta\}_{\delta \in \mathbb T}  \cup \{\bfrlift^\delta\}_{\delta \in \mathbb T}\cup \{\frlift^0\}$. Note that $\fllift^0=\bfllift^0$ and $\frlift^0=\bfrlift^0$. 
   
    \item The universal lifts to the free product are classified into the family 
$\{(\fllift^\gamma, \frlift^\delta)\}_{\gamma,\delta \in\mathbb T}\cup \{(\bfllift^\gamma, \bfrlift^\delta)\}_{\gamma,\delta \in \mathbb T}\cup \{(\fllift^0, \frlift^0)\}$. 
\end{enumerate}
We will call $(\fllift^0, \frlift^0) = (\bfllift^0, \bfrlift^0)$ the \emph{boolean lift} to the free product. 
\end{thm}

The proof is divided into two parts: Propositions \ref{prop:free_lift} and \ref{prop:class-free-lift}. The easier and straightforward part is that the proposed lifts are indeed universal lifts. 

\begin{prop}\label{prop:free_lift}\ 
\begin{enumerate}[label=\rm(\arabic*)]
\item\label{item:easier1} Let $\gamma,\delta\in \mathbb T \cup \{0\}$. Then $\fllift^\gamma$ and $\bfllift^\gamma$ are left universal lifts and $\frlift^\delta$ and $\bfrlift^\delta$ are right universal lifts to the free product.
\item\label{item:easier2} Let $(\gamma,\delta)\in \mathbb T^2 \cup\{(0,0)\}$. Then $(\fllift^\gamma, \frlift^\delta)$ and $(\bfllift^\gamma, \bfrlift^\delta)$ are universal lifts to the free product. 
\item\label{item:easier3} Let $\gamma,\delta \in\mathbb C$, $(\gamma,\delta)\notin \mathbb T^2 \cup\{(0,0)\}$. Then neither $(\fllift^\gamma, \frlift^\delta)$ nor $(\bfllift^\gamma, \bfrlift^\delta)$ is  a universal lift to the free product.
\end{enumerate}
\end{prop}

\begin{proof} 
\subcasenps{\ref{item:easier1}} By symmetry, it suffices to work only on $\fllift^\gamma$. 
The mapping $\fllift^\gamma_{H_1,H_2}\colon \Lop(H_1)\to \Lop(H_1 \AST H_2)$ is a $\ast${\hyp}homomorphism because $T\mapsto T \otimes P_\Omega + T_\gamma \otimes P_{\Omega^\perp}$ is a $\ast$-homomorphism due to Theorem \ref{thm:tensor_lift}.

\subcase{Left restriction} 
Note that for $h\in H_1$ we have $U_{H_1,H_2} h=h\otimes\Omega$ and $U_{H_1,H_2}^* h\otimes \Omega=h$. Therefore, 
for $T \in \Lop(H_1)$ and $h \in H_1$ we have 
\begin{align*}
\fllift^\gamma_{H_1,H_2}&(T)h= U_{H_1,H_2}^* (T \otimes P_\Omega + T_\gamma \otimes P_{\Omega^\perp})U_{H_1,H_2} (h) = U_{H_1,H_2}^* Th\otimes \Omega =Th.
\end{align*}


\subcase{Left universality of pre-Hilbert spaces} 
Let $\UU_k\colon H_k \to G_k$ be a possibly non-adjointable $\Omega$-preserving isometry for $k=1,2$ and let $S \in \Lop(G_1)$ and $T \in \Lop(H_1)$ be such that $\UU_1 T = S \UU_1$; the other condition $\UU_1 T^* = S^* \UU_1$ is unnecessary below. 
Recall that $\UU_1\AST \UU_2$ is a linear operator from $H_1 \AST H_2$ into $G_1\AST G_2$ which applies $\UU_1$ for all $\hat H_1$ components and $\UU_2$ for all $\hat H_2$ components simultaneously.

Note that the formulas 
\[
\UU_1 T_\gamma = S_\gamma \UU_1, \qquad P_\Omega (\UU_1 \AST \UU_2) =  (\UU_1 \AST \UU_2) P_\Omega, \qquad P_{\Omega^\perp} (\UU_1 \AST \UU_2) = (\UU_1 \AST \UU_2) P_{\Omega^\perp} 
\]
hold similarly to \eqref{eq:key} and, additionally, the formulas
\begin{align*}
U_{G_1,G_2} (\UU_1\AST \UU_2) 
&= [\UU_1 \otimes (\UU_1 \AST \UU_2)] U_{H_1,H_2} \quad \text{and}  \\
U_{G_1,G_2}^* [\UU_1 \otimes (\UU_1 \AST \UU_2)] 
&=  (\UU_1 \AST \UU_2)U_{H_1,H_2}^* 
\end{align*}
hold. Then we have 
\begin{align*}
\fllift^\gamma_{G_1,G_2}(S) (\UU_1 \AST \UU_2) 
&= U_{G_1,G_2}^* [S\otimes P_\Omega + S_\gamma\otimes P_{\Omega^\perp}] [\UU_1 \otimes (\UU_1 \AST \UU_2)] U_{H_1,H_2} \\
&= U_{G_1,G_2}^* [S \UU_1\otimes (P_\Omega (\UU_1 \AST \UU_2)) + S_\gamma \UU_1 \otimes (P_{\Omega^\perp} (\UU_1 \AST \UU_2))]  U_{H_1,H_2} \\
&= U_{G_1,G_2}^* [\UU_1 T \otimes ((\UU_1 \AST \UU_2) P_\Omega) + \UU_1 T_\gamma \otimes ((\UU_1 \AST \UU_2)P_{\Omega^\perp}) ] U_{H_1,H_2} \\ 
&=U_{G_1,G_2}^* [\UU_1 \otimes (\UU_1 \AST \UU_2)] [T \otimes P_\Omega +  T_\gamma \otimes P_{\Omega^\perp}] U_{H_1,H_2} \\ 
&=(\UU_1\AST \UU_2)\fllift^\gamma_{H_1,H_2}(T),  
\end{align*}
as desired.

\subcase{Left associativity}\label{prop:free_lift:left-associativity}
    As a first step, we convince ourselves that it is enough to prove the associativity condition on operators $a^*_x$ with $x\in\hat H_1$. 
    Those operators generate the $*$-algebra of adjointable finite-rank operators on $H_1$. Now let $T\in \Lop(H_1)$ have arbitrary rank. Any element $w\in H_1 * H_2 * H_3$ lies in a subspace of the form $ G_1 * H_2 * H_3$ with $G_1$ finite dimensional. Let $W\colon G_1\hookrightarrow H_1$ be the embedding, which is an adjointable isometry, $P=WW^*$ the projection onto $G_1$, and $P^\perp=1-P$ the projection onto $G_1^\perp$. Then 
    \[\fllift^\gamma_{H_1,H_2*H_3}(T)w=\fllift^\gamma_{H_1,H_2*H_3}(TP)w+\fllift^\gamma_{H_1,H_2*H_3}(TP^\perp)w\]
    and, by universality, $P^\perp W=0_{G_1\to H_1} = W 0_{G_1\to G_1}$ implies
    \[\fllift^\gamma_{H_1,H_2*H_3} (P^\perp)w=\fllift^\gamma_{H_1,H_2\AST H_3} (P^\perp)(W*{\id}*{\id})w=(W*{\id}*{\id})\fllift^\gamma_{G_1,H_2\AST H_3}(0_{G_1\to G_1})w=0.\]
    Therefore, $\fllift^\gamma_{H_1,H_2*H_3}(T)w=\fllift^\gamma_{H_1,H_2*H_3}(TP)w$. 
    A similar reasoning can be applied to show \[
    \fllift^\gamma_{H_1*H_2,H_3}(\fllift^\gamma_{H_1,H_2}(T))w=\fllift^\gamma_{H_1*H_2,H_3}(\fllift^\gamma_{H_1,H_2}(TP))w\] because
    \begin{align*}
        &P^\perp W=0_{G_1\to H_1}=W0_{G_1\to G_1}\\ 
        &\implies \fllift^\gamma_{H_1,H_2}(P^\perp) (W*{\id})=(W*{\id}) \fllift^\gamma_{G_1,H_2}(0_{G_1\to G_1})\\
        &\implies \fllift^\gamma_{H_1*H_2,H_3}(\fllift^\gamma_{H_1,H_2}(P^\perp))(W*{\id}*{\id})=(W*{\id}*{\id})\fllift^\gamma_{G_1*H_2,H_3}(\fllift^\gamma_{G_1,H_2}(0_{G_1\to G_1}))=0. 
    \end{align*}
    In conclusion, left associativity on finite-rank operators implies
    \[\fllift^\gamma_{H_1*H_2,H_3}(\fllift^\gamma_{H_1,H_2}(T))w=\fllift^\gamma_{H_1*H_2,H_3}(\fllift^\gamma_{H_1,H_2}(TP))w=\fllift^\gamma_{H_1,H_2*H_3}(TP)w=\fllift^\gamma_{H_1,H_2*H_3}(T)w\]
    and, therefore general left associativity.
    
    For the rest of the proof, we will write the parameters again explicitly, but abbreviate $\fllift_{H_1,H_2}^\gamma$ to $\fllift_{1,2}^\gamma$, $\fllift_{H_1*H_2,H_3}^\gamma$ to $\fllift_{12,3}^\gamma$ and so on.
    \revise{It turns out that the computation of $\fllift_{12,3}^\gamma(\fllift_{1,2}^\gamma (a^*_x))w$ and $\fllift_{1,23}^\gamma (a^*_x)w$ for words $w$ varies depending on the space to which the initial tensor factor of $w$ belongs, and moreover, if the initial factor belongs to $\hat H_2$ then the computation also depends on whether the word $w$ contains a factor from $\hat H_3$ or not. It is therefore convenient to discuss separately the words of the forms}
    \[
       \Omega,\quad
       x'\otimes w_{123},\quad
       y\otimes w_{12},\quad
       z\otimes w_{123},\quad
       y\otimes w_{12}\otimes z \otimes w_{123}
    \]
    \revise{with $x'\in \hat H_1, y\in \hat H_2,z\in \hat H_3$, and $w_{12}, w_{123}$ alternating words of arbitrary length (possibly of length 0) with tensor factors from those $\hat H_k$ that are indicated by the subscripts (of course $w_{12}$ must not start with a factor from $\hat H_2$, and there are similar constraints on the first and second $w_{123}$).}  
    For elements of the last type, we get
\footnote{Figure \ref{figure:free} will be helpful for seeing how $\gamma$'s appear in the main calculations. Also, it might be helpful for the understanding to examine the first equality in a bit more detail, as the underlying principle is important when using Figure \ref{figure:free} in the calculations. The map $\fllift^\gamma_{1,2}(a^*_x)$ has a matrix decomposition of the form $\begin{pmatrix}
        0 & 0 \\
        t & \hat T 
    \end{pmatrix}$ with respect to $H_1*H_2=\mathbb C \Omega \oplus ((H_1*H_2)\ominus \mathbb C \Omega)$.
    Since $y\otimes w_{12}\in (H_1*H_2)\ominus \mathbb C \Omega$, $[\fllift^\gamma_{1,2}(a^*_x)]_\gamma (y\otimes w_{12})= |\gamma| \fllift^\gamma_{1,2}(a^*_x) (y\otimes w_{12})$. Roughly speaking, in order to determine the correct result with parameter, it is enough to check from which part of the space pre-image and image come when one ignores the parameter.}
    \begin{align*}
        \fllift_{12,3}^\gamma(\fllift_{1,2}^\gamma (a^*_x)) \underbrace{y \otimes w_{12}\otimes z \otimes w_{123}}_{\in \widehat{H_1*H_2}\otimes \hat H_3 \otimes \cdots} &= |\gamma| \underbrace{[\fllift_{1,2}^\gamma (a^*_x) \underbrace{y\otimes w_{12}}_{\in\hat H_2\otimes \cdots}]\otimes z\otimes w_{123}}_{\in \widehat{H_1*H_2}\otimes \hat H_3 \otimes \cdots}\\
        &=|\gamma|\gamma \underbrace{x\otimes y\otimes w_{12}}_{\in \hat H_1\otimes \hat H_2\otimes \cdots}\otimes z \otimes w_{123} \qquad \text{and} \\
        \fllift_{1,23}^\gamma (a^*_x)   \underbrace{y \otimes w_{12}\otimes z \otimes w_{123}}_{\in \widehat{H_2*H_3} \otimes \cdots} 
        &=\gamma \underbrace{x\otimes y\otimes w_{12}\otimes z \otimes w_{123}}_{\in \hat H_1\otimes \widehat{H_2*H_3}\otimes \cdots}, 
    \end{align*}
so they coincide.     

    The operators $\fllift_{12,3}^\gamma(\fllift_{1,2}^\gamma (a^*_x))$ and $\fllift_{1,23}^\gamma (a^*_x)$ are easily seen to agree on elements of the first four types. 
    This shows the left associativity.

\subcasenp{\ref{item:easier2}, \ref{item:easier3}}  We only need to see whether middle associativity holds or not. It is enough to check it on operators $a^*_y$ with $y\in\hat H_2$; this can be done analogously to the corresponding statement in the proof of left associativity. 
 \revise{The computation of $\fllift^\gamma_{12, 3}(\frlift^\delta_{1,2}(a^*_y))w$ and of $\frlift^\delta_{1,23}(\fllift^\gamma_{2,3}(a^*_y))w$ for words $w$ varies depending on the space to which the initial tensor factor of $w$ belongs, and moreover, if the initial factor belongs to $\hat H_1$ or $\hat H_3$ then the computation also depends on whether or not the word $w$ contains factors from both $\hat H_1$ and $\hat H_3$. It is therefore convenient to discuss separately the words of the forms}
\[
    \Omega,\quad
    y'\otimes w_{123},\quad
    x\otimes w_{12},\quad
    z\otimes w_{23},\quad
    x\otimes w_{12} \otimes z \otimes w_{123},\quad
    z\otimes w_{23}\otimes x\otimes w_{123}
\]
with $x\in \hat H_1, y'\in \hat H_2, z\in\hat H_3$, and $w_{i_1\ldots i_k}$ an alternating word of arbitrary length with tensor factors from $\hat H_{i_1},\ldots, \hat H_{i_k}$ (length $0$ allowed, the obvious restrictions on the first letter needed to obtain an alternating word are assumed).

When we compare $\fllift^\gamma_{12, 3}(\frlift^\delta_{1,2}(a^*_y))$ with $\frlift^\delta_{1,23}(\fllift^\gamma_{2,3}(a^*_y))$ on $x\otimes w_{12} \otimes z \otimes w_{123}$, we obtain (cf.\ Figure \ref{figure:free})

\begin{align*}
    \fllift_{12, 3}^\gamma(\frlift_{1,2}^\delta(a^*_y))\underbrace{z\otimes w_{23}\otimes x\otimes w_{123}}_{\in \hat H_3\otimes \widehat{H_1*H_2}\otimes\cdots}
    &=\gamma \underbrace{(\frlift_{1,2}^\delta(a^*_y)\Omega)\otimes z\otimes w_{23}\otimes x\otimes w_{123}}_{\in \widehat{H_1*H_2}\otimes \hat H_3 \otimes \widehat{H_1*H_2}\otimes\cdots}\\
    &=\gamma y\otimes z\otimes w_{23}\otimes x\otimes w_{123}\qquad \text{and} \\
    \frlift_{1,23}^\delta(\fllift_{2,3}^\gamma(a^*_y)) \underbrace{z\otimes w_{23}\otimes x\otimes w_{123}}_{\in \widehat{H_2*H_3}\otimes \hat H_1 \otimes \cdots}
    &=|\delta|\underbrace{[\fllift_{2,3}^\gamma(a^*_y)\underbrace{z\otimes w_{23}}_{\in \hat H_3\otimes\cdots}]\otimes x\otimes w_{123}}_{\in \widehat{H_2*H_3}\otimes \hat H_1 \otimes \cdots}\\
    &=|\delta|\gamma \underbrace{y\otimes z\otimes w_{23}}_{\in \hat H_2\otimes \hat H_3\otimes\cdots}\otimes x\otimes w_{123}. 
\end{align*}
The two terms agree if and only if $\gamma=|\delta|\gamma$, i.e.\ if and only if $|\delta|=1$ or $\gamma=0$.
By symmetry, the results on $x\otimes w_{12} \otimes z \otimes w_{123}$ agree if and only if $|\gamma|=1$ or $\delta=0$. In all other listed cases, the two operators are easily seen to agree without any restriction on the choice of parameters. Therefore, middle associativity is equivalent to $|\gamma|=|\delta|=1$ or $\gamma=\delta=0$.
\end{proof}

\begin{ex}\label{ex:monotone_lift}
We make a comment on an instructive counterexample. Let $W^{\antimon}_{H_1,H_2}\colon H_1 \otimes H_2 \hookrightarrow H_1 \AST H_2$ be the obvious embedding as the subspace $\mathbb C \Omega \oplus \hat H_1 \oplus \hat H_2 \oplus (\hat H_1\otimes \hat H_2)$ and $P^{\antimon}_{H_1,H_2}=W^{\antimon}_{H_1,H_2}(W^{\antimon}_{H_1,H_2})^*$ the projection onto that subspace. One might think of defining the left and right lifts  
\begin{align*}\label{eq:naive_lift}
   \allift_{H_1,H_2}^\gamma(T)&= P^{\antimon}_{H_1,H_2}(U_{H_1,H_2})^*(T \otimes P_\Omega + T_\gamma \otimes P_{\Omega^\perp})U_{H_1,H_2}P^{\antimon}_{H_1,H_2}\\
   \arlift_{H_1,H_2}^\delta(T)&= P^{\antimon}_{H_1,H_2}(V_{H_1,H_2})^*(T \otimes P_\Omega + T_\delta \otimes P_{\Omega^\perp})V_{H_1,H_2}P^{\antimon}_{H_1,H_2}
\end{align*}
for $\gamma,\delta \in \mathbb T \cup \{0\}$. However, this is not a universal lift to the free product in the sense of Definition \ref{defi:universal_lift} unless $\gamma=\delta=0$; more precisely,  $\allift^\gamma_{H_1,H_2\AST H_3}(\id_{H_1})$ and $\allift^\gamma_{H_1\AST H_2, H_3}(\allift^\gamma_{H_1,H_2}(\id_{H_1}))$ are different on $\hat H_3 \otimes \hat H_2$ unless $\gamma=0$ (the latter vanishes while the former does  not) and analogously $\arlift^\delta$ is not right associative unless $\delta=0$; in the case $\gamma=\delta=0$, the lift $(\allift^0, \arlift^0)$ coincides with the boolean lift $(\fllift^0,\frlift^0)$ and hence is universal. If $\gamma\ne0$, then the left lift $\allift^\gamma$ satisfies the left restriction property and left universality of pre-Hilbert spaces, but does not satisfy the left associativity. This means that we cannot drop the left associativity axiom, by contrast to the tensor case (see Remark \ref{rem_lifts} and the proof of Theorem \ref{thm:tensor_lift}).

Actually, the lift $(\allift^1,\arlift^0)$ coincides with the lift which Liu proves in \cite[Proposition 2.7]{Liu19} to implement antimonotone independence, and which he suggests to use for a definition of free-antimonotone independence in \cite[Definition 3.4]{Liu19}.\footnote{To be precise, Liu suggests to use the opposite version to define free-monotone independence, but we prefer to make our argument in the antimonotone setting for notational reasons. Liu's approach is to regard Muraki's ``monotone product Hilbert space'' (i.e.\ the tensor product with the monotone lifts) \cite[pages 4,5]{Muraki00p} naturally as a subspace of the free product Hilbert space. If you compare with Liu's paper, be aware of a confusing typo in his definition of the spaces $\mathcal X({\mon},i)$ and $\mathcal X({\antimon},i)$; the direct sum has to range over $i\geq i_1 > \cdots > i_n$ and $i_1<\cdots < i_n \leq i$, respectively.}
That it is not associative means that the corresponding product of representations is not associative. As the antimonotone product of states is associative, one could hope that the product of states on two-faced algebras defining Liu's free-antimonotone independence is also associative, but this is not the case. 
Let $H=\mathbb C\Omega \oplus \mathbb C x$ be a 2-dimensional Hilbert space with $x$ a unit vector. 
Denote by $A^{(f)},A^{(m)}$ two copies of the $*$-algebra $\Lop(H)$. 
Then $A=A^{(f)}\sqcup A^{(m)}$ is a two-faced $*$-algebra. Let $\pi={\id}\sqcup{\id}\colon A\to \Lop(H)$ denote the obvious $*$-representation and $\varphi(\cdot)=\langle\Omega, \pi(\cdot) \Omega\rangle$ the corresponding state. 
We compare the two states $(\varphi \odot  \varphi)\odot \varphi$ with $\varphi\odot (\varphi\odot \varphi)$ on $A\sqcup A\sqcup A=A^{(f)}\sqcup A^{(m)}\sqcup A^{(f)}\sqcup A^{(m)}\sqcup A^{(f)}\sqcup A^{(m)}$. We denote the operator $a^*_x$ in the $f$-face or $m$-face of the $k$'th factor $A$ by $x_k^{(f)}$ and $x_k^{(m)}$, respectively, and the vector $x$ in the $k$'th copy of $H$ in $H*H*H$ by $x_k$, $k=1,2,3$. Then we find
\begin{align*}
    \left((\varphi \odot  \varphi)\odot \varphi\right)\left( (x_1^{(m)} x_3^{(f)} x_2^{(f)})^*(x_1^{(m)} x_3^{(f)} x_2^{(f)})\right)
    =\left\|((\pi\odot\pi)\odot \pi) (x_1^{(m)} x_3^{(f)} x_2^{(f)})\Omega\right\|^2, \\
    \left(\varphi \odot(\varphi\odot \varphi)\right)\left( (x_1^{(m)} x_3^{(f)} x_2^{(f)})^*(x_1^{(m)} x_3^{(f)} x_2^{(f)})\right)
    =\left\|(\pi\odot(\pi\odot \pi)) (x_1^{(m)} x_3^{(f)} x_2^{(f)})\Omega\right\|^2
\end{align*}
and
\begin{align*}
    ((\pi\odot\pi)\odot \pi) (x_1^{(m)} x_3^{(f)} x_2^{(f)})\Omega
    = ((\pi\odot\pi)\odot \pi) (x_1^{(m)}) x_3\otimes x_2
    = \allift_{12,3}^1(\allift_{1,2}^1(x_1^{(m)})) x_3\otimes x_2=0, 
\end{align*}
while
\begin{align*}
    (\pi\odot(\pi\odot \pi)) (x_1^{(m)} x_3^{(f)} x_2^{(f)})\Omega
    = (\pi\odot(\pi\odot \pi)) (x_1^{(m)}) x_3\otimes x_2
    = \allift_{1,23}^1(x_1^{(m)}) x_3\otimes x_2=x_1\otimes x_3\otimes x_2. 
\end{align*}
This shows that the free-antimonotone product of states is not associative.
\end{ex}

The preceding example emphasizes the value of our definition of associative universal lifts and Theorem \ref{thm:repr_state} which assures that universal lifts can be combined to produce well-behaved multi-faced independences.

We formulate the remaining part of Theorem \ref{thm:free_lift} for further reference.

\begin{prop}\label{prop:class-free-lift}\ 
\begin{enumerate}[label=\rm(\arabic*)]
\item The only left universal lifts to the free product are those discussed in Theorem \ref{thm:free_lift}, namely $\fllift^\gamma$ and $\bfllift^\gamma$ with $\gamma\in\mathbb T \cup\{0\}$.
\item  The only universal lifts to the free product are those discussed in Theorem \ref{thm:free_lift}, namely $(\fllift^\gamma,\frlift^\delta)$ and $(\bfllift^\gamma,\bfrlift^\delta)$ with $(\gamma,\delta)\in\mathbb T^2$ and the boolean lift $(\fllift^0,\frlift^0)=(\bfllift^0,\bfrlift^0)$.
\end{enumerate}
\end{prop}

The proof is based on several lemmas below. The following result might be well-known.
\begin{lm}\label{lm:orthonormal}
  Let $n\in\N$ and let $H$ be a pre-Hilbert space of dimensions at least $n$. Then 
    \[H^{\otimes n}= \operatorname{span}(e_1\otimes \cdots \otimes e_n : \{e_k\}_{k=1}^n \text{\rm~is an orthonormal system of $H$}).
    \]
\end{lm}

\begin{proof}
    The simple key observation is that for two orthonormal vectors $e$ and $f$, it holds that 
    \begin{align}\label{eq:orthonormal2}e\otimes e = \frac{1}{2}(e+f)\otimes (e-f) + \frac{1}{2}(e+\mathrm{i}f)\otimes (e-\mathrm{i}f) + \frac{1+\mathrm{i}}{2}e\otimes f -\frac{1+\mathrm{i}}{2} f\otimes e\end{align}
    is a linear combination of tensors whose factors are orthonormal and lie in the span of $e$ and $f$.
    With this in mind, we prove the claim by induction on $n$. For $n=1$, the statement is trivial. 
    Suppose the statement holds for $n$ and let $x_1\otimes \cdots \otimes x_{n+1}\in H^{\otimes n+1}$, $\dim H\geq {n+1}$. Since $\dim H\geq n$, it follows from the induction hypothesis that $x_1\otimes\cdots \otimes x_n=\sum_{i\in I} \alpha_i e_1^{(i)}\otimes \cdots \otimes e_n^{(i)}$ for a finite number of orthonormal systems $(e_1^{(i)},\ldots, e_n^{(i)})$ in $H$ and $\alpha_i\in\mathbb C$. Since $\dim H\geq{n+1}$, for each $i\in I$, there is a unit vector $e_{n+1}^{(i)}$ orthogonal to all  $e_k^{(i)}$ , $k\leq n$ with $x_{n+1}= \sum_{k=1}^{n+1} \beta_k^{(i)} e_k^{(i)}$. 
     
    Then,
    \[x_1\otimes \cdots \otimes x_{n+1}=\sum_{i\in I} \sum_{k=1}^{n} \alpha_i\beta_k^{(i)} e_1^{(i)} \otimes \cdots \otimes e_n^{(i)}\otimes e_k^{(i)} +
    \sum_{i\in I}\alpha_i\beta_{n+1}^{(i)} e_1^{(i)} \otimes \cdots \otimes e_n^{(i)}\otimes e_{n+1}^{(i)}.\]
    Each $e_1^{(i)}\otimes \cdots \otimes e_n^{(i)}\otimes e_k^{(i)}$ with $k\in\{1,\ldots, n\}$ can be written in the desired form using \eqref{eq:orthonormal2} with $e=e_k^{(i)}$ and $f=e_{n+1}^{(i)}$, so the proof is finished.
\end{proof}

\begin{lm}\label{lm:permutations} 
  For $k_0\in\{1,2\}$, $n\ge 1$, and \revise{$\mathbf k= (k_1, \dots, k_n) \in \{1,2\}^n_{\mathrm{alt}}$, 
  }
    we denote by $\perm{k_0}{\mathbf k}$ the set of permutations $\sigma$ of $\{0,1,\dots,n\}$ such that 
    $(k_{\sigma(0)}, k_{\sigma(1)} , \cdots, k_{\sigma(n)})\in \{1,2\}^{n+1}_{\mathrm{alt}}$.

Let $\lambda$ be a left universal lift to the free product. Then  
there exists a family of universal coefficients $(\coef{1}{\mathbf k}{\sigma})_{\sigma\in \perm{1}{\mathbf k}}$ such that 
for all pre-Hilbert spaces $H_1,H_2$ and all $z_i \in \hat H_{k_i}, 0 \le i \le n$ ($k_0\colonequals1$) we have 
\begin{equation}\label{eq:univ}
\lambda_{H_1, H_2}(a_{z_0}^*)z_1\otimes  z_2 \otimes \cdots \otimes z_n = \sum_{\sigma \in \perm{1}{\mathbf k}}\coef{1}{\mathbf k}{\sigma} z_{\sigma(0)}\otimes z_{\sigma(1)}\otimes  \cdots \otimes z_{\sigma(n)}.   
\end{equation}
\end{lm}

Note that the left associativity axiom will not be used in the proof.

\begin{proof} 

We fix $n \ge 1$, \revise{$\mathbf k= (k_1, \dots, k_n) \in \{1,2\}^n_{\mathrm{alt}}$} 
and $k_0\colonequals1$. Also, fix the pre-Hilbert spaces $K_1, K_2$ with 
\begin{align}\label{eq:K-ONB}
    \mathfrak E=\{e_i:i\in\{0,\ldots, n\},k_i=1\},\quad \mathfrak F=\{f_i:i\in\{0,\ldots, n\},k_i=2\} 
\end{align} 
orthonormal bases of $\hat K_1, \hat K_2$, respectively. Our plan is to reduce the general statement to the special case of $\lambda_{K_1,K_2}$ (step 1), then establish the coefficients by calculating $\lambda_{K_1,K_2}(a^*_{e_0})g_1\otimes \cdots \otimes g_n$, $g_i=e_i$ or $g_i=f_i$ according to whether $k_i=1$ or $k_i=2$ (steps 2 and 3) and finally show that the formula holds for all tensors (step 4).

For notational simplicity we will sometimes assume that $\mathbf{k}$ is of the form $(2,1,2,1,\ldots, 2)$. This allows to write $e_i$ or $f_i$ instead of $g_i$ and also has the nice effect that $k_i=1$ if and only if $i$ is even. The arguments are however not specialized to this situation and work for arbitrary $\mathbf k$ with obvious adjustments.

In this proof, we will abbreviate the notation $\lambda_{H_1,H_2}(a_z^*)$ to $\lambda_{H_1,H_2}(z)$.

\subcase{Step 1: reduction to specific pre-Hilbert spaces}  
    Suppose that the statement has been proved for $K_1, K_2$. More precisely, suppose that there exists a family of coefficients $(c_{K_1,K_2}(\sigma))_{\sigma \in \perm{1}{\mathbf k}}$ such that 
\[
\lambda_{K_1, K_2}(u_0)u_1\otimes  u_2 \otimes \cdots \otimes u_n = \sum_{\sigma \in \perm{1}{\mathbf k}}c_{K_1,K_2}(\sigma) u_{\sigma(0)}\otimes u_{\sigma(1)}\otimes  \cdots \otimes u_{\sigma(n)}  
\]
for all $u_i \in \hat K_{k_i}, 0 \le i \le n$. 

Let $H_1,H_2$ be general pre-Hilbert spaces and let $z_i \in \hat H_{k_i}, 0 \le i \le n$. Let $G_j \colonequals \Span(\Omega, z_i: 0 \le i \le n \text{~with~} k_i=j) \subset H_j$ for $j=1,2$. Since $\dim (G_j) \le \dim(K_j)$, there are $\Omega$-preserving isometries $W_j \colon G_j \to K_j, j=1,2$.  Let $\tilde z_i$ be the image of $z_i$ by the embedding $W_{k_i}$. 
 Since $W_1 a_{z_0}^* = a_{\tilde z_0}^* W_1$, the left universality implies 
 \begin{align*}
 (W_1 \ast W_2)\lambda_{G_1,G_2}(z_0)  z_1 \otimes \cdots  \otimes z_n 
 &=  \lambda_{K_1,K_2}(\tilde  z_0) (W_1 \ast W_2) z_1 \otimes \cdots   \otimes z_n  \\
 &=  \lambda_{K_1,K_2}(\tilde z_0) \tilde z_1 \otimes \cdots  \otimes \tilde z_n \\
 &= \sum_{\sigma \in \perm{1}{\mathbf k}} c_{K_1,K_2}(\sigma) \tilde z_{\sigma(0)}\otimes \tilde z_{\sigma(1)}\otimes  \cdots \otimes \tilde z_{\sigma(n)} \\
 &= (W_1\ast W_2)\sum_{\sigma \in \perm{1}{\mathbf k}}c_{K_1,K_2}(\sigma) z_{\sigma(0)}\otimes z_{\sigma(1)}\otimes  \cdots \otimes z_{\sigma(n)}. 
 \end{align*}
 Because $W_1 \ast W_2$ is an isometry, we obtain
 \begin{equation}\label{eq:univ1}
  \lambda_{G_1,G_2}(z_0)  z_1 \otimes \cdots  \otimes z_n= \sum_{\sigma \in \perm{1}{\mathbf k}}c_{K_1,K_2}(\sigma) z_{\sigma(0)}\otimes z_{\sigma(1)}\otimes  \cdots \otimes z_{\sigma(n)}. 
\end{equation}
Let $W_k'\colon G_k \hookrightarrow H_k$ be the embeddings. Because of $W_1' a_{z_0}^* = a_{z_0}^* W_1'$, using the left universality of pre-Hilbert spaces and \eqref{eq:univ1} we obtain
\begin{align*}
    \lambda_{H_1,H_2}(z_0) z_1 \otimes \cdots  \otimes z_n 
    &=\lambda_{H_1,H_2}(z_0) (W_1' \ast W_2') z_1 \otimes \cdots   \otimes z_n \\
    &= (W_1' \ast W_2')\lambda_{G_1,G_2}(z_0)  z_1 \otimes \cdots  \otimes z_n \\
    &= \sum_{\sigma \in \perm{1}{\mathbf k}}c_{K_1,K_2}(\sigma) z_{\sigma(0)}\otimes z_{\sigma(1)}\otimes  \cdots \otimes z_{\sigma(n)},   
\end{align*}
as desired.

\subcase{Step 2: existence of universal coefficients} 
From the preceding arguments, we may and do work only on the pre-Hilbert spaces $K_1,K_2$ from now on. As announced above, we also assume $\mathbf{k}=(2,1,2,\dots,2)$; in particular, $n$ is odd and $k_i=1$ if and only if $i$ is even. 

In order to describe a basis for $K_1 \AST K_2$, let ${\mathbf W}$ be the set of words $w = w_1 w_2 \cdots w_m, m \ge0,$ with $w_1,w_2,\dots, w_m\in\{0,1,\ldots, n\}$, such that either 
\begin{itemize}
\item $w_1,w_3,w_5, \dots \in 2\Z$ and $w_2,w_4,w_6, \dots \in 2\Z+1$ or 
\item $w_1,w_3,w_5, \dots \in 2\Z+1$ and $w_2,w_4,w_6, \dots \in 2\Z$.  
\end{itemize}
We include the empty word (corresponding to $m=0$) as an element in ${\mathbf W}$. 
For $w =w_1 w_2 \dots w_m\in {\mathbf W}$,  we define $w(\mathfrak E,\mathfrak F) \in K_1 \AST K_2$ for orthonormal bases as in \eqref{eq:K-ONB} to be 
\[
w(\mathfrak E,\mathfrak F)= 
\begin{cases}
e_{w_1} \otimes f_{w_2} \otimes e_{w_3} \otimes f_{w_4}\otimes  \cdots, & \text{if $w_1$ is even}, \\
f_{w_1} \otimes e_{w_2} \otimes f_{w_3} \otimes e_{w_4}\otimes  \cdots, & \text{if $w_1$ is odd}, \\ 
\Omega, & \text{if}~w=\emptyset. 
\end{cases}
\]

Because $\{w(\mathfrak E,\mathfrak F)\}_{w \in \mathbf W}$ is a basis of $K_1\ast K_2$ (as a vector space), there exists a family of coefficients $\alpha(w) \in \mathbb C$, which is finitely supported as a function of $w$, such that
\[
\lambda_{K_1,K_2}(e_0) f_1 \otimes e_2 \otimes \cdots \otimes e_{n-1} \otimes f_{n}  = \sum_{w \in {\mathbf W}} \alpha (w) w(\mathfrak E,\mathfrak F). 
\]
The coefficients $\alpha(w)$ are independent of the choice of the orthonormal bases $\mathfrak E$ and $\mathfrak F$. To see this, we take different orthonormal bases $\mathfrak E' = \{e_{i}':k_i=1\}$ and $\mathfrak F'=\{f_{j}':k_j=2\}$ and define the $\Omega$-preserving isometries (unitaries) $W_1 \colon K_1 \to K_1$ and $W_2\colon K_2\to K_2$ by $W_1(e_i) = e_i'$ and $W_2(f_j)=f_j'$. Then we have $W_1 a_{e_0}^* = a_{e_0'}^* W_1$, and hence, by the left universality of pre-Hilbert spaces, we get 
\begin{align*}
\lambda_{K_1,K_2}(e_0') f_{1}' \otimes e_{2}' \otimes \cdots \otimes e_{n-1}' \otimes f_{n}'   
&= \lambda_{K_1,K_2}(e_0') (W_1 \AST W_2) f_1 \otimes e_{2} \otimes \cdots \otimes e_{n-1} \otimes f_{n} \\ 
&=  (W_1 \AST W_2) \lambda_{K_1,K_2}(e_{0})  f_{1} \otimes e_{2} \otimes \cdots \otimes e_{n-1} \otimes f_{n} \\
&=  (W_1 \AST W_2) \sum_{w \in {\mathbf W}} \alpha(w) w(\mathfrak E,\mathfrak F) \\
&= \sum_{w \in {\mathbf W}} \alpha(w) w(\mathfrak E',\mathfrak F'). 
\end{align*}

\subcase{Step 3: vanishing of irrelevant coefficients}  Let 
$\mathfrak E= \{e_{i}:k_i=1\}$ and $\mathfrak F=\{f_{j}:k_j=2\}$ be orthonormal bases of $\hat K_1$ and $\hat K_2$ as before. Fix for some time $m \in \{0,2,\dots, n-1\}$, i.e.\ $m$ is such that $k_m=1$. We define an orthonormal basis $\mathfrak E^\theta= \{e_{i}^\theta:k_i=1\}$ of $\hat K_1$ for $\theta \in[-\pi,\pi]$ by setting $e_{m}^\theta\colonequals e^{\mathrm{i}\theta} e_{m}$ and $e_{i}^\theta\colonequals  e_{i}$ for $i \ne m$.  
Using the established fact in Step 2 above, we have 
\[
\lambda_{K_1,K_2}(e_{0}^\theta) f_{1} \otimes e_{2}^\theta \otimes \cdots \otimes e_{n-1}^\theta \otimes f_{n}   = \sum_{w \in {\mathbf W}} \alpha(w) w(\mathfrak E^\theta,\mathfrak F). 
\]
It is obvious that 
\[
\lambda_{K_1,K_2}(e_{0}^\theta) f_{1} \otimes e_{2}^\theta \otimes \cdots \otimes e_{n-1}^\theta \otimes f_{n} = e^{\mathrm{i}\theta} \lambda_{K_1,K_2}(e_{0}) f_{1} \otimes e_{2} \otimes \cdots \otimes e_{n-1} \otimes f_{n}. 
\]
On the other hand, we can see that 
\[
\sum_{w \in {\mathbf W}} \alpha(w) w(\mathfrak E^\theta,\mathfrak F) = h_0 +  e^{\mathrm{i}\theta}h_1  + e^{2\mathrm{i}\theta}h_2 + \cdots  
\]
for some $h_j\in K_1 \AST K_2, j \ge0$ independent of $\theta$. By the uniqueness of Fourier series expansion, $h_j=0$ for all $j\ne1$. This means that $\alpha(w)\ne0$ only if the even number $m$ appears in $w$ exactly once. Similar arguments  for $\mathfrak F$ show that  $\alpha(w)\ne0$ only if a fixed odd number $m$ appears in $w$ exactly once. Since $m$ can be an arbitrary number from $\{0,\ldots, n\}$, we conclude that $\alpha(w)$ can be non-zero only if $w$ is a permutation of the word $012\dots n$: 
\begin{equation}\label{eq:perm}
\lambda_{K_1,K_2}(e_{0}) f_{1} \otimes e_{2} \otimes \cdots \otimes e_{n-1} \otimes f_{n}   = \sum_{\sigma \in \perm{1}{\mathbf k}} \alpha(\sigma) w_\sigma (\mathfrak E,\mathfrak F),  
\end{equation}
where $w_\sigma \colonequals  \sigma(0) \sigma(1)\dotsm \sigma(n) \in {\mathbf W}$.

\subcase
{Step 4: extending the formula to all vectors}  
The multilinear mapping $\hat K_{k_0} \times \hat K_{k_1}  \times \cdots \times \hat K_{k_{n}} \ni (z_0, z_1 \dots, z_{n})\mapsto \lambda_{K_1,K_2}(z_{0}) z_{1} \otimes \cdots  \otimes z_{n} \in K_1\AST K_2$ can be regarded as a linear mapping $\Lambda \colon \hat K_{k_0} \otimes \hat K_{k_1}  \otimes \cdots \otimes \hat K_{k_{n}} \to K_1\AST K_2$. 

Let $(z_0, z_1 \dots, z_{n}) \in \hat K_{k_0} \times \hat K_{k_1}  \times \cdots \times \hat K_{k_{n}}$. 
Applying Lemma \ref{lm:orthonormal} to the $\hat K_1$ components and to the $\hat K_2$ components separately, we can find finitely many orthonormal bases $\mathfrak E^{(m)}= \{e_{i}^{(m)}:k_i=1\}$ and $\mathfrak F^{(\ell)}=\{f_{j}^{(\ell)}:k_j=2\}$, $m \in M, \ell \in L$ and some coefficients $\beta_{m}, \gamma_\ell \in \mathbb C$    such that  
\begin{equation}\label{eq:orthonormal}
z_0 \otimes z_1 \otimes \cdots \otimes z_{n} = \sum_{m\in M, \ell \in L} \beta_{m}\gamma_\ell e_0^{(m)} \otimes  f_1^{(\ell)} \otimes e_2^{(m)} \otimes \cdots \otimes f_{n}^{(\ell)}. 
\end{equation}

Combining \eqref{eq:perm} and \eqref{eq:orthonormal} we have
\begin{align*}
\lambda_{K_1,K_2}(z_{0}) z_{1} \otimes \cdots \otimes z_{n}   
&=  \Lambda (z_0 \otimes z_1 \otimes \cdots \otimes z_{n}) \\
&= \sum_{m\in M, \ell \in L}\beta_{m}\gamma_\ell  \Lambda (e_0^{(m)} \otimes  f_1^{(\ell)} \otimes e_2^{(m)} \otimes \cdots \otimes f_{n}^{(\ell)}) \\
&= \sum_{m\in M, \ell \in L} \beta_{m}\gamma_\ell \lambda_{K_1,K_2}(e_0^{(m)})  f_1^{(\ell)} \otimes e_2^{(m)} \otimes \cdots \otimes f_{n}^{(\ell)} \\
&= \sum_{\sigma \in \perm{1}{\mathbf k}} \alpha(\sigma) \sum_{m\in M, \ell \in L}  \beta_{m}\gamma_\ell \sigma (\mathfrak E^{(m)},\mathfrak F^{(\ell)})  \\
&=  \sum_{\sigma \in \perm{1}{\mathbf k}} \alpha(\sigma) z_{\sigma(0)} \otimes \cdots \otimes z_{\sigma(n)}, 
\end{align*}
as desired. 
\end{proof}

\begin{lm}\label{lm:class-left-free-lift}
  Let $\lambda$ be a left universal lift to the free product with with universal coefficients $(\coef{1}{\mathbf k}{\sigma}:\mathbf k\in\{1,2\}^*_{\mathrm{alt}},\sigma\in \perm{1}{\mathbf k})$. Assume that  $\lambda_{H_1,H_2}(a^*_x)y=\gamma x\otimes y$ for all pre-Hilbert spaces $H_1,H_2$ and all $x \in \hat H_1, y\in \hat H_2$; in other words, $c^1_{(2)}(\mathrm{id}_{\{0,1\}})=\gamma$ and $c^1_{(2)}(\tau_{\{0,1\}})=0$, $\tau_{\{0,1\}}$ the transposition of $0$ and $1$. Then $\gamma \in \mathbb T \cup\{0\}$ and for all $n\in\mathbb N_0$, $\mathbf k=(k_1,\ldots, k_n)\in \{1,2\}^n_{\mathrm{alt}}$ and $\sigma\in \perm{1}{\mathbf k}$,
  \begin{align}\label{eq:coeff-lambda}
  c_{\mathbf k}^{1}(\sigma)=
    \begin{cases}
      \gamma & \text{if $\sigma=\mathrm{id}_{\{0,\ldots,n\}}$ and $k_1=2$},\\
      0 &\text{else},
    \end{cases}
  \end{align}
  i.e.\ $\lambda = \fllift^\gamma$ as defined in Theorem \ref{thm:free_lift}.
\end{lm}

\begin{proof}
In the proof, $H_1,H_2,H_3$ will always denote arbitrary pre-Hilbert spaces of dimensions at least $2$. 
  To shorten notation a bit, we will write elements of free product pre-Hilbert spaces without the tensor signs between factors from $\hat H_{k_i}$, for example $x_1x_2x_3\in H_1*H_2$ means that $x_i\in \hat H_{k_i}$ with $\mathbf k=(k_1,k_2,k_3)\in \{1,2\}^3_{\mathrm{alt}}$. Also, for the lift, we will write $\lambda=\lambda_{1,2}=\lambda_{H_1,H_2}$, $\lambda_{1,23}=\lambda_{H_1,H_2\ast H_3}$ and $\lambda_{12,3}=\lambda_{H_1\ast H_2,H_3}$.

\subcase{Proof of $\gamma \in \mathbb T \cup\{0\}$}    
The free product space $H_1*H_2$ has a natural $\N$-grading given by word length. From Lemma \ref{lm:permutations}, it easily follows that $\lambda(a^*_x)$ is a homogeneous operator of degree $1$. Therefore,
    \[\langle \lambda(a^*_x)y,x' y'\rangle= \overline \gamma \langle x,x'\rangle \langle y,y'\rangle \quad \text{and}\quad \langle \lambda(a^*_x)x'',x' y'\rangle= 0 \quad (\forall\ x,x',x''\in\hat H_1, y,y'\in\hat H_2) \]
    are enough to conclude $\lambda(a_x)x' y'=\overline\gamma \langle x,x'\rangle y'$.
    Accordingly, for unit vectors $x\in\hat H_1$ and $y\in\hat H_2$, we have $\lambda(a_x)x y=\overline \gamma y$ and $\lambda (a_xa^*_x)y = |\gamma|^2y$ for $y\in\hat H_2$. Since $\lambda$ is a $*$-homomorphism and $a_x a_x^*$ is the projection onto $\mathbb C\Omega$, $\lambda(a_x a_x^*)$ is a projection too, so we conclude that its eigenvalue $|\gamma|^2$ is $0$ or $1$, i.e.\ $\gamma\in\mathbb T\cup\{0\}$. In particular, $\fllift^\gamma$ is a left universal lift by Proposition \ref{prop:free_lift}\ref{item:easier1}.

\subcase{Proof of \eqref{eq:coeff-lambda}}
Our plan is to prove \eqref{eq:coeff-lambda} by induction on $n$,  the length of $\mathbf k$.
    For $n=1$, one half is the assumption that $\lambda(a_{x_0}^*)y=\gamma x_0y$ and the other half is the obvious observation that 
    \[\lambda(a^*_{x_0}) x = 0.\]
    (This follows either from the left restriction axiom or from the fact that $x=\lambda(a^*_{x})\Omega$ and the homomorphism property or from  Lemma \ref{lm:permutations}.)

    \revise{Now assume that  \eqref{eq:coeff-lambda} holds for for all $\mathbf k \in\{1,2\}^n_{\mathrm{alt}}$ of length $n\ge1$.  Note that this means $\lambda(a_x^*)$ and $\fllift^\gamma(a_x^*)$ coincide on $\hat H_{\mathbf k}$ for every $x \in \hat H_1$ and $\mathbf k \in\{1,2\}^{n}_{\mathrm{alt}}$. We divide the remaining arguments to show that this implies  \eqref{eq:coeff-lambda} for all $\mathbf k$ of length $n+1$ into some steps.}
 

   \subcase{Step 1: some preparatory formulas} We convince ourselves that %
    \begin{align}\label{eq:free-lift-a_x}
        \lambda(a_{x_0})x_1x_2 \dotsm x_{n+1}&=\fllift^\gamma(a_{x_0})x_1x_2 \dotsm x_{n+1} \quad\text{and}\\
        \label{eq:free-lift-hat_T}
        \lambda(\hat T)x_1x_2 \dotsm x_{n+1}&=\fllift^\gamma(\hat T) x_1x_2 \dotsm x_{n+1}
    \end{align} for all $H_1,H_2$, all $x_0\in\hat H_1$, all $\hat T\in \Lop(\hat H_1)$, all $\mathbf k=(k_1,\ldots,k_{n+1})\in \{1,2\}^{n+1}_{\mathrm{alt}}$, and all alternating words $x_1\cdots x_{n+1}\in\hat H_{\mathbf k}$. Using again that $\fllift^\gamma(a^*_{x_0})$ and $\lambda(a^*_{x_0})$ are homogeneous operators of degree 1, \eqref{eq:free-lift-a_x} follows from
    \begin{align*}
        \langle x_1x_2 \dotsm x_{n},\fllift^\gamma (a_{x_0}) y_1y_2 \dotsm y_{n+1}\rangle &=\langle \fllift^\gamma (a^*_{x_0}) x_1x_2 \dotsm x_{n}, y_1y_2 \dotsm y_{n+1}\rangle \\
        &= \langle \lambda (a^*_{x_0}) x_1x_2 \dotsm x_{n}, y_1y_2 \dotsm y_{n+1}\rangle\\
        &=\langle  x_1x_2 \dotsm x_{n},\lambda (a_{x_0}) y_1y_2 \dotsm y_{n+1}\rangle 
    \end{align*}
    for all pairs of alternating words $x_1x_2 \dotsm x_{n}$ and $y_1y_2 \dotsm y_{n+1}$ in $H_1*H_2$ of length ${n}$ and ${n+1}$, respectively.
    
    If $\hat T$ is a finite-rank operator, \eqref{eq:free-lift-hat_T} follows easily from \revise{the induction hypothesis}
    and \eqref{eq:free-lift-a_x} by writing $\hat T=\sum a^*_x a_{x'}$ as a sum of rank one operators:  
    \begin{align*}
    \lambda(\hat T)x_1x_2 \dotsm x_{n+1}
    &=\sum\lambda(a^*_x)\lambda(a_{x'})x_1x_2 \dotsm x_{n+1}\\
    &=\sum\lambda(a^*_x)\fllift^\gamma(a_{x'})x_1x_2 \dotsm x_{n+1}\\
    &=\sum\fllift^\gamma(a^*_x)\fllift^\gamma(a_{x'})x_1x_2 \dotsm x_{n+1}\\
    &=\fllift^\gamma(\hat T) x_1x_2 \dotsm x_{n+1}
    \end{align*}
  because both $\lambda$ and $\fllift^\gamma$ are $\ast$-homomorphisms and
    \[\fllift^\gamma(a_{x'})x_1x_2 \dotsm x_{n+1}=\begin{cases}
        \overline{\gamma}\langle x',x_1\rangle x_2\dotsm x_{n+1},& x_1\in\hat H_1\\
        0,& x_1\in \hat H_2
    \end{cases}\]
    either has length $n$ or vanishes.
    For general $\hat T$, \eqref{eq:free-lift-hat_T} follows from the arguments used in the proof of \revise{Proposition \ref{prop:free_lift} (\hyperref[prop:free_lift:left-associativity]{Left associativity})}.

 \subcase{Step 2: a key consequence of left associativity}  
We define $G_1:=H_1 \AST  H_2$ and $G_2:=H_3$, so that $\lambda_{12,3}=\lambda_{G_1,G_2}$. Let $x_0 \in \hat H_1$ and denote $y_0:=x_0$ as an element in $G_1=H_1 \AST  H_2$. Let $y_1\cdots y_{n+1}\in \hat H_{\mathbf j}\subset\hat G_{\mathbf k}$ be an alternating word, where $\mathbf j=(j_1,\ldots, j_{n+1})\in\{2,3\}^{n+1}_{\mathrm{alt}}$
and $\mathbf k:=(j_1-1,\ldots,j_{n+1}-1)\in\{1,2\}^{n+1}_{\mathrm{alt}}$. Furthermore, we assume that 
\begin{equation}\label{eq:key_assumption}
    \fllift_{12,3}^\gamma(1-P)y_1y_2 \dotsm y_{n+1}=0
\end{equation}
 where $P \in \Lop(H_1*H_2)$ is the projection onto $\mathbb C\Omega$. Then we get 
    \begin{align*}
        \lambda_{12,3}(\lambda_{1,2}(a^*_{x_0})) y_1y_2 \dotsm y_{n+1}
        = \lambda_{12,3} (\lambda_{1,2}(a^*_{x_0})P) y_1y_2 \dotsm y_{n+1} + \lambda_{12,3} (\lambda_{1,2}(a^*_{x_0})(1-P)) y_1y_2 \dotsm y_{n+1}.
    \end{align*}
    Clearly, $\lambda_{1,2}(a^*_{x_0})P=a^*_{y_0}$. By \eqref{eq:free-lift-hat_T} applied to $\hat T\colonequals 1-P$, \[\lambda_{12,3}(1-P)y_1y_2 \dotsm y_{n+1}
      = \fllift_{12,3}^\gamma(1-P)y_1y_2 \dotsm y_{n+1},\]
      which vanishes by assumption. 
    Since $\lambda_{12,3}$ is a homomorphism, we conclude
    \begin{align*}\label{eq:a^*_x-a^*_y}
    \lambda_{12,3}(\lambda_{1,2}(a^*_{x_0})) y_1y_2 \dotsm y_{n+1}
        = \lambda_{12,3} (a^*_{y_0}) y_1y_2 \dotsm y_{n+1}.
    \end{align*}
Moreover, according to Lemma \ref{lm:permutations}, the RHS equals $\sum_{\sigma\in \perm{1}{\mathbf k}} \coef{1}{\mathbf k}{\sigma} y_{\sigma(0)}y_{\sigma(1)} \dotsm  y_{\sigma(n+1)}$
. 

 On the other hand, the left associativity implies    
 \[ \lambda_{12,3}(\lambda_{1,2}(a^*_{x_0})) y_1y_2 \dotsm y_{n+1} = \lambda_{1,23}(a^*_{x_0})(y_1y_2 \dotsm y_{n+1}) = \gamma y_0 y_1 \dotsm y_{n+1}.  
 \]
Note that the last equality follows from the very assumption of this lemma because  the element $y_1y_2 \dotsm y_{n+1}$ is regarded as a word of length 1 in the space $H_1 \AST(H_2 \AST H_3)$. 

 Combining the above together yields the identity
 \begin{equation} \label{eq:key_identity}
     \sum_{\sigma\in \perm{1}{\mathbf k}} \coef{1}{\mathbf k}{\sigma} y_{\sigma(0)}y_{\sigma(1)} \dotsm  y_{\sigma(n+1)} = \gamma y_0 y_1 \dotsm y_{n+1}.   
 \end{equation} 
The universal coefficients can be determined if we are allowed to compare the coefficients.  The details are as follows. 

\subcase{Step 3 -- Case 1: $j_1=3~(k_1=2)$}
In this case the key assumption \eqref{eq:key_assumption} of Step 2 holds because
    \[\fllift_{12,3}^\gamma(1-P)y_1y_2 \dotsm y_{n+1} 
    =|\gamma| ((1-P)\Omega) y_1y_2 \dotsm y_{n+1} =0,\]
    so that \eqref{eq:key_identity} holds for all alternating words $y_1 y_2 \cdots y_{n+1}$ of $\hat H_{\mathbf j}$. Because the pre-Hilbert spaces $H_1,H_2,H_3$ were arbitrary, we can assume that their dimensions are sufficiently large (dimensions $n+1$ are enough). We then take a unit vector $y_0\in \hat H_1$, an orthonormal system $\{y_1, y_3, y_5, \dots\}$ in $\hat H_3$ and an orthonormal system $\{y_2, y_4, y_6, \dots\}$ in $\hat H_2$. Because the set $\{y_{\sigma(0)}y_{\sigma(1)} \dotsm  y_{\sigma(n+1)} \}_{\sigma \in \perm{1}{\mathbf k}}$ is orthonormal in $H_1 \AST H_2 \AST H_3$ and hence is linearly independent, one can compare the coefficients of identity \eqref{eq:key_identity} to conclude $\coef{1}{\mathbf k}{\sigma} =0$ for all $\sigma \ne \id$ and $\coef{1}{\mathbf k}{\id}= \gamma$, i.e.\  \eqref{eq:coeff-lambda} for $k=(2,1,2,1,\dots)$ of length $n+1$.

\subcase{Step 3 -- Case 2: $j_1=2~(k_1=1), \gamma =0$} 
   Using $\gamma=0$ we again obtain \[\fllift_{12,3}^\gamma(1-P)y_1y_2 \dotsm y_{n+1} 
    =|\gamma| ((1-P) y_1)y_2 \dotsm y_{n+1} =0,\]
    so that \eqref{eq:key_identity} holds for all alternating words $y_1 y_2 \cdots y_{n+1}$ of $\hat H_{\mathbf j}$. From the same reasoning as in the previous case, we can compare the coefficients to conclude 
       $\coef{1}{\mathbf k}{\sigma}=0$ for all permutations $\sigma\in\mathcal P^1_{\mathbf k}$, i.e.\  \eqref{eq:coeff-lambda} for $k=(1,2,1,\dots)$ of length $n+1$.

\subcase{Step 3 -- Case 3: $j_1=2~(k_1=1), |\gamma|=1$} Although the key assumption \eqref{eq:key_assumption} of Step 2 is not satisfied,   \eqref{eq:coeff-lambda} can be proved more directly and easily in this case. For all pre-Hilbert spaces $H_1,H_2$, $x_0\in H_1$ and alternating words $x_1\cdots x_{n+1}\in \hat H_{\mathbf k}$, \revise{using the induction hypothesis and the homomorphism property of $\lambda$ yields
    \begin{align*}
    \lambda(a^*_{x_0})x_1x_2 \dotsm x_{n+1}&=\gamma^{-1}\lambda(a^*_{x_0})\fllift^\gamma(a^*_{x_1}) x_2\dotsm x_{n+1}
    =\gamma^{-1}\lambda(a^*_{x_0})\lambda(a^*_{x_1})x_2\dotsm x_{n+1}\\
    &=\gamma^{-1}\lambda(a^*_{x_0}a^*_{x_1})x_2\dotsm x_{n+1}=0,
    \end{align*}}
    showing $c^1_{\mathbf k}(\sigma)=0$ for all $\sigma\in \mathcal P_{\mathbf k}^1$, i.e.\  \eqref{eq:coeff-lambda} for $k=(1,2,1,\dots)$ of length $n+1$. 

Through the above Case 1 -- Case 3, the desired  \eqref{eq:coeff-lambda} is fully proved.   
\end{proof}

We are now ready to finish the proof of Theorem \ref{thm:free_lift}, i.e. prove Proposition \ref{prop:class-free-lift}. 
\begin{proof}[Proof of Proposition \ref{prop:class-free-lift}]
    Let $\lambda$ be a left universal lift. We will follow the notation in Lemma \ref{lm:class-left-free-lift}. We denote elements of pre-Hilbert spaces $H_1,H_2,H_3$ by, $x,x',\ldots\in \hat H_1$, $y,y',\ldots\in\hat H_2$, $z,z',\ldots\in\hat H_3$, respectively, without further mentioning; all those vectors are assumed nonzero. 
    
    \revise{By Lemma \ref{lm:permutations}, 
    there are some universal constants  $\alpha,\beta,\gamma,\delta\in\mathbb C$ such that 
    \begin{align}
     \lambda(a^*_x)y&=\gamma x y + \delta y x,  \label{eq:gamma_delta} \\
     \lambda(a^*_x) x'y &= \alpha x y x' + \beta x' y x. \label{eq:alpha_beta}
    \end{align}}
    We will need to understand how an operator $\hat T\in \Lop(\hat H_1)$ is lifted. Without loss of generality assume $\hat T=\sum a^*_x a_{x'}$ is a finite-rank operator \revise{--- see the corresponding arguments in the proof of Lemma \ref{lm:unique}}. Similar to the first part in the proof of Lemma \ref{lm:class-left-free-lift}, one concludes $\lambda(a_x)x' y'=\overline \gamma\langle x,x'\rangle y'$, so that 
    \begin{align}
    \lambda(\hat T) x''y 
    &= \sum \lambda(a^*_x)\lambda(a_{x'})x''y 
    =\sum \lambda(a^*_x)\overline\gamma \langle{x',x''}\rangle y \label{eq:lift-T_hat-2}  \\
    &=\sum\langle x',x''\rangle (|\gamma|^2 xy + \overline \gamma \delta yx) 
    = |\gamma|^2 (\hat Tx'')y + \overline\gamma\delta y(\hat Tx''). \notag
    \end{align}
    
    Now we want to use associativity to show that at least one of the parameters $\gamma,\delta$ must vanish. 
    Of course,
\begin{equation}\label{eq:associativity1}
    \lambda_{1,23}(a^*_x)yz= \gamma xyz + \delta yzx.
\end{equation}
    It takes a bit more effort to evaluate 
    \[\lambda_{12,3}(\lambda_{1,2}(a^*_x))yz =
    \lambda_{12,3}(\lambda_{1,2}(a^*_x)P)yz+ \lambda_{12,3}(\lambda_{1,2}(a^*_x)(1-P))yz. \]
    As in the proof of Lemma \ref{lm:class-left-free-lift}, $\lambda_{1,2}(a^*_x)P=a^*_x\in \Lop(H_1*H_2)$. By \eqref{eq:alpha_beta} we get 
    \[\lambda_{12,3}(\lambda_{1,2}(a^*_x)P)yz= \lambda_{12,3}(a^*_x)yz=\alpha xzy + \beta yzx.\] 
  By Lemma \ref{lm:permutations}, the operator $\lambda_{1,2}(a^*_x)(1-P)$ lies in $\Lop(\widehat{H_1*H_2})$, so we can use \eqref{eq:lift-T_hat-2} to obtain
    \begin{align*}
    \lambda_{12,3}(\lambda_{1,2}(a^*_x)(1-P))yz
    &= |\gamma|^2 (\lambda_{1,2}(a^*_x)y)z + \overline\gamma\delta z(\lambda_{1,2}(a^*_x)y)  \\
    &= |\gamma|^2 \gamma xyz + |\gamma|^2 \delta yxz + |\gamma|^2\delta zxy + \overline\gamma\delta^2 zyx. \notag
    \end{align*}
    Putting both parts together, we find
    \begin{align}
        \label{eq:associativity2}
        \lambda_{12,3}(\lambda_{1,2}(a^*_x))yz
        =\alpha xzy + \beta yzx
        +|\gamma|^2 \gamma xyz + |\gamma|^2 \delta yxz + |\gamma|^2\delta zxy + \overline\gamma\delta^2 zyx.
    \end{align}
    \revise{We are allowed to compare the coefficients of $yxz$ in \eqref{eq:associativity1} and \eqref{eq:associativity2}, see the arguments in [Step 3 -- Case 1] of the previous lemma, and can conclude $|\gamma|^2\delta=0$.} 
    
    By symmetry, the analogous result for right lifts holds as well.
    
    From Lemma \ref{lm:class-left-free-lift} with its obvious variations for left lifts with $\lambda(a_x^*)y=\gamma yx$, right lifts with $\rho(a_x^*)y=\delta xy$ and right lifts with $\rho(a_x^*)y=\delta yx$, we now established a complete classification of left and right universal lifts to the free product. 
    
    A combination $(\fllift^\gamma,\bfrlift^\delta)$ with $\gamma,\delta\in\mathbb T$ cannot be associative, e.g.
    \[\bfrlift^\delta_{1,23}(\fllift^\gamma_{2,3}(a^*_y))xzxz=\gamma xzxyz \neq \delta xyzxz = \fllift^\gamma_{12,3}(\bfrlift^\delta_{1,2}(a^*_y))xzxz.\]
    Analogously $(\bfllift^\gamma,\frlift^\delta)$ is not associative.

    Note that $\fllift^0=\bfllift^0$ and $\frlift^0=\bfrlift^0$. Therefore,
    Theorem \ref{thm:free_lift}, which lists all associative combinations of $\fllift^\gamma$ with $\frlift^\delta$ and of $\bfllift^\gamma$ with $\bfrlift^\delta$, indeed lists all associative universal lifts to the free product.  
 \end{proof}

\subsection{New multi-faced independence arising from universal lifts to the free product}\label{subsec:free_lift2}

We can construct new multi-faced universal products of states using the universal lifts $(\fllift^{\gamma},\frlift^{\delta})$ and $(\bfllift^{\gamma},\bfrlift^{\delta})$. Especially, we can construct a continuous deformation of bi-freeness which actually depends on two parameters in $\mathbb T$. 
First we investigate the basic 1-faced case. 

\begin{prop}\label{prop:free_single_face} For $(\gamma,\delta) \in \mathbb T^2\cup\{(0,0)\}$, let $\ufprod{\gamma}{\delta}$ and  $\ubprod{\gamma}{\delta}$ be the universal products of states associated with the universal lifts to the free product $(\fllift^\gamma,\frlift^\delta)$ and $(\bfllift^\gamma,\bfrlift^\delta)$, 
respectively. 
\begin{enumerate}[label=\rm(\arabic*)]
    \item $\ufprod{\gamma}{\delta}$ is the free product if and only if $\ubprod{\gamma}{\delta}$ is the free product if and only if $\gamma,\delta \in\mathbb T$. 
     \item $\ufprod{\gamma}{\delta}$ is the boolean product if and only if $\ubprod{\gamma}{\delta}$ is the boolean product if and only if $\gamma=\delta =0$. 
\end{enumerate}
\end{prop}
\begin{proof}
By symmetry, it suffices to discuss $\ufprod{\gamma}{\delta}$. 
\revise{We} know that $\ufprod{\gamma}{\delta}$ is one of the five universal products, and hence, checking some moments will suffice to identify it. \revise{(Similar to the tensor case, one could circumvent the use of Muraki's Theorem by using a single-faced variant of Proposition \ref{prop:free-free-parameters} below, namely proving that $\ufprod{\gamma}{\delta}=\ufprod{|\gamma|}{|\delta|}=\ubprod{\gamma}{\delta}$.)}

For later reference, let us take general four parameters $\gamma_1, \gamma_2,\delta_1,\delta_2 \in \mathbb T \cup \{0\}$. Let $X_1, X_2 \in \Lop(H)$ and $Y_1,Y_2 \in \Lop(G)$. We denote $\fllift^\gamma=\fllift^\gamma_{H,G}, \frlift^\delta=\frlift^\delta_{H,G}$ and $\langle T \rangle = \langle \Omega, T \Omega\rangle$ for shorter notation.
We will compute 
\begin{equation}\label{eq:vacuum_free_free}
\langle \fllift^{\gamma_1}(X_1)\frlift^{\delta_1}(Y_1)\fllift^{\gamma_2}(X_2)\frlift^{\delta_2}(Y_2) \rangle.  
\end{equation}
To compute this, we begin with
\begin{align}
    \fllift^{\gamma_2}(X_2)\frlift^{\delta_2}(Y_2)\Omega 
    &=\fllift^{\gamma_2}(X_2)(\langle Y_2\rangle \Omega + P_{\hat G}Y_2 \Omega) \label{eq:fllift-frlift-Omega}\\
    &=\langle X_2\rangle \langle Y_2\rangle\Omega + \langle Y_2\rangle P_{\hat H} X_2 \Omega + \langle \Omega, (X_2)_{\gamma_2} \Omega \rangle  P_{\hat G}Y_2 \Omega + (P_{\hat H} (X_2)_{\gamma_2} \Omega)\otimes(P_{\hat G}Y_2 \Omega)\notag \\
    &= \langle X_2\rangle \langle Y_2\rangle\Omega + \langle Y_2\rangle P_{\hat H} X_2 \Omega + |\gamma_2|\langle X_2\rangle  P_{\hat G}Y_2 \Omega + \gamma_2(P_{\hat H} X_2\Omega)\otimes(P_{\hat G}Y_2 \Omega),  \notag
\end{align}
and then we get 
\begin{align*}
    &\frlift^{\delta_1}(Y_1)\fllift^{\gamma_2}(X_2)\frlift^{\delta_2}(Y_2)\Omega 
    \\
    &\qquad=\langle Y_1 \rangle \langle X_2 \rangle \langle Y_2 \rangle \Omega 
    + \langle X_2 \rangle \langle Y_2 \rangle P_{\hat G} Y_1\Omega  
    + |\delta_1| \langle Y_1 \rangle \langle Y_2 \rangle P_{\hat H} X_2\Omega
    \\
    &\qquad\quad+\delta_1 \langle Y_2 \rangle (P_{\hat G} Y_1 \Omega)\otimes (P_{\hat H} X_2 \Omega)
    + |\gamma_2| \langle X_2 \rangle (\langle Y_1Y_2 \rangle - \langle Y_1\rangle \langle Y_2\rangle)\Omega
    +|\gamma_2|\langle X_2\rangle  P_{\hat G}Y_1P_{\hat G}Y_2 \Omega
    \\
    &\qquad\quad+ \gamma_2 |\delta_1| \langle Y_1\rangle (P_{\hat H} X_2 \Omega)\otimes(P_{\hat G} Y_2 \Omega) + \gamma_2\delta_1 (P_{\hat G} Y_1\Omega)\otimes (P_{\hat H} X_2 \Omega)\otimes(P_{\hat G} Y_2 \Omega).  
\end{align*}
Applying the operator $\fllift^{\gamma_1}(X_1)$ to the above, we see that the coefficient of $\Omega$ will be 
\begin{align}
    &\langle \fllift^{\gamma_1}(X_1)\frlift^{\delta_1}(Y_1)\fllift^{\gamma_2}(X_2)\frlift^{\delta_2}(Y_2) \rangle  
    \label{eq:free_free}
    \\
    &\qquad= \langle X_1 \rangle\langle Y_1 \rangle \langle X_2 \rangle \langle Y_2 \rangle 
    + |\delta_1| (\langle X_1X_2 \rangle - \langle X_1\rangle \langle X_2\rangle)\langle Y_1 \rangle \langle Y_2\rangle 
    \notag
    \\
    &\qquad\qquad\qquad\qquad\qquad\quad~+|\gamma_2|\langle X_1 \rangle \langle X_2 \rangle (\langle Y_1Y_2 \rangle - \langle Y_1\rangle \langle Y_2\rangle),  \notag 
\end{align}
which coincides with  
$
\langle X_1 \rangle\langle Y_1 \rangle \langle X_2 \rangle \langle Y_2 \rangle
$
 if  $\gamma_1=\gamma_2=\delta_1=\delta_2=0$ and with
\[
\langle X_1X_2 \rangle\langle Y_1\rangle \langle Y_2\rangle + \langle X_1 \rangle \langle X_2 \rangle \langle Y_1Y_2 \rangle  - \langle X_1 \rangle\langle Y_1 \rangle \langle X_2 \rangle \langle Y_2 \rangle
\]
if $\gamma_1 =\gamma_2 \in \mathbb T, \delta_1 = \delta_2 \in \mathbb T$. Among the five universal products,  only the boolean product yields the former result and only the free product yields the latter, so the proof is completed. Note that we only worked on operators on pre-Hilbert spaces but we can treat  $\ast$-algebras analogously to the proof of Proposition \ref{prop:classification_product_of_states}.  
\end{proof}

\begin{rmq}
  Note that there is a slight shortcut from \eqref{eq:fllift-frlift-Omega} to \eqref{eq:free_free}. Analogously to  \eqref{eq:fllift-frlift-Omega} (or by symmetry arguements) one obtains
  \[\frlift^{\delta_1}(Y_1^*)\fllift^{\gamma_1}(X_1^*)\Omega
    = \langle X_1^*\rangle \langle Y_1^*\rangle\Omega + \langle X_1^*\rangle P_{\hat G} Y_1^* \Omega + |\delta_1|\langle Y_1^*\rangle  P_{\hat H}X_1^* \Omega + \delta_1(P_{\hat G} Y_1^*\Omega)\otimes(P_{\hat H}X_1^* \Omega),  \notag\]
  and, therefore,
  \begin{align*}
    \MoveEqLeft\langle \fllift^{\gamma_1}(X_1)\frlift^{\delta_1}(Y_1)\fllift^{\gamma_2}(X_2^*)\frlift^{\delta_2}(Y_2^*) \rangle\\
    &=\langle\frlift^{\delta_1}(Y_1^*)\fllift^{\gamma_1}(X_1^*)\Omega,\fllift^{\gamma_2}(X_2)\frlift^{\delta_2}(Y_2)\Omega\rangle\\
    &= \langle X_1 \rangle\langle Y_1 \rangle \langle X_2 \rangle \langle Y_2 \rangle 
      + |\delta_1| (\langle X_1X_2 \rangle - \langle X_1\rangle \langle X_2\rangle)\langle Y_1 \rangle \langle Y_2\rangle
      + |\gamma_2|\langle X_1 \rangle \langle X_2 \rangle (\langle Y_1Y_2 \rangle - \langle Y_1\rangle \langle Y_2\rangle). 
  \end{align*}
  Similarly, having calculated $\frlift^{\delta_1}(Y_1)\fllift^{\gamma_2}(X_2)\frlift^{\delta_2}(Y_2)\Omega$, it would be quite easy to deduce mixed moments up to length 6.  
\end{rmq}

\revise{In our search for two-faced universal products of states coming from lifts to the free product,} 
we begin with the case where face 1 and face 2 are equipped with the universal lifts $(\fllift^{\gamma_1},\frlift^{\delta_1})$  and $(\fllift^{\gamma_2},\frlift^{\delta_2})$, respectively, where $(\gamma_1,\delta_1),(\gamma_2,\delta_2)$ are taken from $\mathbb T^2\cup\{(0,0)\}$. Let us denote by $\uffprod{\gamma_1}{\delta_1}{\gamma_2}{\delta_2}$ the associated universal product of states; the construction is analogous to the tensor case described in detail in Subsection \ref{subsec:tensor_lift2}. 

\begin{ex}\label{ex:free-free_moments} Let $\phi$ and $\psi$ be states on $A=A^{(1)}\sqcup A^{(2)}$ and $B=B^{(1)}\sqcup B^{(2)}$, respectively, and let $a_1\in A^{(1)}, a_2 \in A^{(2)}, b_1 \in B^{(1)}, b_2 \in B^{(2)}$. For $a \in A, b \in B$ let $\langle a \rangle =\phi(a), \langle b \rangle = \psi(b)$ for shorter notation.  
Equation \eqref{eq:free_free} implies that 
\begin{align}\label{eq:free-free_moment:a_1b_1a_2b_2}
    \MoveEqLeft \phi \uffprod{\gamma_1}{\delta_1}{\gamma_2}{\delta_2}\psi(a_1b_1a_2b_2)
    \\
    &= \langle a_1 \rangle\langle b_1 \rangle \langle a_2 \rangle \langle b_2 \rangle 
    + |\delta_1| (\langle a_1a_2 \rangle - \langle a_1\rangle \langle a_2\rangle)\langle b_1 \rangle \langle b_2\rangle
    +|\gamma_2|\langle a_1 \rangle \langle a_2 \rangle (\langle b_1b_2 \rangle - \langle b_1\rangle \langle b_2\rangle),   \notag
\end{align}
which only depends on $\{0,1\}$-parameters. On the other hand, we can compute
\begin{align}
    \MoveEqLeft\phi \uffprod{\gamma_1}{\delta_1}{\gamma_2}{\delta_2}\psi(b_2a_2a_1b_1)
    \label{eq:free-free_moment}\\
    &= \langle b_2 \rangle\langle a_2 \rangle \langle a_1 \rangle \langle b_1 \rangle 
    +\langle b_2 \rangle \langle b_1 \rangle (\langle a_2a_1 \rangle - \langle a_2\rangle \langle a_1\rangle)   
    + |\gamma_1\gamma_2|\langle a_2 \rangle \langle a_1\rangle (\langle b_2b_1 \rangle - \langle b_2\rangle \langle b_1\rangle)
    \notag \\
    &\quad+\gamma_1\overline{\gamma_2}(\langle b_2b_1 \rangle - \langle b_2\rangle \langle b_1\rangle)(\langle a_2a_1 \rangle - \langle a_2\rangle \langle a_1\rangle),  \notag
\end{align}
which depends on the continuous parameter $\gamma_1\overline{\gamma_2}$.
\end{ex}
\begin{ex}\label{ex:free-free} Let $(a_1,a_2)$ and $(b_1,b_2)$ be  $\uffprod{1}{1}{1}{1}$-independent pairs of elements in a $\ast$-probability space $(A,\phi)$ in the sense of Definition \ref{defi:multi-faced_independence}. This actually means that $\Alg{a_1,a_2}$ and $\Alg{b_1,b_2}$ are free, which can be confirmed from the canonical operator model.  
\end{ex}

In order to see how the universal product $\uffprod{\gamma_1}{\delta_1}{\gamma_2}{\delta_2}$ depends on the parameters, we need to examine mixed moments such as \revise{\eqref{eq:free-free_moment:a_1b_1a_2b_2} and \eqref{eq:free-free_moment}.} 
The way parameters are gained in the course of applying lifted operators to the vacuum vector can be summarized in Figure \ref{figure:free_free}.

\begin{figure}[t]
\[  
    \xymatrix{
   & \ldarrow{1} \hat H \rarrow{\text{$\delta_1$ or $\delta_2$}} \ar@(u,l)_{1,|\delta_1| \text{~or~} |\delta_2|}  &  \larrow{\text{$\overline{\delta_1}$ or $\overline{\delta_2}$}}  \hat G \otimes \hat H  \rarrow{\text{$\gamma_1$ or $\gamma_2$}}  \ar@(ur,ul)_{|\gamma_1|,|\gamma_2|,|\delta_1|\text{~or~}|\delta_2|~~~} & \larrow{\text{$\overline{\gamma_1}$ or $\overline{\gamma_2}$}} \hat H \otimes \hat G \otimes \hat H \rarrow{\text{$\delta_1$ or $\delta_2$}}  \ar@(ur,ul)_{|\gamma_1|,|\gamma_2|,|\delta_1|\text{~or~}|\delta_2|}  & \larrow{\text{$\overline{\delta_1}$ or $\overline{\delta_2}$}} \hat G \otimes \hat H \otimes \hat G \otimes \hat H  \rarrow{\text{$\gamma_1$ or $\gamma_2$}} \ar@(ur,ul)_{|\gamma_1|,|\gamma_2|,|\delta_1|\text{~or~}|\delta_2|}  &\larrow{\text{$\overline{\gamma_1}$ or $\overline{\gamma_2}$}} \cdots 
    \\
    \mathbb C \Omega \ruarrow{1} \rdarrow{1} \ar@(ul,dl)_{1}& & &  
    \\ 
    & \hat G \rarrow{\text{$\gamma_1$ or $\gamma_2$}} \luarrow{1}  \ar@(l,d)_{1,|\gamma_1|\text{~or~}|\gamma_2|}  &  \larrow{\text{$\overline{\gamma_1}$ or $\overline{\gamma_2}$}}  \hat H \otimes \hat G  \rarrow{\text{$\delta_1$ or $\delta_2$}}  \ar@(dl,dr)_{|\gamma_1|,|\gamma_2|,|\delta_1|\text{~or~}|\delta_2|~~~} & \larrow{\text{$\overline{\delta_1}$ or $\overline{\delta_2}$}} \hat G \otimes \hat H \otimes \hat G \rarrow{\text{$\gamma_1$ or $\gamma_2$}} \ar@(dl,dr)_{|\gamma_1|,|\gamma_2|,|\delta_1|\text{~or~}|\delta_2|}  & \larrow{\text{$\overline{\gamma_1}$ or $\overline{\gamma_2}$}} \hat H \otimes \hat G \otimes \hat H \otimes \hat G \rarrow{\text{$\delta_1$ or $\delta_2$}} \ar@(dl,dr)_{|\gamma_1|,|\gamma_2|,|\delta_1|\text{~or~}|\delta_2|} &\larrow{\text{$\overline{\delta_1}$ or $\overline{\delta_2}$}} \cdots 
    } 
\] 
\caption{The way operators lifted by $(\fllift^{\gamma_k},\frlift^{\delta_k}),k=1,2$ create parameters as they act on $H\AST G$}\label{figure:free_free}
\end{figure}

\revise{
The following mixed moments are especially useful later. 

  \begin{lm}\label{lem:mixed_moments} Let $(\gamma_i,\delta_i)\in \mathbb T^2\cup\{(0,0)\}, i=1,2$.
    Let $\phi$ and $\psi$ be states on $A=A^{(1)}\sqcup A^{(2)}$ and $B=B^{(1)}\sqcup B^{(2)}$, respectively, and let $a_1\in A^{(1)}, a_2 \in A^{(2)}, b_1 \in B^{(1)}, b_2 \in B^{(2)}$ such that $\phi(a_i)=\psi(b_i)=0$ for $i=1,2$ and $\phi(a_k^*a_\ell)=\psi(b_k^*b_\ell)=1$ for all $k,\ell\in\{1,2\}$.\footnote{For example, (with apologies for a slight conflict of notation) $A^{(i)}= B^{(i)}= \Lop(\mathbb C\Omega\oplus \mathbb C\xi)$ and $a_1=a_\xi^*\in A^{(1)}, b_1=a^*_\xi\in B^{(1)},a_2=a_\xi^*\in A^{(2)}, b_2=a_\xi^*\in B^{(2)}$ for a unit vector $\xi$. }
    Then, for all $i,k,\ell,j\in\{1,2\}$,
    \begin{equation}
      \phi \uffprod{\gamma_1}{\delta_1}{\gamma_2}{\delta_2}\psi (b_i^*a_k^*a_\ell b_j)=\gamma_\ell\overline{\gamma_k} \qquad \text{and} \qquad 
     \phi \uffprod{\gamma_1}{\delta_1}{\gamma_2}{\delta_2}\psi (a_i^*b_k^*b_\ell a_j)=\delta_\ell\overline{\delta_k}.\label{eq:free-free_moment3}
    \end{equation}
  \end{lm}

  \begin{proof}
    Direct verification. (One of these formulas is a  special case of \eqref{eq:free-free_moment}.) 
  \end{proof}

}


\revise{
  \begin{prop}\label{prop:free-free-parameters} Let $(\gamma_i,\delta_i), (\gamma_i',\delta_i') \in \mathbb T^2\cup\{(0,0)\}, i=1,2$. 
    The universal products of states $\uffprod{\gamma_1}{\delta_1}{\gamma_2}{\delta_2}$ and $\uffprod{\gamma_1'}{\delta_1'}{\gamma_2'}{\delta_2'}$ agree if and only if there are $\alpha,\beta\in\mathbb T$ with $\gamma_i'=\alpha \gamma_i$ and $\delta_i'=\beta\delta_i$ for $i=1,2$.
   \end{prop}

   \begin{proof}
     It easily follows from Figure \ref{figure:free_free} that each term of a mixed moment for $\uffprod{\gamma_1}{\delta_1}{\gamma_2}{\delta_2}$ is a monomial on the variables 
     \[
|\gamma_i|,|\delta_i|~(i=1,2), \gamma_1 \overline{\gamma_2}, \overline{\gamma_1} \gamma_2, \delta_1 \overline{\delta_2},  \overline{\delta_1} \delta_2. 
     \]
   This readily verifies the ``if''-part. The ``only if''-part follows from \eqref{eq:free-free_moment3}. Indeed, for $\gamma_1,\gamma_2\in\mathbb T\cup\{0\}$,  $\gamma_\ell\overline{\gamma_k}=\gamma_\ell'\overline{\gamma_k'}$ holds for all $k,\ell\in\{1,2\}$ (note that this implies $|\gamma_k|=|\gamma_k'|$).  The desired $\alpha$ can be defined as follows. 
     \begin{itemize}
     \item Case $\gamma_1=\gamma_1'=\gamma_2=\gamma_2'=0$: $\alpha\gamma_i=\gamma_i'$ ($i=1,2$) holds for an arbitrary $\alpha\in\mathbb T$,
     \item Case $\gamma_1=\gamma_1'=0$, $\gamma_2,\gamma_2'\in\mathbb T$: $\alpha\gamma_i=\gamma_i'$ ($i=1,2$) holds for $\alpha:=\gamma_2'\overline{\gamma_2}\in\mathbb T$,
     \item Case $\gamma_1,\gamma_1'\in\mathbb T$, $\gamma_2=\gamma_2'=0$: $\alpha\gamma_i=\gamma_i'$ ($i=1,2$) holds for $\alpha:=\gamma_1'\overline{\gamma_1}\in\mathbb T$,
     \item Case  $\gamma_1,\gamma_1',\gamma_2,\gamma_2'\in\mathbb T$: $\alpha\gamma_i=\gamma_i'$ ($i=1,2$) holds for $\alpha:=\gamma_1'\overline{\gamma_1}=\gamma_2'\overline{\gamma_2}\in\mathbb T$.
     \end{itemize}
     Finding $\beta\in\mathbb T$ with $\delta_i'=\beta\delta_i$ for $i=1,2$ works analogously, using the second equality of \eqref{eq:free-free_moment3}.
   \end{proof}

 }



Next, we investigate the case where face 1 and face 2 are equipped with the universal lifts $(\fllift^{\gamma_1},\frlift^{\delta_1})$ and $(\bfllift^{\gamma_2},\bfrlift^{\delta_2})$,  respectively, with general parameters $ (\gamma_1,\delta_1),(\gamma_2,\delta_2) \in \mathbb T^2\cup\{(0,0)\}$, with $\ubfprod{\gamma_1}{\delta_1}{\gamma_2}{\delta_2}$ the associated universal product of states. This is identical to the bi-free product when all the parameters are equal to one. In fact, this is the bi-free product if and only if $\gamma_1 = \delta_2 \in \mathbb T$ and $\delta_1 = \gamma_2 \in \mathbb T$
. 

\begin{ex}\label{ex:mixed_moments_bifree} For two states $\phi$ on $A$ and $\psi$ on $B$, take $\ast$\hyp{}representations $\pi\colon A \to \Lop(H)$ and $\sigma\colon B \to \Lop(G)$ such that $\langle a \rangle\colonequals \phi(a) = \langle \Omega, \pi(a)\Omega\rangle$ and $\langle b \rangle\colonequals \psi(b) = \langle \Omega, \sigma(b)\Omega\rangle$. As before, we employ the notation $X_k=\pi(a_k), Y_k=\sigma(b_k)$, and for the lifts,  $\bfllift^{\gamma}=\bfllift^{\gamma}_{H,K}$ and so on. 
We begin with
\begin{align*}
    \bfllift^{\gamma_2}(X_2)\bfrlift^{\delta_2}(Y_2)\Omega 
    &=\bfllift^{\gamma_2}(X_2)(\langle Y_2\rangle \Omega + P_{\hat G}Y_2 \Omega) \\
    &=\langle X_2\rangle \langle Y_2\rangle\Omega + \langle Y_2\rangle P_{\hat H} X_2 \Omega + \langle \Omega, (X_2)_{\gamma_2} \Omega \rangle  P_{\hat G}Y_2 \Omega + (P_{\hat G}Y_2 \Omega)\otimes (P_{\hat H} (X_2)_{\gamma_2} \Omega)\\
    &= \langle X_2\rangle \langle Y_2\rangle\Omega + \langle Y_2\rangle P_{\hat H} X_2 \Omega + |\gamma_2|\langle X_2\rangle  P_{\hat G}Y_2 \Omega + \gamma_2(P_{\hat G}Y_2 \Omega)\otimes (P_{\hat H} X_2 \Omega).  
\end{align*}
This implies that 
\begin{align*}
    &\frlift^{\delta_1}(Y_1)\bfllift^{\gamma_2}(X_2)\bfrlift^{\delta_2}(Y_2)\Omega 
    \\
    &\qquad=\langle Y_1 \rangle \langle X_2 \rangle \langle Y_2 \rangle \Omega 
    + \langle X_2 \rangle \langle Y_2 \rangle P_{\hat G} Y_1\Omega  
    + |\delta_1| \langle Y_1 \rangle \langle Y_2 \rangle P_{\hat H} X_2\Omega
    \\
    &\qquad\quad+\delta_1 \langle Y_2 \rangle (P_{\hat G} Y_1 \Omega)\otimes (P_{\hat H} X_2 \Omega)
    + |\gamma_2| \langle X_2 \rangle (\langle Y_1Y_2 \rangle - \langle Y_1\rangle \langle Y_2\rangle)\Omega
    +|\gamma_2|\langle X_2\rangle  P_{\hat G}Y_1P_{\hat G}Y_2 \Omega
    \\
    &\qquad\quad+ \gamma_2 \overline{\delta_1} \langle \Omega, Y_1P_{\hat G}Y_2\Omega \rangle P_{\hat H} X_2 \Omega + \gamma_2
    |\delta_1| (P_{\hat G}Y_1P_{\hat G} Y_2\Omega)\otimes (P_{\hat H} X_2 \Omega).  
\end{align*}
Applying the operator $\fllift^{\gamma_1}(X_1)$ and collecting the coefficients of $\Omega$ we obtain
\begin{align}
\phi\ubfprod{\gamma_1}{\delta_1}{\gamma_2}{\delta_2}\psi(a_1b_1a_2b_2)
&=\langle a_1\rangle \langle b_1\rangle \langle a_2\rangle \langle b_2\rangle 
+ |\delta_1|(\langle a_1a_2\rangle - \langle a_1\rangle \langle a_2 \rangle ) \langle b_1\rangle \langle b_2\rangle 
\label{eq:bifree_moment1}\\
&\quad +|\gamma_2|\langle a_1\rangle \langle a_2\rangle( \langle b_1b_2\rangle - \langle b_1\rangle \langle b_2 \rangle ) 
\notag\\
&\quad+ \gamma_2 \overline{\delta_1} ( \langle a_1a_2\rangle - \langle a_1\rangle \langle a_2 \rangle )( \langle b_1 b_2\rangle - \langle b_1\rangle \langle b_2 \rangle).\notag 
\end{align}
\end{ex}

The way lifted operators create parameters in the process of computing mixed moments is visualized in Figure \ref{figure:bifree}. 

\revise{
To parametrize the universal products of states $\ubfprod{\gamma_1}{\delta_1}{\gamma_2}{\delta_2}$ with a minimal set of parameters, the following mixed moments are useful. 


}

 \revise{\begin{lm} \label{lem:bi-free_moment3}
    In the setting of Lemma \ref{lem:mixed_moments}, for all $i,k,\ell,j\in\{1,2\}$ with $k \ne \ell$, 
    \begin{align}
      \phi \ubfprod{\gamma_1}{\delta_1}{\gamma_2}{\delta_2}\psi (a_i^*b_k^*a_\ell b_j) &=\gamma_\ell\overline{\delta_k}, \\  
    \phi \ubfprod{\gamma_1}{\delta_1}{\gamma_2}{\delta_2}\psi (b_i^*a_k^* a_k b_j)=|\gamma_k|^2 \qquad &\text{and}\qquad 
    \phi \ubfprod{\gamma_1}{\delta_1}{\gamma_2}{\delta_2}\psi (a_i^*b_k^* b_k a_j)=|\delta_k|^2. 
    \end{align}
  \end{lm}
   \begin{proof}
    Direct verification. (One of these formulas is a special case of \eqref{eq:bifree_moment1}.)
  \end{proof}
  }

\begin{figure}[tb]
    \[\xymatrix{
   & \ldarrow{1} \hat H \rarrow{\delta_1} \ar@(u,l)_{1,|\delta_1| \text{~or~} |\delta_2|} \rddarrow{\delta_2} &  \larrow{\overline{\delta_1}}  \hat G \otimes \hat H  \rarrow{\gamma_1}  \ar@(ur,ul)_{|\gamma_1|,|\gamma_2|,|\delta_1|\text{~or~}|\delta_2|~~~} \rddarrow{\delta_2} \lddarrow{\overline{\gamma_2}} & \larrow{\overline{\gamma_1}} \hat H \otimes \hat G \otimes \hat H \rarrow{\delta_1}  \ar@(ur,ul)_{|\gamma_1|,|\gamma_2|,|\delta_1|\text{~or~}|\delta_2|} \rddarrow{\delta_2}\lddarrow{\overline{\gamma_2}} & \larrow{\overline{\delta_1}} \hat G \otimes \hat H \otimes \hat G \otimes \hat H  \rarrow{\gamma_1} \ar@(ur,ul)_{|\gamma_1|,|\gamma_2|,|\delta_1|\text{~or~}|\delta_2|} \rddarrow{\delta_2} \lddarrow{\overline{\gamma_2}} &\larrow{\overline{\gamma_1}} \lddarrow{\overline{\gamma_2}}\cdots 
    \\
    \mathbb C \Omega \ruarrow{1} \rdarrow{1} \ar@(ul,dl)_{1}& & &  
    \\ 
    & \hat G \rarrow{\gamma_1} \luarrow{1}  \ar@(l,d)_{1,|\gamma_1|\text{~or~}|\gamma_2|} \ruuarrow{\gamma_2} &  \larrow{\overline{\gamma_1}}  \hat H \otimes \hat G  \rarrow{\delta_1}  \ar@(dl,dr)_{|\gamma_1|,|\gamma_2|,|\delta_1|\text{~or~}|\delta_2|~~~} \ruuarrow{\gamma_2} \luuarrow{\overline{\delta_2}}& \larrow{\overline{\delta_1}} \hat G \otimes \hat H \otimes \hat G \rarrow{\gamma_1} \ar@(dl,dr)_{|\gamma_1|,|\gamma_2|,|\delta_1|\text{~or~}|\delta_2|} \ruuarrow{\gamma_2} \luuarrow{\overline{\delta_2}}& \larrow{\overline{\gamma_1}} \hat H \otimes \hat G \otimes \hat H \otimes \hat G \rarrow{\delta_1} \ar@(dl,dr)_{|\gamma_1|,|\gamma_2|,|\delta_1|\text{~or~}|\delta_2|} \ruuarrow{\gamma_2}\luuarrow{\overline{\delta_2}} &\larrow{\overline{\delta_1}} \luuarrow{\overline{\delta_2}}\cdots 
    } \]
    \caption{The way operators lifted by $(\fllift^{\gamma_1},\frlift^{\delta_1}),(\bfllift^{\gamma_2},\bfrlift^{\delta_2})$ create parameters, acting on $H\AST G$}\label{figure:bifree}
\end{figure}



\revise{
  \begin{prop}\label{prop:bi-free-parameters} Let $(\gamma_i,\delta_i), (\gamma_i',\delta_i') \in \mathbb T^2\cup\{(0,0)\}, i=1,2$. 
     The universal products of states $\ubfprod{\gamma_1}{\delta_1}{\gamma_2}{\delta_2}$ and $\ubfprod{\gamma_1'}{\delta_1'}{\gamma_2'}{\delta_2'}$ agree if and only if there are $\alpha,\beta\in\mathbb T$ with $\gamma_1'=\alpha \gamma_1$, $\delta_2'=\alpha \delta_2$, $\gamma_2'=\beta \gamma_2$,  and $\delta_1'=\beta\delta_1$. 
   \end{prop}
   }


   \begin{proof}
     \revise{The ``only if''-part follows from Lemma \ref{lem:bi-free_moment3} analogously to the proof of Proposition \ref{prop:free-free-parameters}. For the ``if''-part, it suffices to prove that each term of a mixed moment with respect to $\ubfprod{\gamma_1}{\delta_1}{\gamma_2}{\delta_2}$ is a monomial on the variables 
   \[
   |\gamma_i|, |\delta_i|~(i=1,2), \gamma_1 \overline{\delta_2}, \overline{\gamma_1} \delta_2, \gamma_2 \overline{\delta_1}, \overline{\gamma_2} \delta_1. 
   \]} What needs to be proved is that, for any directed walk on Figure \ref{figure:bifree} (a digraph having pre-Hilbert spaces as vertices and complex weights on edges) whose initial and terminal vertices are both $\mathbb C \Omega$, the product of all the weights on the walk (called the total weight below) must be such a monomial. We will actually show that the total weight of every cycle must be of that form.
  
  First, observe that the total weight of every cycle of length four is such a monomial; the weights repeat periodically, so there are only finitely many cycles to check, namely those lying completely within Figure \ref{figure:bifree}. Cycles of length less than four are  easily seen to have total weight 1.
  
  Suppose that the claim is true for all cycles of length less than $k$, and let us take a cycle $w=(V_0,\ldots, V_k)$, $V_0=V_k$ of length $k>4$. For a vertex $V$, we denote by $d(V)$ the number of tensor factors, i.e. $d(\mathbb C \Omega)=0, d(\hat G\otimes \hat H)=2$ etc., and call $d(V)$ the \emph{degree} of $V$. Note that there is an edge from $V$ to $V'$ only if $|d(V)-d(V')| \leq 1$. If $w$ contains a loop, removing the \revise{loop does not change the total weight up to a factor $|\gamma_i|$ or $|\delta_i|$, and we are done.} If $w=(V_0,\ldots, V_k)$ contains no loop, let $V_h$ be a vertex of highest degree in $w$ and consider the subwalk $(V_{h-2},V_{h-1},V_h,V_{h+1})$ of $w$ (assume w.l.o.g.\ that $2\le h\le k-1$ for notational convenience). It follows that \revise{$d(V_{h+1})=d(V_h)-1=d(V_{h-1})=d(V_{h-2})\pm 1$} and, therefore, there are edges from $V_{h-2}$ to $V_{h+1}$ and back (with inverse weights). The total weight of $w$ is the same as the total weight of the cycle of length $k-2$ $w':=(V_0,\ldots, V_{h-2},V_{h+1}, \ldots, V_k)$ multiplied with the total weight of the cycle of length four $(V_{h+1}, V_{h-2},V_{h-1},V_{h},V_{h+1})$ \revise{up to a factor $|\gamma_i|^2$ or $|\delta_i|^2$.} By the induction hypothesis, the total weight of $w$ is of the desired form. 
\end{proof}

\revise{
The universal products $\ubfbfprod{\gamma_1}{\delta_1}{\gamma_2}{\delta_2}$ and $\ubffprod{\gamma_1}{\delta_1}{\gamma_2}{\delta_2}$, defined analogously, are easily seen to agree with the ones we discussed so far, $\uffprod{\gamma_1}{\delta_1}{\gamma_2}{\delta_2}=\ubfbfprod{\gamma_1}{\delta_1}{\gamma_2}{\delta_2}$ and $\ubfprod{\gamma_1}{\delta_1}{\gamma_2}{\delta_2}=\ubffprod{\gamma_1}{\delta_1}{\gamma_2}{\delta_2}$. This follows by symmetry of the construction, or more precisely by reversing the order of tensor components in the free product of pre-Hilbert spaces. 
Also one can see from Figures \ref{figure:free_free} and \ref{figure:bifree} that $\uffprod{0}{0}{1}{1}=\ubfprod{0}{0}{1}{1}$ and $\uffprod{1}{1}{0}{0}=\ubfprod{1}{1}{0}{0}$ because the two figures are ``isomorphic'' as soon as the edges with zero weights are deleted.  

We combine this observation with Propositions \ref{prop:free-free-parameters} and \ref{prop:bi-free-parameters} to obtain a completely parametrized family which contains all universal products of states coming from lifts to the free product. Note that the free-boolean product is the one introduced by Liu \cite{Liu19} and the boolean-boolean product is the same trivial two-faced extension of the boolean product as in Table \ref{table:tensor_lift-products}. 
 
\begin{thm}  \label{thm:classification_universal_product_states_free} 
Every element of the set
\[
\left\{{\uffprod{\gamma_1}{\delta_1}{\gamma_2}{\delta_2}},{ \ubfprod{\gamma_1}{\delta_1}{\gamma_2}{\delta_2}}, {\ubffprod{\gamma_1}{\delta_1}{\gamma_2}{\delta_2}},{\ubfbfprod{\gamma_1}{\delta_1}{\gamma_2}{\delta_2}}: (\gamma_1,\delta_1), (\gamma_2,\delta_2) \in \mathbb T^2\cup\{(0,0)\}\right\}
\]
coincides with exactly one of the following list. 
\begin{enumerate}[label=\rm(\arabic*)]

\item $\deffreeprod{\zeta}{\theta}:=\uffprod{\zeta}{\theta}{1}{1} =\ubfbfprod{\zeta}{\theta}{1}{1}~(\zeta,\theta\in \mathbb T)$ (a deformation of the free-free product),

\item $\defbifreeprod{\zeta}{\theta}:=\ubfprod{\zeta}{\theta}{1}{1}= \ubffprod{\zeta}{\theta}{1}{1}~(\zeta,\theta\in \mathbb T)$  (a deformation of the bi-free product), 

\item $\freeboolprod:=\uffprod{1}{1}{0}{0}=\ubfprod{1}{1}{0}{0}= \ubffprod{1}{1}{0}{0}=\ubfbfprod{1}{1}{0}{0}$ (the free-boolean product),

\item $\boolfreeprod:=\uffprod{0}{0}{1}{1}=\ubfprod{0}{0}{1}{1}= \ubffprod{0}{0}{1}{1}=\ubfbfprod{0}{0}{1}{1}$ (the boolean-free product), 

\item $\boolboolprod:=\uffprod{0}{0}{0}{0}$ (the boolean-boolean product), 

\end{enumerate}

\end{thm} 
Table \ref{table:free_lift-products} summarizes which combination of lifts gives rise to which universal product in the list of Theorem \ref{thm:classification_universal_product_states_free}.  

 }

\begin{table}[tb]\caption{Two-faced universal products of states from lifts to the free product.}
  \label{table:free_lift-products}
  \centering
  \begin{tabular}{|c|c|c|c|} \hline
   \multirow{2}{*}{\diagbox[width=3cm]{face 1}{face 2}} & left free  & right free  & boolean \\ 
    & $\Asterisk$ $(\fllift^{\alpha_2},\frlift^{\beta_2})$ & $\Asterisk$ $(\bfllift^{\alpha_2},\bfrlift^{\beta_2})$& $\bool$ $(\fllift^0,\frlift^0)$\\ \hline
     left free & \multirow{2}{*}{$\deffreeprod{\zeta}{\theta}$}  & \multirow{2}{*}{$\defbifreeprod{\zeta}{\theta}$}  & \multirow{2}{*}{$\freeboolprod$}  \\
   $\Asterisk$ $(\fllift^{\alpha_1},\frlift^{\beta_1})$ & &   &  \\ \hline
     right free & \multirow{2}{*}{$\defbifreeprod{\zeta}{\theta}$}  & \multirow{2}{*}{$\deffreeprod{\zeta}{\theta}$}  & \multirow{2}{*}{$\freeboolprod$}  \\
    $\Asterisk$ $(\bfllift^{\alpha_1},\bfrlift^{\beta_1})$
                                                          & &   & \\ \hline
    boolean & \multirow{2}{*}{$\boolfreeprod$}  & \multirow{2}{*}{$\boolfreeprod$}  &\multirow{2}{*}{$\boolboolprod$}  \\
    $\bool$ $(\fllift^0,\frlift^0)$ & &   & \\ \hline 
  \end{tabular}\\[.5em]
  ($\zeta,\theta,\alpha_1,\alpha_2,\beta_1,\beta_2\in\mathbb T, \zeta=\alpha_1\overline{\alpha_2}, \theta=\beta_1\overline{\beta_2}$) 
   
\end{table}

\revise{We can classify the symmetric two-faced universal products listed in Theorem \ref{thm:classification_universal_product_states_free}. 
}

\revise{
\begin{lm}\label{lm:symmetry-free}
For all $(\gamma_1,\delta_1), (\gamma_2,\delta_2) \in \mathbb T^2\cup\{0,0\}$ and any two states $\phi$ on $A$ and $\psi$ on $B$,  we have $\phi \uffprod{\gamma_1}{\delta_1}{\gamma_2}{\delta_2} \psi =\psi \uffprod{\delta_1}{\gamma_1}{\delta_2}{\gamma_2} \phi$ and $\phi \ubfprod{\gamma_1}{\delta_1}{\gamma_2}{\delta_2} \psi =\psi \ubfprod{\delta_1}{\gamma_1}{\delta_2}{\gamma_2} \phi$. 
\end{lm}

\begin{proof}
  The proof works very similar to that of Lemma \ref{lm:symmetry} in the tensor case. If we identify $H*G$ with $G*H$, we find $\fllift_{H,G}^\gamma=\frlift_{G,H}^\gamma$ and $\bfllift_{H,G}^\gamma=\bfrlift_{G,H}^\gamma$ for all pre-Hilbert spaces $H,G$ and all $\gamma\in\mathbb T\cup \{0\}$. For the second claim also note that $\ubfprod{\delta_1}{\gamma_1}{\delta_2}{\gamma_2}=\ubffprod{\delta_1}{\gamma_1}{\delta_2}{\gamma_2}$. We leave the details to the reader.
\end{proof}
}

\begin{prop}\label{prop:symmetry_free}
Among the two-faced universal products listed in Theorem \ref{thm:classification_universal_product_states_free}, the symmetric ones are exactly $\deffreeprod{\zeta}{\zeta}, \defbifreeprod{\zeta}{\zeta}~(\zeta\in\mathbb T),\freeboolprod,\boolfreeprod,\boolboolprod$. 
\end{prop}

\begin{proof}
 \revise{ By Lemma \ref{lm:symmetry-free}, the symmetry of the product $\uffprod{\gamma_1}{\delta_1}{\gamma_2}{\delta_2}$ is equivalent to the condition $\uffprod{\gamma_1}{\delta_1}{\gamma_2}{\delta_2} =\uffprod{\delta_1}{\gamma_1}{\delta_2}{\gamma_2}$. Hence, by Proposition \ref{prop:free-free-parameters}, $\deffreeprod{\zeta}{\theta}$ is symmetric if and only if $\zeta=\theta$. Similarly, the symmetry of the product $\defbifreeprod{\zeta}{\theta}$ is equivalent to the condition  $\zeta=\theta$ by Proposition \ref{prop:bi-free-parameters}, and  the other three products are symmetric.}  
\end{proof}

\begin{rmq}\label{rmq:convolution} The reader might wonder whether the two-faced universal products of states 
\begin{align}
& \left\{\uttprod{\gamma_1}{\delta_1}{\gamma_2}{\delta_2}: (\gamma_1,\delta_1), (\gamma_2,\delta_2) \in J_\otimes \right\}, \label{eq:two-faced1} \\
& \left\{\uffprod{\gamma_1}{\delta_1}{\gamma_2}{\delta_2}: (\gamma_1,\delta_1), (\gamma_2,\delta_2) \in \mathbb T^2 \cup \{(0,0)\} \right\} \quad \text{and} \\
& \left\{\ubfprod{\gamma_1}{\delta_1}{\gamma_2}{\delta_2}: (\gamma_1,\delta_1), (\gamma_2,\delta_2) \in \mathbb T^2 \cup \{(0,0)\} \right\}   \label{eq:two-faced3}
\end{align}
 give rise to convolutions of probability measures on $\R^2$. To make this question clearer, note that a pair $\mathbf a=(a_1,a_2)$ of commuting bounded self-adjoint operators in a $W^\ast$-probability space $(A,\phi)$ associates the Borel probability measure $\mu_{\mathbf a}$ on $\R^2$ defined by $\mu_{\mathbf a}(\cdot) = \phi(E_{\mathbf a}(\cdot))$, where $\{E_{\mathbf a}(B)\}_{B \in \mathcal B(\R^2)}$ is the spectral decomposition of $\mathbf a$. 
If we have two $\odot$-independent pairs $\mathbf a=(a_1,a_2)$ and $\mathbf b=(b_1,b_2)$ of bounded self-adjoint operators such that $[a_1,a_2]=[a_1,b_2]=[b_1,a_2]=[b_1,b_2]=0$, then $\mathbf a + \mathbf b$ also consists of commuting bounded self-adjoint operators, so that the probability measure $\mu_{\mathbf a + \mathbf b}$ can be defined and called the $\odot$-convolution of $\mu_{\mathbf a}$ and $\mu_{\mathbf b}$. The convolution thus defined has been studied in the literature for the two cases $\odot = \uttprod{1}{1}{1}{1}$ (the standard tensor-tensor product) and $\odot = \ubfprod{1}{1}{1}{1}$ (the bi-free product), see e.g.\ \cite{MS01,Sato13} for the former and \cite{BBGS18,GuHuangMingo16,HuangWang16,HasebeHuangWang18} for the latter.   

In fact, such a definition of convolution for a reasonably wide class of probability measures works only for the above mentioned cases $\odot = \uttprod{1}{1}{1}{1}$ and $\odot = \ubfprod{1}{1}{1}{1}$ among the family  \eqref{eq:two-faced1}--\eqref{eq:two-faced3}. For example, this can be confirmed for the family $\uttprod{\gamma_1}{\delta_1}{\gamma_2}{\delta_2}$ from Example \ref{ex:tensor_mixed_moments} where we have obtained  
\begin{align}
&\langle a_1b_1a_2b_2\rangle \label{eq:commutativityA} \\
&=\langle a_1 \rangle \langle b_1\rangle  \langle a_2\rangle \langle b_2\rangle + |\delta_1| \bigl( \langle a_1a_2\rangle-\langle a_1\rangle\langle a_2\rangle\bigr)\langle b_1\rangle \langle b_2\rangle+|\gamma_2| \langle a_1\rangle\langle a_2\rangle \bigl(\langle b_1b_2\rangle-\langle b_1\rangle\langle b_2\rangle\bigr) \notag \\
&\quad
+ \gamma_2\overline{\delta_1} \bigl( \langle a_1a_2\rangle-\langle a_1\rangle\langle a_2\rangle\bigr)\bigl( \langle b_1b_2\rangle-\langle b_1\rangle\langle b_2\rangle\bigr).  \notag
\end{align}
If the commutativity $[b_1,a_2]=0$ holds, then some calculations show that the LHS of \eqref{eq:commutativityA} coincides with 
$\langle a_1 a_2 \rangle \langle b_1b_2\rangle$, which will imply $\gamma_2\overline{\delta_1}=1$ under a moderate condition (e.g.\ the random vectors $\mathbf a$ and $\mathbf b$ have mean $(0,0)$ and non-vanishing covariances). A similar argument by symmetry also implies $\gamma_1\overline{\delta_2}=1$. Then we conclude $\gamma_1=\delta_1 = \gamma_2 = \delta_2 \in\mathbb T$, which implies $\uttprod{\gamma_1}{\delta_1}{\gamma_2}{\delta_2}= \uttprod{1}{1}{1}{1}$ thanks to 
\revise{\eqref{eq:rotation2}}. Similar reasonings apply to $\uffprod{\gamma_1}{\delta_1}{\gamma_2}{\delta_2}$ and $\ubfprod{\gamma_1}{\delta_1}{\gamma_2}{\delta_2}$. 
\end{rmq}


\section{Conclusion and future research}
\revise{In this paper, we were}
able to prove the scheme \eqref{new_scheme} connecting the concepts of universal lift, universal product of representations and universal product of states. Moreover, the universal lifts on both the tensor product and the free product were entirely classified. By doing that, and through the connection with universal products of states exhibited by scheme \eqref{new_scheme}, new multi-faced products of states were found.

Nevertheless, important questions remain open. The most obvious one is that our result is in no way a classification of all multi-faced products of states. The problem is two-fold: first, scheme \eqref{new_scheme} does not provide an equivalence between universal products of representations and universal products of states; second, even if we had such an equivalence, we would need to know if there are other monoidal products to consider and, if so, to classify universal lifts for those products.

The possibility of completing scheme \eqref{new_scheme} into a complete equivalence is not clear. Given a universal product of states, one should be able to build both a monoidal product and a universal product of representations on this monoidal product. It seems not clear to know how this monoidal product should be constructed from the universal product of states.

Thus, we still need a complete classification of multi-faced independences, which would be a true generalization of Muraki's result for single-faced independences. This article can, nevertheless, be seen as a first attempt in this direction, clearing the way for further research in this area.

\section*{Acknowledgements} The authors are grateful to Philipp Var\v{s}o for useful discussions.

\bibliographystyle{myalphaurl}
\bibliography{sample}

\begin{thebibliography}{BBGS18}

\bibitem[Avi82]{Avitzour82}
D.~Avitzour.
\newblock Free products of {$C^{\ast} $}-algebras.
\newblock {\em Trans. Amer. Math. Soc.}, 271(2):423--435, 1982.
\newblock \href {https://doi.org/10.2307/1998890} {\path{doi:10.2307/1998890}}.

\bibitem[BBGS18]{BBGS18}
S.~T. Belinschi, H.~Bercovici, Y.~Gu, and P.~Skoufranis.
\newblock Analytic subordination for bi-free convolution.
\newblock {\em J. Funct. Anal.}, 275(4):926--966, 2018.
\newblock \href {https://doi.org/10.1016/j.jfa.2018.03.003}
  {\path{doi:10.1016/j.jfa.2018.03.003}}.

\bibitem[BGS02]{BenGhorbalSchurmann}
A.~Ben~Ghorbal and M.~Sch{\"u}rmann.
\newblock Non-commutative notions of stochastic independence.
\newblock {\em Math. Proc. Cambridge Philos. Soc.}, 133(3):531--561, 2002.
\newblock \href {https://doi.org/10.1017/S0305004102006072}
  {\path{doi:10.1017/S0305004102006072}}.

\bibitem[BGS05]{BenGhorbalSchurmann05}
A.~Ben~Ghorbal and M.~Sch{\"u}rmann.
\newblock Quantum {L}\'evy processes on dual groups.
\newblock {\em Math. Z.}, 251(1):147--165, 2005.
\newblock \href {https://doi.org/10.1007/s00209-005-0793-x}
  {\path{doi:10.1007/s00209-005-0793-x}}.

\bibitem[Bo{\.z}86]{Bozejko86}
M.~Bo{\.z}ejko.
\newblock Positive definite functions on the free group and the noncommutative
  {R}iesz product.
\newblock {\em Boll. Un. Mat. Ital. A (6)}, 5(1):13--21, 1986.

\bibitem[BLS96]{BozejkoLeinertSpeicher}
M.~Bo\.{z}ejko, M.~Leinert, and R.~Speicher.
\newblock Convolution and limit theorems for conditionally free random
  variables.
\newblock {\em Pacific J. Math.}, 175(2):357--388, 1996.
\newblock Available from \url{http://projecteuclid.org/euclid.pjm/1102353149}.

\bibitem[Fra01]{Franz01}
U.~Franz.
\newblock Monotone independence is associative.
\newblock {\em Infin. Dimens. Anal. Quantum Probab. Relat. Top.},
  4(3):401--407, 2001.
\newblock \href {https://doi.org/10.1142/S0219025701000565}
  {\path{doi:10.1142/S0219025701000565}}.

\bibitem[Fra06]{Franz06}
U.~Franz.
\newblock L\'evy processes on quantum groups and dual groups.
\newblock In {\em Quantum independent increment processes. {II}}, volume 1866
  of {\em Lecture Notes in Math.}, pages 161--257. Springer, Berlin, 2006.
\newblock \href {https://doi.org/10.1007/11376637_3}
  {\path{doi:10.1007/11376637_3}}.

\bibitem[Ger17]{Gerhold17p}
M.~Gerhold.
\newblock Bimonotone brownian motion.
\newblock Pre\-print, t.a. in Proceedings of QP 38 (38th International
  Conference on Quantum Probability and Related Topics, Tokyo, Japan, October
  1-7, 2017), 2017.
\newblock \href {https://arxiv.org/abs/1708.03510} {\path{arXiv:1708.03510}}.

\bibitem[Ger21]{Gerhold21p}
M.~Gerhold.
\newblock Schoenberg correspondence for multifaced independence.
\newblock Pre\-print, 2021.
\newblock \href {https://arxiv.org/abs/2104.02985} {\path{arXiv:2104.02985}}.

\bibitem[GL15]{GerholdLachs15}
M.~Gerhold and S.~Lachs.
\newblock Classification and {GNS}-construction for general universal products.
\newblock {\em Infin. Dimens. Anal. Quantum Probab. Relat. Top.},
  18(1):1550004, 29 pages, 2015.
\newblock \href {https://doi.org/10.1142/S0219025715500046}
  {\path{doi:10.1142/S0219025715500046}}.

\bibitem[GLS22]{GerholdLachsSchurmann22}
M.~Gerhold, S.~Lachs, and M.~Sch\"{u}rmann.
\newblock Categorial independence and {L}\'{e}vy processes.
\newblock {\em SIGMA Symmetry Integrability Geom. Methods Appl.}, 18(075):27
  pages, 2022.
\newblock \href {https://doi.org/10.3842/SIGMA.2022.075}
  {\path{doi:10.3842/SIGMA.2022.075}}.

\bibitem[GV23]{GerholdVarso23p}
M.~Gerhold and P.~Var\v{s}o.
\newblock Towards a classification of multi-faced independences: a
  combinatorial approach.
\newblock Pre\-print, 2023.
\newblock \href {https://arxiv.org/abs/2301.01816} {\path{arXiv:2301.01816}}.

\bibitem[GHS20]{GuHasebeSkoufranis20}
Y.~Gu, T.~Hasebe, and P.~Skoufranis.
\newblock Bi-monotonic independence for pairs of algebras.
\newblock {\em J. Theoret. Probab.}, 33(1):533--566, 2020.
\newblock \href {https://doi.org/10.1007/s10959-019-00884-2}
  {\path{doi:10.1007/s10959-019-00884-2}}.

\bibitem[GHM16]{GuHuangMingo16}
Y.~Gu, H.-W. Huang, and J.~A. Mingo.
\newblock An analogue of the {L}\'{e}vy-{H}in\v{c}in formula for bi-free
  infinitely divisible distributions.
\newblock {\em Indiana Univ. Math. J.}, 65(5):1795--1831, 2016.
\newblock \href {https://doi.org/10.1512/iumj.2016.65.5911}
  {\path{doi:10.1512/iumj.2016.65.5911}}.

\bibitem[GS19]{GuSkoufranis}
Y.~Gu and P.~Skoufranis.
\newblock Bi-{B}oolean independence for pairs of algebras.
\newblock {\em Complex Anal. Oper. Theory}, 13(7):3023--3089, 2019.
\newblock \href {https://doi.org/10.1007/s11785-017-0750-9}
  {\path{doi:10.1007/s11785-017-0750-9}}.

\bibitem[HHW18]{HasebeHuangWang18}
T.~Hasebe, H.-W. Huang, and J.-C. Wang.
\newblock Limit theorems in bi-free probability theory.
\newblock {\em Probab. Theory Related Fields}, 172(3-4):1081--1119, 2018.
\newblock \href {https://doi.org/10.1007/s00440-017-0825-6}
  {\path{doi:10.1007/s00440-017-0825-6}}.

\bibitem[HL17]{HasebeLehner17}
T.~Hasebe and F.~Lehner.
\newblock Cumulants, spreadability and the {C}ampbell-{B}aker-{H}ausdorff
  series.
\newblock Pre\-print, 2017.
\newblock \href {https://arxiv.org/abs/1711.00219} {\path{arXiv:1711.00219}}.

\bibitem[HO07]{HoraObata07}
A.~Hora and N.~Obata.
\newblock {\em Quantum probability and spectral analysis of graphs}.
\newblock Theoretical and Mathematical Physics. Springer, Berlin, 2007.
\newblock With a foreword by Luigi Accardi.
\newblock \href {https://doi.org/10.1007/3-540-48863-4}
  {\path{doi:10.1007/3-540-48863-4}}.

\bibitem[HW16]{HuangWang16}
H.-W. Huang and J.-C. Wang.
\newblock Analytic aspects of the bi-free partial {$R$}-transform.
\newblock {\em J. Funct. Anal.}, 271(4):922--957, 2016.
\newblock \href {https://doi.org/10.1016/j.jfa.2016.04.026}
  {\path{doi:10.1016/j.jfa.2016.04.026}}.

\bibitem[JL20]{JekelLiu_tree_independence}
D.~Jekel and W.~Liu.
\newblock An operad of non-commutative independences defined by trees.
\newblock {\em Dissertationes Math.}, 553:100, 2020.
\newblock \href {https://doi.org/10.4064/dm797-6-2020}
  {\path{doi:10.4064/dm797-6-2020}}.

\bibitem[Lac15]{Lachs15}
S.~Lachs.
\newblock {\em A new family of universal products and aspects of a non-positive
  quantum probability theory}.
\newblock PhD thesis, EMAU Greifswald, 2015.
\newblock \url{http://ub-ed.ub.uni-greifswald.de/opus/volltexte/2015/2242/}.

\bibitem[Leh04]{Lehner04}
F.~Lehner.
\newblock Cumulants in noncommutative probability theory. {I}. {N}oncommutative
  exchangeability systems.
\newblock {\em Math. Z.}, 248(1):67--100, 2004.
\newblock \href {https://doi.org/10.1007/s00209-004-0653-0}
  {\path{doi:10.1007/s00209-004-0653-0}}.

\bibitem[Len10]{Lenczewski10}
R.~Lenczewski.
\newblock Matricially free random variables.
\newblock {\em J. Funct. Anal.}, 258(12):4075--4121, 2010.
\newblock \href {https://doi.org/10.1016/j.jfa.2010.03.010}
  {\path{doi:10.1016/j.jfa.2010.03.010}}.

\bibitem[Liu18]{Liu_preprint}
W.~Liu.
\newblock Free-free-boolean independence for triples of algebras.
\newblock Pre\-print, 2018.
\newblock \href {https://arxiv.org/abs/1801.03401} {\path{arXiv:1801.03401}}.

\bibitem[Liu19]{Liu19}
W.~Liu.
\newblock Free-{B}oolean independence for pairs of algebras.
\newblock {\em J. Funct. Anal.}, 277(4):994--1028, 2019.
\newblock \href {https://doi.org/10.1016/j.jfa.2019.05.005}
  {\path{doi:10.1016/j.jfa.2019.05.005}}.

\bibitem[Mac71]{MacLane}
S.~MacLane.
\newblock {\em Categories for the working mathematician}.
\newblock Graduate Texts in Mathematics, Vol. 5. Springer-Verlag, New
  York-Berlin, 1971.

\bibitem[MS17]{ManzelSchurmann17}
S.~Manzel and M.~Sch\"{u}rmann.
\newblock Non-commutative stochastic independence and cumulants.
\newblock {\em Infin. Dimens. Anal. Quantum Probab. Relat. Top.},
  20(2):1750010, 38, 2017.
\newblock \href {https://doi.org/10.1142/S0219025717500102}
  {\path{doi:10.1142/S0219025717500102}}.

\bibitem[MS01]{MS01}
M.~M. Meerschaert and H.-P. Scheffler.
\newblock {\em Limit distributions for sums of independent random vectors}.
\newblock Wiley Series in Probability and Statistics: Probability and
  Statistics. John Wiley \& Sons, Inc., New York, 2001.
\newblock Heavy tails in theory and practice.

\bibitem[MS17]{MingoSpeicher}
J.~A. Mingo and R.~Speicher.
\newblock {\em Free probability and random matrices}, volume~35 of {\em Fields
  Institute Monographs}.
\newblock Springer, New York; Fields Institute for Research in Mathematical
  Sciences, Toronto, ON, 2017.
\newblock \href {https://doi.org/10.1007/978-1-4939-6942-5}
  {\path{doi:10.1007/978-1-4939-6942-5}}.

\bibitem[M{\l}o04]{Mlotkowski04}
W.~M{\l}otkowski.
\newblock {$\Lambda$}-free probability.
\newblock {\em Infin. Dimens. Anal. Quantum Probab. Relat. Top.}, 7(1):27--41,
  2004.
\newblock \href {https://doi.org/10.1142/S0219025704001517}
  {\path{doi:10.1142/S0219025704001517}}.

\bibitem[Mur00]{Muraki00p}
N.~Muraki.
\newblock Monotonic convolution and monotonic {L}\'evy-{H}in\v{c}in formula.
\newblock Pre\-print, 2000.
\newblock Available from
  \url{https://www.math.sci.hokudai.ac.jp/~thasebe/Muraki2000.pdf}.

\bibitem[Mur01]{Muraki01}
N.~Muraki.
\newblock Monotonic independence, monotonic central limit theorem and monotonic
  law of small numbers.
\newblock {\em Infin. Dimens. Anal. Quantum Probab. Relat. Top.}, 4(1):39--58,
  2001.
\newblock \href {https://doi.org/10.1142/S0219025701000339}
  {\path{doi:10.1142/S0219025701000339}}.

\bibitem[Mur02]{Muraki02}
N.~Muraki.
\newblock The five independences as quasi-universal products.
\newblock {\em Infin. Dimens. Anal. Quantum Probab. Relat. Top.},
  5(1):113--134, 2002.
\newblock \href {https://doi.org/10.1142/S0219025702000742}
  {\path{doi:10.1142/S0219025702000742}}.

\bibitem[Mur03]{Muraki03}
N.~Muraki.
\newblock The five independences as natural products.
\newblock {\em Infin. Dimens. Anal. Quantum Probab. Relat. Top.},
  6(3):337--371, 2003.
\newblock \href {https://doi.org/10.1142/S0219025703001365}
  {\path{doi:10.1142/S0219025703001365}}.

\bibitem[Mur13a]{Muraki13t}
N.~Muraki.
\newblock On a q-deformation of free independence.
\newblock Talk at Workshop on Combinatorial and Random Matrix Aspects of
  Noncommutative Distributions and Free Probability, Fields Institute, Toronto,
  2013.
\newblock Available from
  \url{http://www.fields.utoronto.ca/talks/q-deformation-free-independence}.

\bibitem[Mur13b]{Muraki13}
N.~Muraki.
\newblock A simple proof of the classification theorem for positive natural
  products.
\newblock {\em Probab. Math. Statist.}, 33(2):315--326, 2013.
\newblock Available from
  \url{https://www.math.uni.wroc.pl/~pms/files/33.2/Article/33.2.12.pdf}.

\bibitem[NS06]{NicaSpeicher}
A.~Nica and R.~Speicher.
\newblock {\em Lectures on the combinatorics of free probability}, volume 335
  of {\em London Mathematical Society Lecture Note Series}.
\newblock Cambridge University Press, Cambridge, 2006.
\newblock \href {https://doi.org/10.1017/CBO9780511735127}
  {\path{doi:10.1017/CBO9780511735127}}.

\bibitem[Sat13]{Sato13}
K.~Sato.
\newblock {\em L\'{e}vy processes and infinitely divisible distributions},
  volume~68 of {\em Cambridge Studies in Advanced Mathematics}.
\newblock Cambridge University Press, Cambridge, 2013.
\newblock Translated from the 1990 Japanese original, Revised edition of the
  1999 English translation.

\bibitem[Sko16]{Skoufranis_AIHP}
P.~Skoufranis.
\newblock Independences and partial {$R$}-transforms in bi-free probability.
\newblock {\em Ann. Inst. Henri Poincar\'{e} Probab. Stat.}, 52(3):1437--1473,
  2016.
\newblock \href {https://doi.org/10.1214/15-AIHP691}
  {\path{doi:10.1214/15-AIHP691}}.

\bibitem[Spe]{SpeicherCours20182019}
R.~Speicher.
\newblock Free probability theory. {L}ecture notes. {W}inter 2018/19.
\newblock Available from
  \url{https://rolandspeicher.files.wordpress.com/2019/08/free-probability.pdf}.

\bibitem[Spe97]{Speicher97}
R.~Speicher.
\newblock On universal products.
\newblock In {\em Free probability theory ({W}aterloo, {ON}, 1995)}, volume~12
  of {\em Fields Inst. Commun.}, pages 257--266. Amer. Math. Soc., Providence,
  RI, 1997.
\newblock \href {https://doi.org/10.1090/fic/012} {\path{doi:10.1090/fic/012}}.

\bibitem[SW97]{SpeicherWoroudi97}
R.~Speicher and R.~Woroudi.
\newblock Boolean convolution.
\newblock In {\em Free probability theory ({W}aterloo, {ON}, 1995)}, volume~12
  of {\em Fields Inst. Commun.}, pages 267--279. Amer. Math. Soc., Providence,
  RI, 1997.
\newblock \href {https://doi.org/10.1090/fic/012} {\path{doi:10.1090/fic/012}}.

\bibitem[Var21]{Varso}
P.~Var\v{s}o.
\newblock {\em Studies on Positive and Symmetric Two-Faced Universal Products}.
\newblock PhD thesis, University of Greifswald, 2021.
\newblock Available from
  \url{https://nbn-resolving.org/urn:nbn:de:gbv:9-opus-75524}.

\bibitem[Voi85]{Voiculescu85}
D.~Voiculescu.
\newblock Symmetries of some reduced free product {$C^\ast$}-algebras.
\newblock In {\em Operator algebras and their connections with topology and
  ergodic theory ({B}u\c{s}teni, 1983)}, volume 1132 of {\em Lecture Notes in
  Math.}, pages 556--588. Springer, Berlin, 1985.
\newblock \href {https://doi.org/10.1007/BFb0074909}
  {\path{doi:10.1007/BFb0074909}}.

\bibitem[Voi14]{Voiculescu14}
D.-V. Voiculescu.
\newblock Free probability for pairs of faces {I}.
\newblock {\em Comm. Math. Phys.}, 332(3):955--980, 2014.
\newblock \href {https://doi.org/10.1007/s00220-014-2060-7}
  {\path{doi:10.1007/s00220-014-2060-7}}.

\bibitem[VDN92]{VoiculescuDykemaNica}
D.~V. Voiculescu, K.~J. Dykema, and A.~Nica.
\newblock {\em Free random variables}, volume~1 of {\em CRM Monograph Series}.
\newblock American Mathematical Society, Providence, RI, 1992.
\newblock A noncommutative probability approach to free products with
  applications to random matrices, operator algebras and harmonic analysis on
  free groups.
\newblock \href {https://doi.org/10.1090/crmm/001}
  {\path{doi:10.1090/crmm/001}}.

\bibitem[vW73]{vWaldenfels73}
W.~von Waldenfels.
\newblock An approach to the theory of pressure broadening of spectral lines.
\newblock In {\em Probability and information theory, {II}}, pages 19--69.
  Lecture Notes in Math., Vol. 296. Springer, Berlin, 1973.
\newblock \href {https://doi.org/10.1007/BFb0059817}
  {\path{doi:10.1007/BFb0059817}}.

\bibitem[vW75]{vWaldenfels75}
W.~von Waldenfels.
\newblock Interval partitions and pair interactions.
\newblock In {\em S\'{e}minaire de {P}robabilit\'{e}s, {IX} ({S}econde
  {P}artie, {U}niv. {S}trasbourg, {S}trasbourg, ann\'{e}es universitaires
  1973/1974 et 1974/1975)}, pages 565--588. Lecture Notes in Math., Vol. 465.
  Springer, Berlin, 1975.
\newblock Available from \url{http://www.numdam.org/item/SPS_1975__9__565_0/}.

\bibitem[Wys10]{Wysoczanski10}
J.~Wysocza\'{n}ski.
\newblock bm-independence and bm-central limit theorems associated with
  symmetric cones.
\newblock {\em Infin. Dimens. Anal. Quantum Probab. Relat. Top.},
  13(3):461--488, 2010.
\newblock \href {https://doi.org/10.1142/S0219025710004115}
  {\path{doi:10.1142/S0219025710004115}}.

\end{thebibliography}

\setlength{\parindent}{0pt}

\end{document}